\documentclass[a4paper,12pt]{article}
\usepackage[dvips]{epsfig}
\usepackage{amsmath,amssymb, amsbsy, amstext,amscd,amsfonts}
\usepackage{color,enumerate,euscript,graphicx,hyperref}
\usepackage{srcltx,tikz}
\input amssym.def
\input amssym.tex

\newtheorem{lem}{Lemma}[section]%
\newtheorem{theorem}[lem]{Theorem}%
\newtheorem{prop}[lem]{Proposition}%
\newtheorem{rem}[lem]{Remark}%

\def\a{\alpha}    
  \def\t{\tau}

 \def\O{\Omega}

 \def\og{\overline G} \def\oh{\overline H}  \def\oc{\overline C}

 \def\ox{\overline X}  \def\o1{\overline 1}

\def\olg{\overline g} \def\ola{\overline a} 
  \def\olc{\overline c} \def\olb{\overline b}
  \def\olx{\overline x} 

\def\o{\overline}   \def\olb{\overline b}
\def\di{\bigm|} \def\lg{\langle} \def\rg{\rangle}

\def\Aut{\hbox{\rm Aut\,}} \def\Inn{\hbox{\rm Inn}} \def\Syl{\hbox{\rm Syl}}
  
  \def\mod{\hbox{\rm mod }}

 \def\AGL{\hbox{\rm AGL}} \def\GL{\hbox{\rm GL}}  \def\P\GL{\hbox{\rm P\GL}}
\def\GF{\hbox{\rm GF}}  \def\FF{{\hbox{\sf F\kern-.43emF}}}
 \def\char{\hbox{\rm char}} 
 
\def\Inn{\hbox{\rm Inn}}

\def\o{\hbox{\rm o}}

\def\char{ \, {\rm char}\,}

\def\ZZ{\mathbb{Z}}  

\def\nd{\mathrel{\bigm|\kern-.7em/}} 
 \def\f{\noindent}
\def\qed{\hfill $\Box$} \def\demo{\f {\bf Proof}\hskip10pt}

\begin{document}
\begin{center}
{\bf\large  The Product of a Generalized Quaternion Group\\And a Cyclic Group$^*$ }
\end{center}


\begin{center}
Shaofei Du{\small \footnotemark},  Hao Yu and Wenjuan Luo\\
\medskip
 {\small
Capital Normal University,\\ School of Mathematical Sciences,\\
Beijing 100048, People's Republic of China
}
\end{center}

\footnotetext{Corresponding author: dushf@mail.cnu.edu.cn. }

\renewcommand{\thefootnote}{\empty}
\footnotetext{*This work is supported in part  by the National Natural Science Foundation of China (12071312).}
\begin{abstract}
  Let $X(Q)=QC$ be a group, where $Q$ is a generalized quaternion group and $C$ is a cyclic group such that
  $Q\cap C=1$. In this paper, $X(Q)$ will be characterized and moreover, a complete classification for that will be given, provided $C$ is core-free.
  For the reason of self-constraint,  in this paper a classification of the group $X(D)=DC$ is also given,    where  $D$ is a dihedral group and $C$ is a cyclic group such that $D\cap C=1$ and $C$ is core-free.
  Remind that  the group  $X(D)$  was recently classified  in \cite{HKK2022},  based on a number of papers on skew-morphisms of dihedral groups.
  In  this paper,  a  different approach from that in \cite{HKK2022} will be used.
  \end{abstract}

\section{Introduction}
A group $G$ is said to be properly {\it factorizable} if $G=AB$ for two proper subgroups $A$ and $B$ of $G$, while the  expression $G=AB$ is called a {\it factorization} of $G$.
Furthermore, if $A\cap B=1$,  then  we say that $G$ has an {\it exact factorization}.

\smallskip

Factorizations of groups naturally arise from the well-known  Frattini's argument, including  its version in permutation groups.
One of the most famous results about factorized groups might be one of theorems of It\^o, saying that  any group is metabelian whenever it is the product of two abelian subgroups (see \cite{Ito1}).
Later, Wielandt and  Kegel showed that the product of two nilpotent subgroups must be soluble
(see \cite{Wie1958} and \cite{Ke1961}).
Douglas showed that  the product of two cyclic groups must be super-solvable (see \cite{D1961}).
The factorizations of the finite almost simple groups were determined in \cite{LPS1990}  and
the  factorizations of almost simple groups with a solvable factor  were determined in  \cite{LX2022}.
There are many other papers related to factorizations, for instance, finite products of soluble groups, factorizations with one nilpotent factor and so on.
Here we are not able to list all references and the readers may refer to a survey paper \cite{AK2009}.
\smallskip

In this paper, we shall focus on the product group $X=GC$, for a finite group $G$ and a cyclic group $C$ such that $G\cap C=1$.
Suppose $C$ is core-free. Then $X$ is  also called a {\it skew product group} of $G$.
Recall that the skew morphism of a group $G$ and a skew product group $X$ of $G$
were introduced by Jajcay and \v{S}ir\'a\v{n} in \cite{JS2002},
which is related to the studies of regular Cayley maps of $G$.
For the reason of the length of the paper, we are not able to explain them in detail.
Recently, there have been a lot of results on skew product groups $X$ of some particular groups $G$.
(1) Cyclic groups: So far there exists no classification of such product groups.
For  partial results, see~\cite{CJT2016,CT,DH,KN1,KN2,Kwo}.
(2) Elementary abelian $p$-groups: a global structure  was characterized in \cite{DYL}.
(3) Finite nonabelian simple group or finite nonabelian characteristically simple groups:
they were classified  in  \cite{BCV2019} and \cite{CDL}, respectively.
(4) Dihedral groups: Based on big efforts of several authors working on regular Cayley maps (see \cite{CJT2016,HKK2022,KKF2006, KMM2013,KK2017,KK2016,RSJTW2005,KK2021,WF2005,WHY2019,Zhang2015,
Zhang20152,ZD2016}),
the final classification of  skew product groups of dihedral groups was given in \cite{HKK2022}.
(5) Generalized quaternion groups: for  partial results, see \cite{HR2022} and \cite{KO2008}.

\smallskip
By $Q$ and $D$, we denote a generalized quaternion group and a dihedral group, respectively.
Let $X(G)=GC$ be a group, where $G\in\{Q,\,D\}$ and $C$ is a cyclic group such that $G\cap C=1$.
In this paper, we shall give a characterization for $X(Q)$ and a complete classification of $X(Q)$ provided $C$ is core-free.
In the most of the above papers dealing with skew product groups of dihedral groups,
the authors adopt some computational technics on skew-morphisms.
Alternatively, in this paper, we shall realize our goals by using classical group theoretical tools
and methods (solvable groups, $p$-groups, permutation groups, group extension theory and so on).
Address that we shall pay attention to the global structures of the group $X(G)$.
Since $X(Q)$ is closely related to $X(D)$ and we shall adopt a completely different approach,
the group $X(D)$ will be considered too, for the reason of self-constraint in this paper.

\smallskip

Throughout this paper, set $C=\lg c\rg $ and
\begin{eqnarray}\label{main0}
\begin{array}{ll}
&Q=\lg a, b\di a^{2n}=1,b^2=a^n,a^b=a^{-1}\rg\cong Q_{4n},\,n\ge 2,\\
&D=\lg a, b\di a^{n}=b^2=1, a^b=a^{-1}\rg \cong D_{2n}, \, n\ge 2.
\end{array}
\end{eqnarray}
Let $G\in\{ Q, D\}$ and $X=X(G)=GC=\lg a, b\rg \lg c\rg $.
Then $\lg a\rg \lg c\rg $ is unnecessarily a subgroup of $X$.
Clearly, $X$ contains a subgroup $M$ of the biggest order such that $\lg c\rg\le M\subseteqq\lg a\rg\lg c\rg$.
This subgroup $M$ will play an important role in this paper.
From now on by $S_X$ we denote the core $\cap_{x\in X}S^x$ of $X$ in a subgroup $S$ of $X$.

\smallskip
There are four main theorems in this manuscript.
In Theorem~\ref{main1}, the global structure of our group $X\in\{X(Q), X(D)\}$ is characterized.

\begin{theorem} \label{main1}
Let $G\in\{ Q,\,D\}$ and $X=G\lg c\rg\in \{X(Q),\,X(D)\}$, where $\o(c)=m\ge 2$ and $G\cap \lg c\rg=1$.
Let $M$ be the subgroup of the biggest order in $X$ such that $\lg c\rg \le M\subseteqq \lg a\rg \lg c\rg$.  Then one of items in Tables \ref{table1} holds.
 \vskip 3mm
\begin{table}
 \center   \caption {The forms of $M$, $M_X$ and $X/M_X$}\label{table1}
 \begin{tabular}{cccc}
  \hline
  Case &$M$ & $M_X$  & $X/M_X$\\
  \hline
   1 & $\lg a\rg\lg c\rg$   & $\lg a\rg\lg c\rg$     &   $\ZZ_2$ \\
   2 &$\lg a^2\rg \lg c\rg$ & $\lg a^2\rg\lg c^2\rg$ &   $D_8$   \\
   3 &$\lg a^2\rg \lg c\rg$ & $\lg a^2\rg\lg c^3\rg$ &   $A_4$   \\
   4 &$\lg a^4\rg \lg c\rg$ & $\lg a^4\rg\lg c^3\rg$ &   $S_4$   \\
   5 &$\lg a^3\rg \lg c\rg$ & $\lg a^3\rg\lg c^4\rg$ &   $S_4$   \\
   \hline
  \end{tabular}
  \end{table}
\end{theorem}

Clearly, $M$ is a product of two cyclic subgroups,
 which  has not been determined so far, as mentioned before,
However, further  properties of our group $X$ is given in Theorem~\ref{main2}.
More powerful properties will be obtained during the proof of Theorems~\ref{main3} and \ref{main4},
see Remarks~\ref{rem1} and~\ref{rem2}.

\begin{theorem} \label{main2}
Let $G\in \{Q, D\}$ and $X\in \{X(Q), X(D)\}$, and $M$  defined as above.
Then we have  $\lg a^2, c\rg \le C_X(\lg c\rg_X)$ and $|X: C_X(\lg c\rg_X)|\le 4$.
Moreover, if $\lg c\rg _X=1$, then $M_X\cap \lg a^2\rg \lhd M_X$.
In particular, if $\lg c\rg _X=1$ and $M=\lg a\rg \lg c\rg$, then $\lg a^2\rg \lhd X$.
\end{theorem}

In Theorem~\ref{main3}, a classification of $X(Q)$ is given, provided that  $C$ is core-free.

\begin{theorem} \label{main3}
Let $X=X(Q)$.
Set $R:=\{a^{2n}=c^m=1,\,b^2=a^n,\,a^b=a^{-1}\}$.
Suppose $\lg c\rg _X=1$.
Then   $X$ is isomorphic to one of the following groups:
\begin{enumerate}
  \item[\rm(1)] $X=\lg a,b,c| R, (a^2)^c=a^{2r},  c^a=a^{2s}c^t,c^b=a^uc^v\rg,$ where
  $$\begin{array}{ll}
  &r^{t-1}-1\equiv r^{v-1}-1\equiv0(\mod n),\, t^2\equiv 1(\mod m),\\
  &2s\sum_{l=1}^tr^{l}+2sr\equiv 2sr+2s\sum_{l=1}^{v}r^l-u\sum_{l=1}^tr^l+ur \equiv 2(1-r)(\mod 2n),\\
  &2s\sum_{l=1}^{w}r^l\equiv u\sum_{l=1}^{w}(1-s(\sum_{l=1}^tr^l+r))^l\equiv 0(\mod 2n)\Leftrightarrow
  w\equiv 0(\mod m).
  \end{array}$$
  Moreover, if $2\di n$, then $u(\sum_{l=0}^{v-1}r^l-1)\equiv 0(\mod 2n)$ and $v^2\equiv 1(\mod m)$;
  if $2\nmid n$, then $u\sum_{l=1}^vr^l-ur\equiv 2sr+(n-1)(1-r)(\mod 2n)$ and $v^2\equiv t(\mod m)$;
  if $t\ne 1$, then $u\equiv0(\mod2)$.
  \item[\rm(2)] $X=\lg a,b,c|R,(a^2)^{c^2}=a^{2r},(c^2)^a=a^{2s}c^{2t},(c^2)^b=a^{2u}c^{2}, a^c=bc^{2w}\rg,$ where either $w=0$ and $r=s=t=u=1$; or
  $$\begin{array}{ll}
  &w\neq0,\,s=u^2\sum_{l=0}^{w-1}r^l+\frac{un}2,\,t=2wu+1,\\
  &r^{2w}-1\equiv (u\sum_{l=1}^{w}r^l+\frac n2)^2-r\equiv0(\mod n),\\
  &s\sum_{l=1}^tr^{l}+sr\equiv 2sr-u\sum_{l=1}^tr^l+ur\equiv 1-r(\mod n),\\
  &2w(1+uw)\equiv nw\equiv 2w(r-1)\equiv0(\mod\frac m2),\\
  &2^{\frac{1+(-1)^u}2}\sum_{l=1}^{i}r^l\equiv 0(\mod n)\Leftrightarrow i\equiv0(\mod\frac m2).
   \end{array}$$
  \item[\rm(3)] $X=\lg a,b,c|R, (a^2)^{c}=a^{2r},(c^3)^a=a^{2s}c^3,(c^3)^b=a^{2u}c^{3}, a^c=bc^{\frac {im}2},b^c=a^xb\rg,$
  where $n\equiv 2(\mod 4)$ and either $i=0$ and $r=x=u=1$; or $i=1$, $6\di m$,
  $r^{\frac m2}\equiv-1(\mod n)$ with $\o(r)=m$, $s\equiv \frac{r^{-3}-1}2(\mod\frac n2)$, $u\equiv\frac{r^3-1}{2r^2}(\mod\frac n2)$ and $x\equiv -r+r^2+\frac n2 (\mod n)$.
  \item[\rm(4)] $X=\lg a,b,c|R, (a_1^2)^{c}=a_1^{2r}, c_1^{a_1}=a_1^{2s}c_1,c_1^b=a_1^{2u}c_1,
  a_1^c=bc^{\frac {im}2},b^c=a_1^xb,c^a=a_1^{1+2z}c^{1+\frac{jm}3}\rg,$ where either
  $i=0,\,r=1,\,x=3$ and $s=u=z=d=0$; or
  $i=1$, $n\equiv 4(\mod 8), m\equiv 0(\mod 6),$
  $r^{\frac m2}\equiv-1(\mod \frac n2),\,\o(r)=m$,
  $s\equiv \frac{r^{-3}-1}2(\mod\frac n4),\,u\equiv\frac{r^3-1}{2r^2}(\mod\frac n4),\,
  x\equiv -r+r^2+\frac n4(\mod \frac n2),\,1+2z\equiv \frac{1-r}{2r}(\mod\frac n2),\,
  j\in\{ 1,2\}$.
  \item[\rm(5)] $X=\lg a,b,c|R,a^{c^4}=a^r, b^{c^4}=a^{1-r}b, (a^3)^{c^{\frac m4}}=a^{-3}, a^{c^{\frac m4}}=bc^{\frac{3m}4}\rg$, where $m\equiv 4(\mod8)$ and $r$ is of order $\frac m4$ in $\ZZ_{2n}^*$.
\end{enumerate}
Moreover, in the families of groups (1)-(5), for any given parameters satisfying the equations,
there exists  $X=X(Q)$.
\end{theorem}

In Theorem~\ref{main4}, a classification of $X(D)$ is given, provided that $C$ is core-free.
Remind that our presentations for $X(D)$ are different form that in \cite{HKK2022}
but they are essentially isomorphic.

\begin{theorem} \label{main4}
Let $X=X(D)$.
Set $R:=\{a^{n}=b^2=c^m=1,\,a^b=a^{-1}\}$.
Suppose $\lg c\rg _X=1.$
Then $X$ is isomorphic to one of the following groups:
\begin{enumerate}
  \item[\rm(1)] $X=\lg a,b,c| R, (a^2)^c=a^{2r},  c^a=a^{2s}c^t,c^b=a^uc^v\rg,$ where
  $$\begin{array}{ll}
  &2(r^{t-1}-1)\equiv 2(r^{v-1}-1)\equiv u(\sum_{l=0}^{v-1}r^l-1)\equiv0(\mod n),\,
  t^2\equiv v^2\equiv1(\mod m)\\
  &2s\sum_{l=1}^tr^{l}+2sr\equiv 2sr+2s\sum_{l=1}^{v}r^l-u\sum_{l=1}^tr^l+ur\equiv 2(1-r)(\mod n),\\
  &{\rm if}\, t\ne 1,\, {\rm then}\, u\equiv0(\mod2),\\
  &2s\sum_{l=1}^{w}r^l\equiv u\sum_{l=1}^{w}(1-s(\sum_{l=1}^tr^{l}+r))^l\equiv0(\mod n)\Leftrightarrow w\equiv0(\mod m).
  \end{array}$$
  \item[\rm(2)] $X=\lg a,b,c|R,(a^2)^{c^2}=a^{2r}, (c^2)^b=a^{2s}c^{2}, (c^2)^a=a^{2u}c^{2v}, a^c=bc^{2w}\rg,$ where either $w=s=u=0$ and $r=t=1$; or
  $$\begin{array}{ll}
  &w\neq0,\,s=u^2\sum_{l=0}^{w-1}r^l,\,t=1+2wu,\\
  &nw\equiv2w(r-1)\equiv2w(1+uw)\equiv0(\mod\frac m2),\\
  &r^{2w}-1\equiv(u\sum_{l=1}^{w}r^l)^2-r\equiv (r^w+1)(1+s\sum_{l=0}^{w-1}r^l) \equiv0(\mod\frac n2),\\
  &\sum_{l=1}^{i}r^l\equiv0(\mod\frac n2)\Leftrightarrow i\equiv0(\mod\frac m2).
  \end{array}$$
  \item[\rm(3)] $X=\lg a,b,c|R, a^{c^3}=a^r,(c^3)^b=a^{2u}c^{3}, a^c=bc^{\frac {im}2},b^c=a^xb\rg,$
  where $n\equiv 2(\mod 4)$ and either $i=u=0$ and $r=x=1$; or $i=1$, $6\di m$,
  $l^{\frac m2}\equiv-1(\mod\frac n2)$ with $\o(l)=m$, $r=l^3$, $u=\frac{l^3-1}{2l^2}$ and $x\equiv -l+l^2+\frac n2(\mod n)$.
  \item[\rm(4)] $X=\lg a,b,c|R, (a^2)^{c^3}=a^{2r}, (c^3)^b=a^{\frac{2(l^3-1)}{l^2}}c^3,
  (a^2)^c=bc^{\frac {im}2}, b^c=a^{2(-l+l^2+\frac n4)}b,c^a=a^{2+4z}c^{2+3d}\rg,$ where either $i=z=d=0$ and $l=1$; or $i=1$, $n\equiv 4(\mod 8),\, m\equiv 0(\mod 6)$,
  $l^{\frac m2}\equiv-1(\mod\frac n4)$ with $\o(l)=m$, $r=l^3$, $z=\frac {1-3l}{4l}$,
  $1+3d\equiv0(\mod\frac m3)$ and $\sum_{i=1}^{j}r^i\equiv0(\mod\frac n2)\Leftrightarrow
  j\equiv0(\mod\frac m3)$.
  \item[\rm(5)] $X=\lg a,b,c|R,a^{c^4}=a^r, b^{c^4}=a^{1-r}b, (a^3)^{c^{\frac m4}}=a^{-3}, a^{c^{\frac m4}}=bc^{\frac{3m}4}\rg,$ where  $m\equiv 4(\mod8)$ and $r$ is of order $\frac m4$ in $\ZZ_{n}^*$.
\end{enumerate}
Moreover, in the families of groups (1)-(5), for any given parameters satisfying the equations,
there exists  $X=X(D)$.
\end{theorem}

\begin{rem}\label{rem1}
From Theorems~\ref{main3} and \ref{main4},
one may observe that $\lg a^4\rg\lhd X$ for groups in (3) and (4); and $\lg a^3\rg\lhd X$ for groups in (5).
\end{rem}

\begin{rem} \label{rem2}
Checking Theorem~\ref{main3},
we know that $\lg a^n\rg\lhd X(Q)$ for all cases (2) and (3) and some cases in (1).
Moreover, corresponding to $D\cong Q/\lg a^n\rg$, we have $X(D)=X(Q)/\lg a^n, c_1\rg$,
where $\lg a^n\rg \lhd X$ and  $\lg a^n, c_1\rg=\lg a^n, c\rg_{X(Q)}$.
\end{rem}

\begin{rem} \label{rem3}
The group $X$  where $\lg c\rg_X\ne 1$  is a cyclic extension of $X/\lg c\rg_X$ which is given in Theorem~\ref{main3} and   Theorem~\ref{main4}, respectively.
Furthermore,  by Theorem~\ref{main2}, the subgroup $C_X(\lg c\rg_X)$ of $X$ is of index at most 4.
So one may classify $X$ by using Theorems~\ref{main3} and \ref{main4}.
This needs complicate computations and we cannot do it in this paper.
\end{rem}

\begin{rem} \label{rem4}
One may determine regular Cayley maps of dihedral groups by Theorem~\ref{main4}
(which were  done in \cite{KK2021} via skew-morphism computations)
and of  generalized quaternion groups by Theorem~\ref{main3}.
\end{rem}

After this introductory section, some  preliminary results will be given in Section 2, Theorems~\ref{main1}--\ref{main4} will be proved in Sections 3-6, respectively.

\section{Preliminaries}
In this section, the notation and elementary facts used in this paper are  collected.

\subsection{Notation}
In this paper, all the groups are supposed to be finite.
We set up the notation below, where $G$ and $H$ are groups, $M$ is a subgroup of $G$,
$n$ is a positive integer and $p$ is a prime number.
\begin{enumerate}
  \setlength{\itemsep}{0ex}
  \setlength{\itemindent}{-0.5em}
  \item[] $|G|$ and $\o(g)$: the order of $G$ and an element $g$ in $G$, resp.;
  \item[] $H\leq G$ and $H<G$: $H$ is a subgroup of $G$ and $H$ is a proper subgroup of $G$, resp.;
  \item[] $[G:H]$: the set of cosets of   $G$ relative to a subgroup $H$;
  \item[] $H\lhd G$ and $H\char~G$: $H$ is a normal and characteristic subgroup of $G$, resp.;
  \item[] $G'$ and $Z(G)$: the derived subgroup and the center of $G$ resp.;
  \item[] $M_G$: the core of $M$ in $G$ which is the maximal normal subgroup of $G$ contained in $M$;
  \item[] $G\rtimes H$: a semidirect product of $G$ by $H$, in which $G$ is  normal;
  \item[] $G.H$:  an extension of $G$ by $H$, where $G$ is normal;
  \item[] $C_M(G)$: centralizer of $M$ in $G$;
  \item[] $N_M(G)$: normalizer of $M$ in $G$;
  \item[] $\Syl_p(G)$: the set of all Sylow $p$-subgroups of $G$;
  \item[] $[a,b]:=a^{-1}b^{-1}ab$, the commutator of $a$ and $b$ in $G$;
  \item[] $\O_1(G)$: the subgroup $\lg g\in G\di g^p=1\rg$ of $G$ where $G$ is a $p$-group;
  \item[] $\mho_n(G)$: the subgroup $\lg g^{p^n}\di g\in G$ of $G$ where $G$ is a $p$-group;
  \item[] $S_n$: the symmetric group of degree $n$ (naturally acting on $\{1, 2, \cdots, n\}$);
  \item[] $A_n$: the alternating group of degree $n$ (naturally acting on $\{1, 2, \cdots, n\}$);
  \item[] $\GF(q)$: finite field of $q$ elements;
  \item[] $\AGL(n,p)$: the affine group on $\GF^n(q)$.
\end{enumerate}
\subsection{Elementary facts}
\begin{prop}\cite[Theorem 1]{M1974}\label{solvable}
The finite group $G=AB$ is solvable,
where both $A$ and $B$ are subgroups with cyclic subgroups of index no more than 2.
\end{prop}

Recall that a group $H$ is said a {\it Burnside group}
if every permutation group containing a regular subgroup isomorphic to $H$ is
either 2-transitive or imprimitive.
The following results are well-known.

\begin{prop} \cite[Theorem 25.3 and Theorem 25.6]{W1964} \label{Burnside}
Every cyclic group of a composite order is a Brunside group.
Every dihedral group is a Burnside group.
\end{prop}

\begin{prop}\cite[Lemma 4.1]{DMM2008}\label{AGL}
Let $n\ge 2$ be an integer and $p$ a prime.  Then $\AGL(n, p)$ contains an element of order $p^n$ if and only if  $(n, p)=(2, 2)$
  and $\AGL(2, 2)\cong S_4$.
\end{prop}

Recall that  our group $X(D)=DC$, where $D$ is a dihedral group of order $2n$  and $C$ is a cyclic group of order $m$ such that $D\cap C=1$, where $n, m\ge 2$.
Then  we have the following results.

\begin{lem}\label{cd}
 Suppose that $X(D)$ is a solvable and has a faithful 2-transitive permutation  representation relative to a subgroup $M$, which is of index  a composite order. 
 Then    $X(D)\le \AGL(k,p)$. Moreover, 
  \begin{enumerate}
  \item[\rm(i)]   if $X(D)$ contains an element of order $p^k$,  then $X(D)=S_4$;
  \item[\rm(ii)]    if the hypotheses holds for $M=C$ where $C$ is core-free,  then   $X(D)=A_4$.
  \end{enumerate}
\end{lem}
\demo Set $\O=[X(D):M]$.  Let $N$ be a minimal normal subgroup of $X(D)$.
Since $X(D)$ is solvable, $N\cong \ZZ_p^k$ for some prime $p$ and integer $k$.
Since $X(D)$ is 2-transitive, it is  primitive, which implies that $N$ is transitive on $\O$ and so is regular on $\O$.
Therefore, $X(D)=N\rtimes X(D)_\a \le \AGL(k,p)$, for some $\a\in \O$.
Since $X(D)$ is 2-transitive and $|\O|=p^k$, we know $|X(D)_\a|\ge p^k-1$ for any $\a\in \O$.

\vskip 3mm
(i) Suppose that $X(D)$ contains an element of order $p^k$. By Proportion~\ref{AGL}, we  get $(k,p)=(2,2)$ so that  $X(D)=S_4$, reminding $|\O|$ is not a prime.

\vskip 3mm (ii)  Let $M=C$ where $C$ is core-free.   Set $C=\lg c\rg $ and $o(c)=m$.
Then  $X(D)=N\rtimes \lg c\rg $, where $\lg c\rg $ is a Singer subgroup of $\GL(k,p)$. Then  both $D$ and $N$ are regular subgroups, that is  $|D|=2n=|\O|=p^k$, which implies $p=2$. Now, we have
   $|X(D)|=2^{k}(2^k-1)=2n\cdot m=p^k\cdot m$ and so $m=2^k-1$. Since both $N$ and $D$ are Sylow 2-subgroups of $X(D)$ and $N\lhd X(D)$, we get $D=N$,  and so $D\cong \ZZ_2^2$.
 Therefore,  $p=2$ and $X(D)=A_4$. \qed

\begin{prop}\label{complement}\cite[Satz 1]{G1952}
Let $N\leq M\leq G$ such that $(|N|,|G:M|)=1$ and $N$ be  an abelian normal subgroup of $G$.
If $N$ has a complement in $M$, then $N$ also has a complement in $G$.
\end{prop}

The Schur multiplier $M(G)$ of a group $G$ is defined as the second integral homology group
$H^2(G; Z)$, where $Z$ is a trivial $G$-module. It plays an important role in the central expansion of groups.
The following result is well-known.

\begin{prop} \cite[ (2.21) of page 301]{S1982}\label{Schur}
The Schur multiplier $M(S_n)$ of $S_n$ is a cyclic group of order 2 if $n\ge4$ and
of order 1 for $n\leq 3$.
\end{prop}

\begin{prop}\cite[Theorem 4.5]{H1967}\label{NC}
Let $H$ be the subgroup of $G$. Then $N_G(H)/C_G(H)$ is isomorphic to a subgroup of $\Aut(H)$.
\end{prop}

\begin{prop}\cite[Theorem]{Luc} \label{cyclic}
If $G$ is a transitive permutation group of degree $n$  with a cyclic point-stabilizer,
then $|G|\le n(n-1)$.
\end{prop}

\begin{prop} \cite[Satz 1 and Satz 2]{Ito1} \label{mateabel}
Let $G=AB$ be a group, where both $A$ and $B$ are abelian subgroups of $G$. Then
\begin{enumerate}
  \item[\rm(1)]  $G$ is meta-abelian, that is, $G'$ is abelian;
  \item[\rm(2)]  if $G\ne 1$, then  $A$ or $B$ contains   a normal subgroup $N\ne 1$ of $G$.
\end{enumerate}
\end{prop}

\begin{prop}\cite[Theorem 11.5]{H1967}\label{ab}\label{matecyclic}
Let $G=\lg a\rg\lg b\rg$ be a group.
If $|\lg a\rg|\leq |\lg b\rg|$, then $\lg b\rg_G\ne1$.
If both $\lg a\rg$ and $\lg b\rg$ are $p$-groups where $p$ is an odd prime, then $G$ is matecyclic.
\end{prop}

\begin{prop}\cite[Corollary 1.3.3]{W1999}\label{sylowp}
Let $G=AB$ be a group, where both $A$ and $B$ are subgroups of $G$.
And let $A_p$ and $B_p$ be Sylow $p$-subgroups of $A$ and $B$ separately, for some prime $p$.
Then $A_pB_p$ is the Sylow $p$-subgroup of $G$.
\end{prop}

\begin{prop}\cite[Theorem 12.5.1]{H1959}\label{pcyclic}
Let $p$ is an odd prime.
Then every finite $p$-group $G$ containing  a cyclic maximal subgroup is isomorphic to
(1)  $\ZZ_{p^n}$;\,   (2) $\lg a,b\di a^{p^{n-1}}=b^p=1,\,[a,b]=1\rg,\,n\geq2$; or
(3) $\lg a,b\di a^{p^{n-1}}=b^p=1,\,[a,b]=a^{p^{n-2}}\rg,\,n\geq3$.
\end{prop}

\section{Proof of Theorem~\ref{main1}}
\vskip 3mm
To prove Theorem~\ref{main1}, let $G=D$ or $Q$, defined in Eq(\ref{main0}).
Let $X=G\lg c\rg \in \{ X(Q), X(D)\}$.
Let $M$ be the subgroup of the biggest order in $X$ such that $\lg c\rg \le M\subseteqq \lg a\rg \lg c\rg$,
and set $M_X=\cap_{x\in X}M^x$.
By Proposition~\ref{solvable}, $X$ is solvable.
Theorem~\ref{main1} will be proved by  dealing with $G=D$ and $Q$, separately in Lemmas~\ref{MD} and \ref{MQ}.
\begin{lem} \label{MD}
 Theorem~\ref{main1} holds, provided  $G=D$ and $X=X(D)$.
\end{lem}
\demo
Let $G=D$ so that $X=X(D)$.
Remind that   $m, n\ge 2$, $|X|$ is even and more than 7.
The lemma is proved by the induction on $|X|$.
All the cases when  $|X|\le 24$ are listed below, which implies that the conclusion holds:
\vskip 3mm
 $M=\lg a \rg \lg c\rg$: $a^n=c^m=1$, where $4\leq nm\leq12$;

 $M=\lg a^2\rg\lg c\rg$: $a^2=c^m=1$, where $m\in\{2,4,6\}$ $M_X=\lg a^2\rg\lg c^2\rg$, $X/M_X=D_8$;

 $M=\lg a^2\rg\lg c\rg$: $a^2=c^m=1$, where $m\in\{3,6\}$, $M_X=\lg a^2\rg\lg c^3\rg$, $X/M_X=A_4$;

 $M=\lg a^3\rg\lg c\rg$: $a^3=c^4=1$, $M_X=\lg a^3\rg\lg c^4\rg=1$,  $X/M_X=S_4$;

 $M=\lg a^4\rg\lg c\rg$: $a^4=c^3=1$, $M_X=\lg a^4\rg\lg c^3\rg=1$,  $X/M_X=S_4$.

\vskip 3mm
\f Assume that the result is true for less than $|X|$ and $|X|\gvertneqq 24$.
Then we shall carry out the proof by the following three steps.
\vskip 3mm
{\it Step 1: $M_X\ne 1$}
\vskip 3mm
Suppose that $M_X\ne 1$. Set $M=\lg a^i\rg\lg c\rg$ for some $i$.
Since $a^i\in \cap_{l_2, l_3} M^{a^{l_2}b^{l_3}}=\cap_{l_1, l_2, l_3} M^{c^{l_1}a^{l_2}b^{l_3}}= M_X$,
we get that $M_X=M_X\cap (\lg a^i\rg\lg c\rg )=\lg a^i\rg \lg c^r\rg$ for some $r$.
Set $\ox:=X/{M_X}=\og\oc$.
Then we claim that $\og\cap\oc=1$.
In fact, for any $\olg=\olc'\in \og\cap \oc$ for some $g\in  G$ and $c'\in C$, we have $gc'^{-1}\in M_X$,
that is $g\in \lg a^i\rg$ and $c'\in \lg c^r\rg$, which implies $\olg=\olc'=1$.
Therefore, $\og\cap\oc=1$.
Let $M_0/M_X=\lg \ola^j\rg \lg \olc\rg $ be the biggest subgroup of $\ox$
containing $\lg \olc\rg$ and contained in the subset $\lg \ola\rg \lg \olc\rg $.
Then $\lg \ola^j\rg \lg \olc\rg =\lg \olc\rg \lg \ola^j\rg$.
Since
$$\lg a^j\rg \lg c\rg M_X=\lg a^j\rg M_X\lg c\rg=\lg a^j\rg \lg a^i\rg \lg c\rg \quad {\rm and}\quad \lg c\rg \lg a^j\rg M_X=\lg c\rg M_X\lg a^j\rg =\lg c\rg \lg a^i\rg \lg a^j\rg ,$$
we get  $\lg a^i, a^j\rg \lg c\rg \le X$.
By the maximality of $M$, we have $\lg a^i, a^j\rg =\lg a^i\rg $ so  that $M_0=M$.

Using the induction hypothesis on $\ox=\og\oc$, noting $M_0/M_X=M/M_X$, which is core-free in $\ox$,
we get $\ox$ is isomorphic to  $\ZZ_2,\,D_8,\,A_4$ or $S_4$, and correspondingly,
$o(\ola)=k$, where $k\in\{1, 2, 3, 4\}$, and so $a^k\in M_X$.
Since $M=\lg a^i\rg \lg c\rg$ and $M_X=\lg a^i\rg \lg c^r\rg$, we know that $\lg a^i\rg =\lg a^k\rg$,
which implies that $i\in \{1, 2, 3, 4\}$.
Clearly,  if $\ox=\ZZ_2$, then $M_X=M$; if $\ox=D_8$ and $\o(\olc)=2$, then $M_X=\lg a^2\rg \lg c^2\rg$;
if $\ox =A_4$ and $\o(\olc)=3$, then $M_X=\lg a^2\rg \lg c^3\rg$;
if $\ox =S_4$ and $\o(\olc)=4$, then $M_X=\lg a^3\rg \lg c^4\rg $;
and  if $\ox=S_4$ and  $\o(\olc)=3$, then $M_X=\lg a^4\rg \lg c^3\rg$.
\vskip 3mm
{\it Step 2:  Show that if $M_X=1$ then $G\in \{D_{2kp}\di k=2, 3,4\}$.}
\vskip 3mm
Suppose that $M_X=1$. Since $\lg a\rg_X, \lg c\rg_X\le M_X$, we get $\lg a\rg _X=\lg c\rg_X=1$.
Now we are showing $G_X=1$. For the  contrary, suppose that $G_X\ne 1$.
If $|G_X|\gneqq 4$, then by $G=\lg a, b\rg \cong D_{2n}$ we get $\lg a\rg _X\ne 1$, a contradiction.
So $|G_X|\le 4$.
Since $G_X\lhd G\cong D_{2n}$, we know that $|G:G_X|\le 2$, which implies $|G|\le 8$,
that is $G\cong D_4$ or $D_8$.
A direct checking shows  that $X$ is $D_8$, $A_4$ or $S_4$.
All cases are impossible, as $|X|\gneqq 24$.
In what follows, we consider the faithful (right multiplication) action of $X$
on the set of right cosets  $\O:=[X:\lg c\rg]$.

Suppose that $X$ is primitive.
By Proposition~\ref{Burnside}, every dihedral group is a Burnside group,
which implies that $X$ is 2-transitive. Since $X$ has a cyclic point-stabilizer $\lg c\rg $.
By Lemma~\ref{cd}.(2), we get $G=D_4$ and $X=A_4$, contradicting with $|X|\gneqq 24$.

Suppose  that  $X$  is imprimitive.
Pick a  maximal subgroup $H$ of $X$ which contains $\lg c\rg$ properly.
Then $H=H\cap X=(H\cap G)\lg c\rg=\lg a^s, b_1\rg \lg c\rg\lneqq X$,
for some $b_1\in G\setminus \lg a\rg$ and some $s$.
Using  the same argument as that in Step 1 (viewing  $H$ as $G$), one has $a^s\in H_X$.  Set $\ox=X/H_X$.
Consider the faithful primitive action of $\ox$ on $\O_1:=[\ox:\oh]$,
with a cyclic regular subgroup of $\ola$, where $|\O_1|=s$.
By Proposition~\ref{Burnside}, a cyclic group of composed order is a Burnside group,
we know that either $s$ is a prime $p$ so that $\ox\le \AGL(1,p)$ or $s$ is of composite order
so that $\ox $ is 2-transitive.
In what follows, we consider these two cases, separately.
\vskip 3mm
Case (1):  $a^s=1$.
\vskip 3mm
In this case, $H=\lg c\rg \rtimes \lg b\rg $ and $X=\lg c, b\rg .\lg a\rg $. Then we  have two cases:

\smallskip
Suppose that $s$ is of composite order so that $\ox $ is 2-transitive.
By Proportion {\ref{AGL}}, $\ox\le \AGL(l,q)$ for some prime $q$,
which contains a cyclic regular subgroup $\lg \ola\rg $ of order $q^l$.
By Lemma~\ref{cd}.(1), $\ox\cong S_4$ and $\o(\ola)=4$ so that $\o(a)=4$ (as $H_X\le \lg b, c\rg $),
which in turn implies $G=D_8$.
In this case, checking by Magma,
we have that either $\o(c)=2, 3$ and $|X|\le 24$;
or $\o(c)=4$, $|X|=32$ but $G_X\ne 1$, contradicting with $G_X=1$.

Suppose that $s$ is a prime $p$  so that  $\ox\le \AGL(1,p)$.
Then $\o(a)=p$ so that $G\cong D_{2p}$, where $p\ge 5$, as $|X|\gvertneqq 24$.
Consider the action of $X$ on the set of blocks of length 2 on $\O=[X:\lg c\rg]$, with the kernel, say $K$.
Then $K\nleqq \lg c\rg $ (as $\lg c\rg_X=1$) so that $K$ interchanges two points $\lg c\rg $ and $\lg c\rg b$, which implies $|K/K\cap \lg c\rg |=2$.
Since $K\cap \lg c\rg $ is cyclic and $K\cap \lg c\rg $ fixes setwise each block of length 2, we get $|K\cap \lg c\rg |=2$.
Therefore, $|K|\le 4$.
Since $K\rtimes \lg a\rg \lhd X$ and $p\ge 5$, we have  $K\rtimes \lg a\rg=K\times \lg a\rg $
so that $\lg a\rg \char (K\times \lg a\rg )\lhd X$, contradicting with $G_X=1$.
\vskip 3mm
Case (2): $a^s\neq1$.
\vskip 3mm
Firstly,  show $s=p$, a prime. To do that, we consider the group $X/H_X$.
Since $a^s\in H_X$, we get $H_X\ne 1$ and of course $H_X\nleqq \lg c\rg$.
Suppose that  $\lg a^j\rg \lg c\rg \le H$. Then $a^j\in M$.
Using the same arguments as that in  the first line of Step 1, we get $a^j\in M_X=1$.
Therefore, there exists an $l$ such that $bc^l\in H_X$, which implies $H/H_X=\lg cH_X\rg$,
a cyclic group so that $X/H_X=(\lg aH_X\rg/H_X) (H/H_X)$, a product of two cyclic subgroups, which cannot be isomorphic to $S_4$.
Suppose that  $X/H_X$ is 2-transitive on $\O_1=[X/H_X, H/H_X]$, with a cyclic regular subgroup.
By Lemma~\ref{cd}.(1), $X/H_X\cong S_4$, a contradiction.
Therefore,  $s=p$, a prime and $X/H_X\le \AGL(1, p)$.

Secondly, we consider the group $\oh:=H/\lg c\rg_H=\overline{\lg c\rg \lg a^p,b\rg}$,
taking into account $s=p$, a  prime.
Then $\lg \olc \rg_{\oh}=1$ and $o(\ola^p)=o(a^p)$.
Let $H_0/\lg c\rg _H=\lg \ola^{pj}\rg \lg \olc\rg $ be the biggest subgroup of $\oh$ containing $\lg \olc\rg$  and contained in the subset $\lg \ola^p\rg \lg \olc\rg$.
Since $|\oh|\le |X|$, by the induction hypothesis on $\oh$,
we know is  $H_0/\lg c\rg _H=\lg \ola^{pk}\rg \lg \olc\rg$,
for one of $k$ in $\{2, 3, 4\}$,
which implies $\lg a^{pk}\rg\lg c\rg\lg c\rg_H=\lg c\rg\lg a^{pk}\rg\lg c\rg_H=\lg c\rg \lg a^{pk}\rg$,
giving $\lg a^{pk} \rg \lg c\rg \le H\le X$.
Therefore, we get $a^{pk}\in M_X$.
Since $M_X=1$ and $a^p=a^s\ne 1$, we have
$$a^{pk}=1,\,  {\rm  for\,one\, of}\,   k\in \{2, 3, 4\}.$$
Therefore, only  the following  three groups are remaining:
        $G=D_{2kp}$, where $k\in \{2, 3,4\}$.

\vskip 3mm
{\it Step 3: Show that  $G$ cannot be $D_{2kp}$,  where $k\in\{2, 3, 4\}$, provided $M_X=1$}
\vskip 3mm

Suppose that $G\cong D_{2kp}$,  $k\in\{2, 3, 4\}$,
reminding that $H=\lg a^p, b\rg \lg c\rg$,  $M_X=1$ and $X$ has  blocks of length $2k$.
Moreover, $\lg a^p\rg _X=1$ and there exists no nontrivial element $a^j\in H$
such that that $\lg a^j \rg \lg c\rg \le H$.
Here, we only give the proof for the case $k=4$, that is $G\cong D_{8p}$.
For $k=2$ or $3$, we have the  same arguments but more easer.

Let $G=D_{8p}$. If $p=2$ or $3$, then $G\cong D_{16}$ or $D_{24}$.
These small cases are directly excluded  by Magma.  So assume $p\ge 5$.
Set $a_1=a^p$ and $a_2=a^4$ so that $H=\lg a_1, b\rg \lg c\rg $, where $\o(a_1)=4$.
Set $C_0=\lg c\rg _H$ and $K=H_X$.
Then  $H$ contains an element $a_1$ of order 4 having two orbits of length 4 on each block of length 8,
where $a_1\le K$.
Consider the action of $\oh:=H/C_0=\lg \ola_1, \olb\rg \lg \olc\rg $ on the block
containing the point $\lg c \rg $,  noting that $\lg \olc\rg $ is core-free.
Remind that $\lg a_1\rg _X=1$.
So $\lg \ola_1,\olb\rg\cong D_8$ and $\oh\cong S_4$.
Moreover, we have $N_{\oh}(\lg \olc \rg )=\lg \olc\rg \rtimes \lg \olb\rg \cong D_6$,
by rechoosing $b$ in $\lg a_1, b\rg $.
Therefore  $L:=\lg c\rg \lg b\rg \le X$.

Now we turn to consider the imprimitive action of $X$ on $[X:L]$, which is of degree $4p$.
Let $K\cap \lg c\rg=\lg c_1\rg$.
Then every orbits of $\lg c_1\rg$ on $[X:L]$ is of length 4.
Observing the cycle-decomposition of $c_1X_L\in X/X_L$ on $[X:L]$, we know that $k_1:=\o(c_1X_L)\di 12$.
Therefore, $c_1^{k_1}$ fixes pointwise $[X:L]$, which implies $c_1^{2k_1}$ fixes pointwise $[X:\lg c\rg]$. Therefore, $c_1^{2k_1}=1$ (as $\lg c\rg_X=1$), that is $|K\cap \lg c\rg|\di 2k_1$
and in particular, $|K\cap C_0|\di 2k_1$.
Moreover, since $\lg a_2K\rg \lhd X/K\cong \ZZ_p\rtimes \ZZ_s$ for some $s\di (p-1)$,
we know that $K\rtimes \lg a_2\rg \lhd X$.

Then $|K\cap C_0|\di 24$, as $k_1\di 12$. Also, $K/(K\cap C_0)\cong KC_0/C_0\lhd H/C_0\cong S_4$.  Since $\ola_1\in KC_0/C_0$ where $\o(\ola_1)=4$, and every normal subgroup of $S_4$
containing an element of order 4 must contain $A_4$,
we know that $K/(K\cap C_0)$ contains a characteristic subgroup $K_1/(K\cap C_0)\cong A_4$.
Suppose that  $K\cap C_0\nleqslant  Z(K_1)$.
Then  as a quotient of $A_4$, we have $3\di |K_1/C_{K_1}(K\cap C_0)|.$
However, $K_1/C_{K_1}(K\cap C_0)\le \Aut(K\cap C_0)$, which does not contain an element of order 3,
noting $K\cap C_0\le \ZZ_3\times \ZZ_8$,
by considering the cycle-decomposition of the generator of $K\cap C_0$.
Therefore,  $C_{K_1}(K\cap C_0)=K_1$, that is $(K\cap C_0)\le Z(K_1)$. Since $K_1/(K\cap C_0)\cong A_4$ and $Z(A_4)=1$, we get  $K\cap C_0=Z(K_1) \char K_1\lhd X$. Therefore, $K\cap C_0\le 1$ (as $\lg c\rg_X=1$) so that $K\cong A_4$ or $S_4$, which implies  $\lg a_2\rg \char (K\times \lg a_2\rg)\lhd X$,  contradicting with $G_X=1$  again.
\qed

\vskip 3mm

To handle $X(Q)$, we need the following result.
\begin{lem}\label{G_X}
Suppose that $\lg c\rg_X=1$ and $X=X(G)$ where either $G=D$  and $M=\lg a\rg \lg c\rg$; or $G=Q$. Then $\lg a\rg _X\ne 1$.
\end{lem}
\demo
Since $\lg c\rg_X=1$,  by Proposition~\ref{cyclic}, we have $m\le |G|$.
So $S:=G\cap G^{c}\ne 1$, otherwise $|X|\ge (2n)^2\gneqq |X|$.

(1) Suppose that $G=Q$. Take a subgroup $T$ of order a prime $p$ of $S$.
Since $o(a^jb)=4$ for any $j$, we know $T\le \lg a\rg$.
Since $S$ has the unique element of order $p$,
we get $T^c=T$, giving $T\lhd X$ and so $\lg a\rg_X\ne 1$, as desired.

\smallskip
(2) Suppose that  that $G=D$. Let $M=\lg a\rg \lg c\rg $, where $\o(a)=n$ and $\o(c)=m$.
If $n\ge m$, then by Proposition~\ref{ab}, $\lg a\rg_M\ne 1$
and then  $\lg a\rg_X\ne 1$.
Suppose  that  that  $n+1\le m$.
Since $\lg c\rg_M\ne 1$, we take  $z:=c^{\frac mp}\le \lg c\rg_M$ for a prime $p$.
Since  $\lg c\rg_X\ne 1$, we know that $\lg z^b\rg \ne \lg z\rg $ so that
$N:=\lg z\rg\times \lg z^b\rg \lhd X$. Set $a_1=a^{\frac np}$.
Then $a_1\in N$.

If $p=2$, then $z\in Z(M)$. suppose   that $p$ is odd. Let $N\le P\in \Syl_p(M)$.
By Proposition~\ref{matecyclic},
$P$ is a  metacyclic group and so  we know that $N=\lg a_1\rg \times \lg z\rg $.
Since $\lg z\rg =\lg c\rg_M$, we may set $z^a=z^i$ and $z^b=a_1^jz^l$, where $j\ne 0$.
Then $(z^a)^b=(z^b)^{a^{-1}}=(a_1^jz^l)^{a^{-1}}=a_1^jz^{li^{-1}}$ and $(z^i)^b=a_1^{ji}z^{li}$,
which implies $a_1^{j(i-1)}=1$, that is $i=1$, and so  $z\in Z(M)$ again.
This implies  $z^b\in Z(M)$ and so $N\le Z(M)$.
Thus  $a_1\in Z(M)$ and then $\lg a_1\rg \lhd X$.
\qed

\begin{lem} \label{MQ}
Theorem~\ref{main1} holds, provided  $G=Q$ and $X=X(Q)$.
\end{lem}
\demo
Now $M=\lg a^i\rg \lg c\rg $ for some $i$.
We shall prove the lemma by the induction on $|X|$.
With the same argument as Lemma~\ref{MD}, we get that the conclusion for $|X|\leq24$ holds.
Assume that the result is true for less than $|X|$ and $|X|\gvertneqq 24$.
Then we shall carry out the proof by the following two cases.

Suppose  that  that $\lg c\rg_X\ne 1$, where set $\ox=X/\lg c\rg_X$.
Let $M_1/\lg c\rg _X$ be the biggest subgroup of $\ox$ containing $\lg \olc\rg $
and contained in $\lg\ola\rg \lg \olc\rg $.
By the induction hypothesis, we get $M_1/\lg c\rg _X=\lg\ola^i\rg\lg\olc\rg $,
where  $i\in \{1, 2,3,4\}$.
This gives that $M=\lg a^k\rg \lg c\rg $ where $k\in \{1, 2, 3, 4\}$, as desired.

Suppose  that  that  $\lg c\rg _X=1$.  Then by Lemma~\ref{G_X},  $\lg a\rg_X\ne 1$.
Set  $\ox:=X/\lg a\rg_X=\og \lg \olc\rg $ and
let $M_2/\lg a\rg_X$ be the biggest subgroup of $\ox$ containing $\lg \olc\rg $ and contained in $\lg\ola\rg \lg \olc\rg $.
If $a^n\not\in \lg a\rg_X$, then $\og$ is a generalized quaternion group;
if $a^n\in \lg a\rg_X$, then $\og$ is a dihedral group.
For the first case, by the induction hypothesis on $\ox$;
and for the second case,  by Lemma~\ref{MD},
we get $M_2/\lg a\rg_X=\lg\ola^i\rg\lg\olc\rg$, where $i\in \{1, 2, 3, 4\}$.
Then $\lg a^i\rg \lg c\rg \lg a\rg _X=\lg c\rg \lg a^i\rg \lg a\rg_X$.
This gives $(\lg a^i\rg \lg a\rg_X)\lg c\rg \le X$,
which implies $M=\lg a^k\rg \lg c\rg $ where $k\in \{1,  2, 3, 4\}$, as desired.
\qed

\section{Proof of Theorem~\ref{main2}}
The proof of Theorem~\ref{main2} consists of the following four lemmas.

\begin{lem}\label{2-group}
Suppose that $G\in \{D, Q\},\,X=X(G),\,M=\lg a\rg\lg c\rg$ and $\lg c\rg_X=1$.
If $G$ is a $2-$group, then $\lg a^2\rg \lhd X$.
\end{lem}
\demo
Suppose that $X$ is a minimal counter-example.
Let $a_0$ be the involution of $\lg a\rg $, for both $G\in \{Q, D\}$.
Since $\lg c\rg _X=1$, by Lemma~\ref{G_X}, we get $\lg a\rg_X\ne 1$, which implies $\lg a_0\rg \leq X$,
as $G$ is a 2-group.
Consider $\ox=X/\lg a_0\rg=\og\lg\olc\rg$. Set $\lg \olc \rg _X=(\lg a_0\rg \times \lg c_0\rg )/\lg a_0\rg $.
Then  $\lg c_0^2\rg \lhd X$, which implies $c_0^2=1$.
If $c_0=1$, then by using the induction on $\ox$, we get $\lg\ola^2\rg\lhd \ox$.
Then $\lg a^2\rg \lhd X$ is a contradiction, noting $X$ is  a minimal counter-example.
Therefore,  $o(c_0)=2$.
By the induction on $X/\lg a_0,c_0\rg$,
we get $(\lg a^2\rg (\lg a_0\rg \lg c_0))\rg/\lg a_0,c_0\rg \lhd  X/\lg a_0,c_0\rg$,
that is $H:=\lg a^2\rg \rtimes \lg c_0\rg\lhd X$.
Then we continue the proof by the following two steps.

\vskip 3mm
{\it Step 1:}  Firstly, we shall show that $X$ is a 2-group in this case.
In fact, noting that $\lg a^4\rg=\mho_1(H)\char H\lhd X$,
relabel $\ox=X/\lg a^4\rg $, where write $\lg \olc\rg _{\ox}=\lg \olc^i\rg$.
Then $\lg a^4\rg \rtimes \lg c^i\rg \lhd X$.
Let $Q$ be the $2'$-Hall subgroup of $\lg c^i\rg $.
Since $\Aut(\lg a^4\rg)$ is a $2-$group, we know that $[Q, a^4]=1$ and so $Q\lhd X$,
contradicting with $\lg c\rg_X=1$.
Therefore,  $\lg c^i\rg$ is also a $2-$group.
Reset $\ox=X/\lg a^4\rg \lg c^i\rg=\og \lg \olc\rg$.
Now, $|\og|=8$ and so $\og\cong D_8$ (clearly, $\og$ cannot be $Q_8$).
Since $\lg \olc\rg _{\ox}=1$ we have $\o(\olc)\di 4$.
Therefore, $\ox$ is  a $2-$group and so is $X$.

\vskip 3mm
{\it Step 2:}  Set $K:=\lg a_0\rg\times \lg c_0\rg\cong\ZZ_2^2$.
Consider the conjugacy of $G$ on $K$.
Since  $\lg c_0\rg \ntrianglelefteq X$, we get $C_G(K)\lneqq G$.
Since $G$ may be generated by some elements of the from $a^jb$,
there exists an element $a^ib\in G\setminus C_G(K)$, exchanging $c_1$ and $c_1a_0$ (as $(a^ib)^2=a_0$).
Since $X=GC=(\lg a\rg\lg c\rg).\lg b\rg$, firstly we write  $c^{b}=a^{s}c^{t}$, where $t\neq0$.
Then
$$c=c^{c_0}=c^{b^2}=(a^sc^t)^{b}=a^{-s}(a^sc^t)^t=c^t(a^sc^t)^{t-1},$$
that is   $(a^sc^t)^{t-1}=c^{1-t}$.
Then we have
  $$(c^{t-1})^{b}=(c^{b})^{t-1}=(a^sc^t)^{t-1}=c^{1-t}.$$
If $t\neq1$, then $c_1^b\in\lg c^{t-1}\rg^b=\lg c^{t-1}\rg$, contradicting with $c_1^b=a_1c_1$.
So $t=1$, that is  $c^b=a^sc$.
Secondly, we write  $c^b=c^{t_1}a^{s_1}$. With the same arguments, we may get $t_1=1$ and  $c^b=ca^{s_1}$.
Therefore,  we have  $a^{s}c=c^{b}=ca^{s_1},$ that is  $(a^{s})^{c}=a^{s_1}$.
Clearly $\lg a^s\rg =\lg a^{s_1}\rg$, that is $c$ normalises $\lg a^{s}, b\rg$.
Then
  $$\lg a^{s}, b\rg\leq  \cap_{c^i\in \lg c\rg} G^{c^i}=\cap_{x\in X} G^x=G_X,$$
which implies $b^a=ba^{-2}\in G_X$  so that $a^2\in G_X$ and then $\lg a^2\rg \lhd X$.
This contradicts the minimal of $X$.
\qed

\begin{lem}\label{Da^2}
Suppose that $G=D,\,X=X(D),\,M=\lg a\rg \lg c\rg$ and $\lg c\rg_X=1$.
Then $\lg a^2\rg\lhd X$.
\end{lem}
\demo
By Lemma~\ref{2-group}, the lemma is true when $G$ is a 2-group.
So assume that $G$ is not a 2-group.
Take a minimal counter-example $X$.
In the following Step 1, we show that the possible groups for $G$ are $D_{2^e p^k}$,
where $p$ is a prime and $e\in\{ 1, 2\}$;
and in Step 2, we show that $G$ cannot be these groups.

\vskip 3mm
{\it Step 1: Show that the possible groups for $G$ are  $D_{2^e p^k}$, where $p$ is an odd prime and $e\in\{ 1, 2\}$.}
\vskip 3mm
By Lemma~\ref{G_X}, let $p$ be the maximal prime divisor of $|\lg a\rg_X|$ and set $a_0=a^{\frac np}$.
Set $\ox=X/\lg a_0\rg=\og\lg \olc \rg$ and $\lg \olc \rg_{\ox}=\lg \olc_0\rg $.
(i) Suppose  that  that $\lg \olc \rg_{\ox}=1$.
Then by the minimality of $X$ we  get $\lg\ola^2\rg\lhd \ox$,
which implies $\lg a^2\rg \lg a_0\rg \lhd X$.
Since $\lg a^2\rg \lg a_0\rg =\lg a^2\rg$ or $\lg a\rg $, we get $\lg a^2\rg \lhd X$, a contradiction.
(ii)Suppose  that  that $\lg \olc\rg \lhd \ox$.
Then $\ox/C_{\ox}(\lg \olc\rg)\le\Aut(\lg \olc\rg)$,
which is abelian and so $\ox'\le C_{\ox}(\lg \olc\rg)$.
Then $\ola^2\in \og'\leq \ox'\leq C_{\ox}(\lg \olc\rg)$,
that is $[a^2, c]\in\lg a_0\rg$, which implies $\lg a^2, a_0\rg \lhd X$,
and again we have $\lg a^2\rg \lhd X$ is a contradiction.
By (i) and (ii), we have  $1\ne \lg \olc\rg_{\ox}=\lg\olc_0\rg \lneqq\lg \olc\rg$.
Reset
$$K=\lg a_0\rg\rtimes\lg c_0\rg,\,\ox=X/K=\og\lg \olc\rg,\,H=\lg a^2\rg\rtimes\lg c_0\rg.$$
If $\o(a_0)<\o(c_0)$, then $\{1\}\subsetneqq \lg c_0^j\rg=Z(K)\lhd X$ is, for some $j$, a contradiction. Therefore, $1< \o(c_0)\le \o(a_0)$.
Then  we have the following two cases:
\vskip 3mm
{\it Case 1: $K=\lg a_0\rg \rtimes \lg c_0\rg \cong \ZZ_p\rtimes \ZZ_r$, a Frobenius group, where $r\ge 2$.}
\vskip 3mm
In this case, $p$ is odd.
Set  $\ox=X/K $.
By the minimality of $X$, we have $H/K=\lg \ola^2\rg\lhd \ox$,
that is $H:=\lg a^2\rg \rtimes \lg c_0\rg \lhd X$.
Since $K\lhd X$, we know that $\lg a^2\rg/\lg a_0\rg$ and $\lg c_0\rg \lg a_0\rg/\lg a_0\rg$ are normal
in $H/\lg a_0\rg$.
Then $[a^2, c_0]\le \lg a_0\rg $.
So one can write
  $$H=\lg a^2,c_0|a^n=c_0^r=1,\,(a^2)^{c_0}=a^2a_0^j\rg.$$
Let $P\in\Syl_p(H)$. Then $P\char H\lhd X$ so that $P\leq \lg a\rg_X$.
Clearly, one can check $Z(H)=\lg a^{2p}\rg$.
Then $\lg a^{2p}\rg\leq \lg a\rg_X$.
Note $\lg a^{2p},P\rg=\lg a^{2p}, a^{n/p^k}\rg =\lg a^2\rg $, where $p^k\di\di n$,
so that  $a^2\in \lg a \rg_X$ is a contradiction again.
\vskip 3mm
{\it Case 2: $K=\lg a_0\rg \times \lg c_0\rg \cong \ZZ_p^2$.}
\vskip 3mm
Set $H=\lg a^2, c_0\rg$ again.
With the same reason as that in Case 1, we have $H\lhd X$.

Suppose that $a_0\notin\lg a^2\rg$. Then $p=2$ and $\frac n2$ is odd.
Noting $\lg a^2\rg$ is $2'-$Hall subgroup of $H$, we have $\lg a^2\rg\char H\lhd X,$
which implies $a^2\in \lg a\rg_X$, a contradiction.

Suppose  that $a_0\in\lg a^2\rg$. Let $H_1$ be the $p'-$Hall subgroup of $H$.
Then $H_1$ is also the $p'-$Hall subgroup of $\lg a^2\rg$.
Then $H_1\lhd X$, which implies $H_1\leq \lg a\rg_X$.
Suppose  that  $H_1\neq1$.
Let $a_2$ be an element of order $q$ in $H_1$, where $q<p$ is a prime as the maximality of $p$.
Consider $\ox:=X/\lg a_2\rg=\og\lg c\rg $.
Similarly, we have $1\ne\lg\olc\rg_{\ox}:=\lg\olc_2\rg\lneqq\oc$ and
$H_0:=\lg a^2\rg\rtimes\lg c_2\rg\lhd X$.
Let $P\in \Syl_p(H_0)$.
Then $P\char H$ and so $P\lhd X$, which implies $P\leq \lg a\rg_X$.
Noting $\lg H_1,P\rg=\lg a^2\rg$, we therefore get $a^2\in \lg a\rg_X$, a contradiction.
So $H_1=1$, which means that $G$ is $D_{2^ep^k}$ where $p$ is an odd prime and $e=1, 2$,
as $G$ is assumed to be not a 2-group.
\vskip 3mm
{\it Step 2: Show that $G$ cannot be $D_{2^ep^k}$, where $e\in \{1, 2\}$ and $p$ is an odd prime.}
\vskip 3mm
Relabel $a^2$ by $a_2$.
Following Proposition~\ref{pcyclic},  we write
$$H=\lg a_2, {c_0}\rg =\lg a_2^{p^k}={c_0}^p=1, a_2^{c_0}=a_2^{1+k'p^{k-1}}\rg,$$
where $k\ge 2$ and $k'$ may be 0.
We shall get a contradiction by  considering two groups.
\vskip 3mm
(1) The cyclic extension $H\lg c\rg =\lg a_2\rg \lg c\rg$ of $H$ by $c$, where $c^l={c_0}$ for some $l$.
\vskip 3mm
Let $\pi$ be  the automorphism  of $H$ by mapping $a_2$ to $a_2^i{c_0}^j$
and $c_0$ to $c_0$, where $i, j\not\equiv 0(\mod p).$
Then we have
\vskip 3mm
(i) $\pi$ preserves the relation of $H$:  that is
$$(a_2^i{c_0}^j)^{p^k}=a_2^{ip^k}=1, \quad (a_2^i{c_0}^j)^{p^{k-1}}=a_2^{ip^{k-1}}\ne 1,$$
$$\pi(a_2)^{1+k'p^{k-1}}=(a_2^i{c_0}^j)^{1+k'p^{k-1}}=
  a_2^i{c_0}^j(a^{ik'p^{k-1}}{c_0}^{p^{k-1}})=a_2^{i(1+k'p^{k-1})}{c_1}^j=(\pi(a_2))^{\pi ({c_0})}.$$
\vskip 3mm
(ii) $\pi^l=\Inn({c_1}):$
$$\pi^l(a)=(a^i{c_0}^j)^l=a^{i^l}{c_0}^{j\sum_{w=0}^{l-1}i^w}a_1^x=\Inn({c_0})(a)=a^{1+k'p^{k-1}},$$
that is
$$i^l+xp^{k-1}\equiv 1+k'p^{k-1}(\mod p^k), \quad \sum_{w=0}^{l-1}i^w\equiv 0(\mod p).$$
So if $i=1(\mod p)$ ($i-1$ is of order a $p$-power), then  $p\di l$;
if $i\ne 1(\mod p)$ then $l=up^v$, where $u\di (p-1).$
Now we have $\o(c)=pl$ and $\o(i)=pl$ or $l$.
\vskip 3mm
(2) The group $X=\lg H\lg c\rg,b\rg$.
\vskip 3mm
Set
 $${c_0}^a=a^{m_2p^{k-1}}{c_0}^{s'},\quad {c_0}^b=a^{m_1p^{k-1}}{c_0}^s,\quad  c^b=a^rc^t,$$
where $m_2\ne 0$.
Then we are showing $t=1$.
In fact, since $b$ preserves $a_2^c=a_2^i{c_0}^j$, we have that in $\ox=X/\lg a_1\rg$, $\ola_2^{-i^t}=\ola_2^{-i},$ that is $i^{t-1}\equiv 1{\pmod {p^{k-1}}}$,
which implies $i^{p(t-1)}\equiv 1(\mod p^k)$ and so either $l\di (t-1)$ or $l\di p(t-1)$.
If $l\di (t-1)$, then $c^{t-1}\le \lg c^l\rg =\lg {c_0}\rg$.
Since $m_1\ne 0$, we have $t=1$.
Suppose that $l\di p(t-1)$ and $l\nmid (t-1)$.
Then ${c_0}\in\lg c^{t-1}\rg$.
Since $$c=c^{b^2}=(a^rc^t)^b=a^{-r}(a^rc^t)^t=c^t(a^rc^t)^{t-1},$$
we have $(a^rc^t)^{t-1}=c^{1-t}$, which implies $(c^{t-1})^b=(a^rc^t)^{t-1}=c^{1-t}.$
Then $b$ normalises $\lg c^{1-t}\rg$, which implies that $b$ normalises $\lg {c_0}\rg$, forcing $t=1$, that is $c^b=a^rc$.

Similarly, we may set $c^b=c^{t_1}a^{r_1}$ and get $t_1=1$ and $c^b=ca^{r_1}$.
Therefore, we have $ca^{r_1}=c^b=a^rc$, that is $(a^{r})^c=a^{r_1}.$
Then
$$\lg a^r, b\rg\leq  \cap_{c^i\in \lg c\rg} G^{c^i}=\cap_{x\in X} G^{x}=G_X.$$
Since any normal subgroup of $G\cong D_{2n}$ containing $b$ is either $\lg a^2, b\rg $ or $G$,
we get $\lg a^2\rg \lhd X$, a contradiction.
\qed

\begin{lem}\label{a}
Suppose that $G\in \{Q, D\},\,X=X(G),\,\lg c\rg _X=1$ and $M=\lg a\rg\lg c\rg$.
If $\lg a\rg\lhd X$, then $G\lhd X$.
\end{lem}
\demo
$X=(\lg a\rg\rtimes\lg c\rg).\lg b\rg,$ and so we may write $a^c=a^i$ and $c^b=a^kc^j$.
If $j=1$, then $G\lhd X$. So assume $j\ne 1$.
Since $b^2=a^n$ ($o(a^n)=1$ or 2, if $G=D$ or $G=Q$, resp.) and $\lg a^n\rg\lhd X$,
we get $a^n\in Z(X)$, which implies $c=c^{b^2}$.
Then
$$c=(c^b)^b=(a^kc^j)^b=a^{-k}(a^kc^j)^j=(c^ja^k)^{j-1}c^j,$$
that is $c^{1-j}=(c^ja^k)^{j-1}=(c^{j-1})^b,$ so that $b$ normalizes $\lg c^{1-j}\rg$.
Since $\ox=X/C_X(\lg a\rg)\le \Aut(\lg a\rg )$ which is abelian,
we get $\olc =\olc^{\olb}=\olc^j$, that is $c^{1-j}\le C_X(\lg a\rg)$ so that $[c^{1-j}, a]=1$.
Thus we get  $\lg c^{1-j}\rg \lhd X$.
It follows from $\lg c\rg _X=1$ that $j=1$, a contradiction.
\qed

\begin{lem}\label{Qa^2}
Suppose that $G=Q,\,X=X(Q),\,M=\lg a\rg \lg c\rg$ and $\lg c\rg_X=1$.
Then $\lg a^2\rg\lhd X$.
\end{lem}
\demo
By Lemma~\ref{2-group}, we just consider the case when $G$ is not a 2-group.
Take a minimal counter-example $X$ and set $a_1:=a^n.$
Similarly, we carry out the proof  by the following two steps:
\vskip 3mm
{\it Step 1: Show that the possible groups for $G$ are $Q_{4p^k}$, where $p$ is an odd prime.}
\vskip 3mm
By Lemma~\ref{G_X}, let $p$ be the maximal prime divisor of $|\lg a\rg_X|$ and $a_0=a^{\frac np}$.
Then $G/\lg a_0\rg $ is either a generalized quaternion group or a dihedral group.
Let $K=\lg a_0\rg\rtimes\lg c_0\rg$ such that $K/\lg a_0\rg $
is the core of $\lg a_0, c\rg/\lg a_0\rg$ in $X/\lg a_0\rg$ and set $H=\lg a^2\rg \rtimes \lg c_0\rg$.
Using  completely same arguments as that in Lemma~\ref{Da^2},
one may get $G\cong Q_{4p^k}$, where $p$ is an odd prime.
Remind that  Lemma~\ref{Da^2} is  used when we make an induction for   $G/\lg a_0\rg $   being a dihedral group.

\vskip 3mm
{\it Step 2: Case $G\cong Q_{4p^k}$, where $p$ is an odd prime.}
\vskip 3mm
Suppose that $a_1\in \lg a\rg_X$.
Let $\lg a_1\rg\rtimes\lg c_1\rg/\lg a_1\rg$ be the core of $\lg a_1,c\rg/\lg a_1\rg$
in $X/\lg a_1\rg $.
Since $\lg c_1^2\rg \lhd X$, we get $c_1^2=1$.  Consider $\ox=X/\lg a_1\rg\lg c_1\rg$.
Since $\og\cong D_{2p^k}$ is a dihedral group, by Lemma~\ref{Da^2},
we get $\lg\ola^2\rg=\lg\ola\rg\lhd \ox$, which implies $\lg a\rg\rtimes\lg c_1\rg\lhd X$.
Then $\lg a^2\rg\lhd X$, a contradiction.
So in what follows, we assume $a_1\notin\lg a\rg_X$, that is $\lg a\rg_X$ is $p$-group.
Then we continue the proof by two substeps.
\vskip  3mm
{\it Substep 2.1: Show that the possible values of $m$ are $pq^e$, for a prime $q$ (may be equal to $p$)  and an integer $e$.}
\vskip 3mm

For the contrary, we assume that $m=pq^em_1$ where $m_1\ne1$, $p\nmid m_1$ and $q\nmid m_1$.
Since $H=\lg a^2, c_0\rg $ and $b^2=a_1$, we get $\ox=X/H=\lg \olb\rg\lg\olc\rg$.
By considering the permutation representation of $\ox$ on the cosets $[\ox:\lg \olc\rg ]$ of size 4,
we know that $\lg\olc^2\rg\lhd X$.
So $\lg b, c^2, H\rg \le X$, that is  $X_1:=\lg a,b\rg\lg c^2\rg=G\lg c^2\rg \le X$.

Firstly, suppose  that $m$ (=$o(c)$) is even.
Then $[X:X_1]=2$.
Let $\lg c_2\rg$ be the Sylow 2-subgroup of $\lg c\rg$.
By the induction on $X_1$, $\lg a^2\rg\lhd X_1$ and in particular, $c^2$ normalizes $\lg a^2\rg $.
Since   $m=pq^em_1$ has other prime divisors distinct with 2 and $p$, we get $X_2:=H(\lg b\rg\lg c_2\rg)\lneqq X.$
By the  induction on $X_2$ again,  $\lg a^2\rg\lhd X_2$, which implies $\lg c_2\rg$ normalises $\lg a^2\rg$. In summary,
$\lg c_2, c^2\rg $ normalizes $\lg a^2\rg $. Since $\lg b\rg\lg c_2\rg$ is a Sylow $2-$subgroup of $X$, we have $\lg c_2, c^2\rg =\lg c\rg$,
  that is $\lg a^2\rg \lhd X$, a contradiction.

Secondly, suppose  that $m$ is odd.
Then  both $q$ and $m_1$ are odd,  so that $X=X_1=(\lg a^2\rg\rtimes\lg c\rg)\rtimes\lg b\rg$.
By the induction on  $X_3:=\lg H, c^{\frac{m}{m_1}}\rg =\lg a,b\rg\lg c^{\frac{m}{m_1}}\rg\lneqq X$
and $X_4:=\lg H, c^{\frac{m}{pq^e}}\rg=\lg a,b\rg\lg c^{\frac{m}{pq^e}}\rg\lneqq X$, respectively,
we get both $\lg c^{\frac{m}{m_1}}\rg$ and $\lg c^{\frac{m}{pq^e}}\rg$ normalise $\lg a^2\rg$.
Noting $\lg c^{\frac{m}{m_1}},c^{\frac{m}{pq^e}}\rg=\lg c\rg$, we get $\lg a^2\rg\lhd X$,
a contradiction again.
\vskip  3mm
{\it Substep 2.2: Exclude the case $m=pq^e$, for a prime $q$  and an integer $e\ge1$.}
\vskip 3mm

Recall $a_1=a^n$, the unique involution in $G$,
$\lg a_0\rg $ is a normal subgroup of order $p$ in $X$,
$\lg a_0\rg\lg c_0\rg/\lg a_0\rg$ is  the core of  $\lg a_0\rg \lg c\rg/\lg a_0\rg$ in $X/\lg a_0\rg $,
$H=\lg a^2\rg\rtimes \lg c_0\rg \lhd X$ (by the induction hypothesis)
and $X=((H.\lg c\rg).\lg a_1\rg).\lg b\rg$.

Noting $\lg a_0\rg\lg c_0\rg=(\lg a_0\rg \lg c\rg)_X$ and $\lg c\rg_X=1$, we get that
$\lg a_0\rg\lg c_0\rg$ is either $\ZZ_p^2$ or  $\ZZ_p\rtimes\ZZ_{q^{e_1}}$
where $q^{e_1}\di p-1$ and $e_1\leq e$.
If $\lg a_0\rg\lg c_0\rg\cong\ZZ_p\rtimes\ZZ_{q^{e_1}}$ for some $q<p$ and $e_1\leq e$, then we get $\lg a^2\rg\char H\lhd X$, a contradiction.
Therefore, we get $\lg a_0, c_0\rg\cong\ZZ_p^2$, which implies $\lg c_0\rg=\lg c^{q^e}\rg\leq\lg c\rg$.
Note that $\lg a^{2p}\rg=\mho_1(H)\char H\lhd X$. Thus $\lg a^{2p}\rg\lhd X$.

Set $X_5=(H.\lg c^q\rg)\rtimes\lg b\rg =\lg a,b\rg\lg c^q\rg < X$.
By the induction on $X_5$, we get $\lg a^2\rg\lhd X_5$,
that is $X_5=(\lg a^2\rg\rtimes\lg c^q\rg)\rtimes\lg b\rg$.
Clearly, $\lg a^2\rg \le  G'\le X_5'\le \lg a^2, c^q\rg$.
So set $X_5'=\lg a^2, c_3\rg $ for $c_3\in \lg c^q\rg$.
By Proportion~\ref{NC}, both $X/C_{X}(\lg a^{2p}\rg)$ and $X_5/C_{X_5}(\lg a^2\rg)$ are abelian,
which implies that $X'\leq C_X(\lg a^{2p}\rg)$ and $X_5'\leq C_{X_5}(\lg a^2\rg)$.
Note that $\lg a^2\rg\leq X_5'$. Thus $X_5'$ is abelian.
The $p'$-Hall subgroup of $X_5'$ is normal, contradicting with $\lg c^q\rg_{X_5}=1$,
meaning that $X_5'$ is an abelian $p$-group.
Set $L:=H\rtimes\lg a_1\rg\lneqq X_5$.

Suppose   that $L\lhd X$.
If $H$ is abelian, then we get that either $\lg a^2\rg =Z(L)\char L\lhd X$, a contradiction;
or $L$ is abelian, forcing $\lg a_1\rg\char L\lhd X$, a contradiction again.
Therefore, $H$ is non-abelian.
Note that $X_5'=\lg a^2, c_3\rg $ for $c_3\in \lg c^q\rg$.
If $c_3\ne1$, then  $c_0\in\lg c_3\rg\leq X_5'$ as $\o(c_0)=p$,
which implies that $H=\lg a^2,c_0\rg$ is abelian, a contradiction.
Therefore, $X_5'=\lg a^2\rg$, which implies $L=\lg a\rg\rtimes\lg c_0\rg$.
Note that $\lg a_1\rg\char L\lhd X$.
Thus $\lg a_1\rg\lhd X$, a contradiction.

Suppose   that $L\ntrianglelefteq X$.
Then in $\ox=X/H=(\lg\olc\rg\rtimes\lg \overline{a_1}\rg).\lg \olb\rg$,
we get that either $\olc^{\overline{a_1}}=\olc^{-1}$ if $q$ is odd; or
$\olc^{\overline{a_1}}$ is either $\olc^{-1}$ or $\olc^{\pm1+2^{e-1}}$ if $q=2$.
Then we shall divide it into the following two cases:
\vskip 3mm
{\it Case 1: $q$ is an odd prime.}
\vskip 3mm
In this case, $q$ is odd.
Since $\olc^{\overline{a_1}}=\olc^{-1}$ in $\ox=X/H$,
we get $\lg a^2, c^{qp}\rg\leq X_5'\leq \lg a^2\rg\lg c^q\rg$.
Note that  $X_5'$ is the abelian $p$-group.
Thus either $q\ne p$ and $e=1$; or $q=p$.
Suppose that $q\ne p$ and $e=1$, that is $\o(c)=pq$.
Consider $M=\lg a\rg\lg c\rg\lhd X$.
Then by Proportion~\ref{mateabel}, $M'$ is abelian.
Note that $\lg c\rg_X=1$ and $G_X$ is the $p-$group.
Thus $M'$ is an abelian $p$-group with the same argument as the case of $X_5'$.
Noting $\lg a_1\rg\lg c^p\rg$ is the $p'$-Hall subgroup of $M$,
we get $[a_1,c^p]\in \lg a_1\rg\lg c^p\rg\cap M'=1$, which implies $\olc^{\overline{a_1}}=\olc$ in $\ox=X/H$, a contradiction.
So in what follows, we assume $q=p$, that is $\o(c)=p^{e+1}$.
Note that $\olc^{\overline{a_1}}=\olc^{-1}$ in $\ox=X/H$ and $\lg a^2\rg\leq X_5'$.
Thus $X_5'$ is either $\lg a^2\rg\lg c^p\rg$ or $\lg a^2\rg$,
noting $X_5'=\lg a^2\rg$ only happens when $e=1$.

Suppose that $X_5'=\lg a^2\rg\lg c^p\rg$. Note that $H\leq X_5'$.
Thus $H=\lg a^2\rg\rtimes\lg c_0\rg$ is abelian.
Note that both $\lg a^2\rg$ and $\lg c\rg$ are $p$-groups and $X=(H.\lg c\rg)\rtimes\lg b\rg$.
Set $(a^2)^{c}=a^{2s}c_0^{t}$ and $c^b=a^{2u}c^v$ where $s\equiv1(\mod p)$ and $p\nmid v$.
Then for an integer $w$, we get $(a^2)^{c^w}=a^{2x_1}c_0^{wt}$ and $c_0^b=a^{x_2}c_0^{v}$
for some integers $x_1$ and $x_2$.
Since $((a^2)^{c})^b=(a^{2s}c_0^{t})^b$, there exist some integers $x$ and $y$ such that
$$((a^2)^{c})^b=(a^{-2})^{c^v}=a^{x}c_0^{-vt}\quad{\rm {and}}\quad
((a^2)^sc_0^{t})^b=a^yc_0^{vt},$$
which gives $t\equiv0(\mod p)$.
Then $\lg a^2\rg\lhd X$, a contradiction.

Suppose that $X_5'=\lg a^2\rg$. Then $\o(c)=p^2$ and $c_0=c^p$.
Note that $X_5=G\rtimes\lg c^q\rg=H\rtimes\lg b\rg$.
Thus $[c_0,a_1]=1$.
Set $c^{a_1}=a^xc^{-1+yp}$ as $\olc^{\overline{a_1}}=\olc^{-1}$ in $\ox=X/H$.
Then
$$c=c^{a_1^2}=a^x(a^xc^{-1+yp})^{-1+yq}=a^xc^{1-yp}a^{-x}c^{yp},$$
which implies $(a^x)^{c^{1-yp}}=a^x$.
Then $[a^x,c]=1$.
Note that $c_0=c^p$ and $c_0=c_0^{a_1}=(c^{a_1})^p=a^{x}c^{-p}$ for some $x$.
Thus we get $c_0^2=1$, contradicting with $\o(c_0)=p$.
\vskip 3mm
{\it Case 2: $q=2$.}
\vskip 3mm
In this case, we know that $X_5=(\lg a^2\rg\lg c^2\rg)\rtimes\lg b\rg\lhd X$
and $\olc^{\overline{a_1}}$ is either  $\olc^{-1}$ or $\olc^{\pm1+2^{e-1}}$ in $X/H$,
noting that $\olc^{\overline{a_1}}=\olc^{\pm1+2^{e-1}}$ only happens when $e\ge2$.
By Proportion~\ref{sylowp},
we get that both $\lg b\rg\lg c^p\rg$ and $\lg c^{2p}\rg\lg b\rg$ are Sylow 2-subgroups of $X$ and $X_5$, separately.
Note that $a^2\in X_5'$, $X_5'\char X_5\lhd X$ and $X_5'$ is the abelian $p$-group.
Thus $X_5'=H$, which implies that $H$ is ablian.
Note that $\lg c^{2p}\rg\lg b\rg$ is the Sylow 2-subgroup of $X_5$.
Thus $[c^{2p},b]\in X_5'\cap \lg c^{2p}\rg\lg b\rg=1$, which implies that $[c^{2p},a_1]=1$ and $\lg c^{2p}\rg\lg b\rg$ is abelian.
Then $\olc^{\overline{a_1}}=\olc^{1+2^{e-1}}$ in $X/H$.
Note that $\lg b\rg\lg c^p\rg=(\lg c^p\rg\rtimes \lg a_1\rg).\lg b\rg$.
Then $(c^p)^{a_1}=c^{p+2^{e-1}p}$, which implies $c^{2^{e-1}p}\in M'$.

Suppose that $e=2$.
Since $|\lg b\rg|\geq|\lg c^p\rg|$, by Proportion~\ref{ab}, we get $\lg a_1\rg\lhd\lg b\rg\lg c^p\rg$.
Note that $b^2=a_1$ is an involution. Thus $[a_1,c^p]=1$, a contradiction.
So in what follows, we assume $e>2$.

Noting $\lg c^p\rg\rtimes\lg a_1\rg=\lg a_1,c^p|a_1^2=c^{2^ep}=1,\,(c^p)^{a_1}=c^{(1+2^{e-1})p}\rg$,
there are only three involutions in $\lg c^p\rg\rtimes\lg a_1\rg$: $a_1,c^{2^{e-1}p}$ and $a_1c^{2^{e-1}p}$.
By considering the permutation representation of $\lg b\rg\lg c^q\rg$ on the cosets $[\lg b\rg\lg c^p\rg:\lg c^p\rg ]$ of size 4, we know that  $\lg c^{2p}\rg\lhd X$.
By Proportion~\ref{mateabel}, $M'$ is abelian.
Let $M_2$ be the Sylow 2-subgroup of $M'$.
Note that $M_2\char M'\char M\lhd X$, $c^{2^{e-1}p}\in M'$, $\lg c\rg_X=1$ and $a_1$ is an involution.
Thus we get $M_2\cong \ZZ_2^2$, which implies $M_2=\lg c^{2^{e-1}p}, a_1\rg$.
Consider $HM_2\leq X$.
Since $M_2\lhd X$, $H\lhd X$, $H\cap M_2=1$ and $p$ is odd prime,  we get $HM_2=H\times M_2\lhd X$,
which implies $a$ normalises $\lg c^{2^{e-1}p}\rg$.
Since  $X=\lg a,b,c\rg$ and $[b, c^{2p}]=1$, we get $\lg c^{2^{e-1}p}\rg\lhd X$, a contradiction.
\qed

\begin{rem} If the readers  are familiar with the theorem of regular maps, then  two known results in  \cite[Theorem 1.3]{ZD2016} and \cite[Lemma 5]{HR2022}  on skew-morphisms can be used so that
 one  may give a short  proof for Theorem~\ref{main2}, see below.
\end{rem}
{\bf Proof:} (1) Suppose that $\lg c\rg_X=1$.
For the  five cases  of Theorem~\ref{main1}, we know  that $M_X$ is
$M,\,\lg a^2\rg\lg c^2\rg,\,\lg a^2\rg\lg c^3\rg,\,\lg a^3\rg\lg c^4\rg$ or $\lg a^4\rg\lg c^3\rg$,
respectively.
Suppose $M_X=\lg a^i\rg\lg c^j\rg$, which is of the above four cases.
Set  $X_1=GM_X=\lg a,b\rg\lg c^j\rg$.
Then  $\lg c^j\rg_{X_1}=1$.
By \cite[Theorem 1.3]{ZD2016} and \cite[Lemma 5]{HR2022},  we get $\lg a^2\rg\lhd X_1$, that is $c^j$ normalizes $\lg a^2\rg$, and so $M_X\cap \lg a^2\rg \lhd M_X$.
\vskip 3mm
(2) Suppose $\lg c_1\rg:=\lg c\rg_X\ne1$.
Then by Proposition~\ref{NC}, we get $\lg c\rg\leq C_X(\lg c_1\rg)\lhd X$ and
$\ox=X/C_X(\lg c_1\rg)=\lg\ola,\olb\rg$ is abelian.
This implies $\lg\ola,\olb\rg \lessapprox D_4$.
Therefore, $\lg a^2, c\rg\leq C_X(\lg c_1\rg)$.
Therefore, $|X:C_X(\lg c_1\rg))|\le 4$.
\qed

\section{Proof of  Theorem~\ref{main3}}

To prove Theorem 1.3, set $R:=\{a^{2n}=c^m=1,\,b^2=a^n,\,a^b=a^{-1}\}.$
Then we shall deal with the five cases in Theorem 1.1  in the following five subsections, separately.
Let $A=G.\lg t\rg $ where $G\lhd A$ and $t^l=g\in G$.
Then $t$ induces an automorphism $\t$ of $G$ by conjugacy.
Recall  that  by the cyclic extension theory of groups,
this extension is valid if and only if
$$\t^l=\Inn(g)\quad  {\rm and} \quad \t(g)=g.$$

\subsection{$M=\lg a\rg \lg c\rg $}
\begin{lem}\label{Qac}
Suppose that $X=X(Q)$, $M=\lg a\rg \lg c\rg $ and $\lg c\rg _X=1$. Then
\begin{eqnarray}\label{m1}
X=\lg a,b,c| R, (a^2)^c=a^{2r},  c^a=a^{2s}c^t,c^b=a^uc^v\rg,
\end{eqnarray}
where
 $$\begin{array}{ll}
  &r^{t-1}-1\equiv r^{v-1}-1\equiv0(\mod n), t^2\equiv 1(\mod m),\\
  &2s\sum_{l=1}^tr^{l}+2sr\equiv 2sr+2s\sum_{l=1}^{v}r^l-u\sum_{l=1}^tr^l+ur
  \equiv 2(1-r)(\mod 2n),\\
  &{\rm if}\, 2\di n,\, {\rm then}\, u(\sum_{l=0}^{v-1}r^l-1)\equiv 0(\mod 2n)\,{\rm and}\,v^2\equiv 1(\mod m),\\
  &{\rm if}\, 2\nmid n,\, {\rm then}\,  u\sum_{l=1}^vr^l-ur\equiv 2sr+(n-1)(1-r)(\mod 2n),\, {\rm and}\,v^2\equiv t(\mod m),\\
  &{\rm if}\, t\ne 1,\, {\rm then}\, u\equiv0(\mod2),\\
  &2s\sum_{l=1}^{w}r^l\equiv u\sum_{l=1}^{w}(1-s(\sum_{l=1}^tr^l+r))^l\equiv 0(\mod 2n)\Leftrightarrow
  w\equiv 0(\mod m).
  \end{array}$$
\end{lem}
\demo
Noting $\lg a^2\rg\lhd X$ and $M=\lg a\rg\lg c\rg\le X$,
we have $X$ may be obtained by three cyclic extension of groups in order:
$$\lg a^2\rg\rtimes \lg c\rg, \quad (\lg a^2\rg\rtimes \lg c\rg).\lg a\rg \quad {\rm {and}}\quad
((\lg a^2\rg\rtimes \lg c\rg).\lg a\rg) \rtimes \lg b\rg.$$
So $X$ has the presentation as in Eq(\ref{m1}).
What we should to determine the parameters $r,s,t,u$ and $v$ by analysing three extensions.

(1) $\lg a^2\rg\rtimes \lg c\rg$, where $(a^2)^c=a^{2r}$.
Set $\pi_1\in \Aut(\lg a^2\rg)$ such that $\pi_1(a^2)=a^{2r}$.

As mentioned before, this extension is valid if and only if
$\o(\pi_1(a^{2}))=\o(a^2)$ and $\pi_1^{m}=1$, that is
\begin{eqnarray}\label{f3.1}
r^m-1\equiv 0(\mod n).
\end{eqnarray}

(2) $(\lg a^2\rg\rtimes \lg c\rg).\lg a\rg$, where $c^a=a^{2s}c^t$. Set $\pi_2\in \Aut((\lg a^2\rg\rtimes \lg c\rg)$: $a^2\to a^2$ and $c\to a^{2s}c^t$.
This extension is valid if and only if the following three equalities hold:

(i) $\pi_2$ preserves $(a^2)^c=a^{2r}$:
\begin{eqnarray}\label{f3.2}
r^{t-1}-1\equiv  0(\mod n).
\end{eqnarray}

(ii) $\o(\pi_2(c))=m$:
  $$(a^{2s}c^t)^m=c^{tm}(a^{2s})^{c^{tm}}\cdots(a^{2s})^{c^t}
  =c^{tm}a^{2s\sum_{l=1}^mr^{tl}}=c^{tm}a^{2s\sum_{l=1}^m r^{l}}=1,$$
that is
\begin{eqnarray}\label{f3.3}
s\sum_{l=1}^{m} r^{l}\equiv  0(\mod n).
\end{eqnarray}

(iii) $\pi_2^2=\Inn(a^2):$
$$ca^{2-2r}=\Inn(a^2)(c)=\pi_2^2(c)=(a^{2s}c^t)^{a}=a^{2s}(a^{2s}c^t)^t
=c^{t^2}a^{2sr+2s\sum_{l=1}^{t} r^{l}},$$
that is
\begin{eqnarray}\label{f3.4}
t^2-1\equiv 0(\mod m) \quad {\rm and}\quad s\sum_{l=1}^{t}r^{l}+rs+r-1\equiv 0(\mod n).
\end{eqnarray}

(3) $((\lg a^2\rg\rtimes \lg c\rg).\lg a\rg)\rtimes \lg b\rg $, where $c^b=a^uc^v$.
Set $\pi_3\in \Aut((\lg a^2\rg\rtimes \lg c\rg).\lg a\rg) :$ $a\to a^{-1}$ and $c\to a^uc^v$.
We divide the proof into two cases according to the parity of $u$, separately.
\vskip 3mm
{\it Case 1: $u$ is even.}
\vskip 3mm

(i)  $\pi_3$ preserves  $(a^2)^c=a^{2r}$:
\begin{eqnarray}\label{f3.5}
r^{v-1}-1\equiv 0(\mod n).
\end{eqnarray}

(ii) $\o(\pi_3(c))=m$:
$$1=(a^uc^v)^m=c^{vm}a^{u\sum_{l=1}^{m}r^l},$$
that is
\begin{eqnarray}\label{f3.6}
u\sum_{l=1}^{m}r^l\equiv0(\mod 2n).
\end{eqnarray}

(iii) $\pi_3$ preserves $c^a=a^{2s}c^t$:  that is
$$a^uc^v=(a^{-2s}(a^uc^v)^t)^a=a^{-2s}(a^u(a^{2s}c^t)^v)^t
=a^{-2s}c^{t^2v}a^{(ru+2s\sum_{l=1}^{v}r^l)\sum_{l=0}^{t-1}r^l},$$
which implies
\begin{eqnarray}\label{f3.7}
r(u+2s)\equiv (ru+2s\sum_{l=1}^{v}r^l)\sum_{l=0}^{t-1}r^l (\mod 2n).
\end{eqnarray}
By  Eq(\ref{f3.4}),  Eq(\ref{f3.7}) is if and only if
\begin{eqnarray}\label{f3.8}
u\sum_{l=1}^{t}r^l-ur-2s\sum_{l=1}^{v}r^l-2sr+2(1-r)\equiv 0 (\mod 2n).
\end{eqnarray}

(iv) $\pi_3^2=\Inn(a^n)$:
If $n$ is even, then $a^n\in\lg a^2\rg$, which implies $(a^n)^c=a^n$.
If $n$ is odd, then $c^{a^n}=(c^{a^{n-1}})^a=c^ta^{(n-1)(1-r)+2sr},$
where $t^2\equiv1(\mod m),\,((n-1)(1-r)+2sr)(1+\sum_{l=0}^{t-1}r^l)\equiv0(\mod 2n)$.
Then there are two subcases:
\vskip 3mm
{\it Subcase 1.1: $n$ is even.}
\vskip 3mm
In this case, $\pi_3^2=\Inn(a^n)=1$.
Suppose that $v=1$. Then $c=c^{b^2}=a^{-u}a^uc=c,$ as desired.
Suppose that $v\neq 1$. Then
$$c=c^{b^2}=a^{-u}(a^uc^v)^v=c^v(a^uc^v)^{v-1}=c^{v^2}a^{u\sum_{l=1}^{v-1}r^l},$$
that is
\begin{eqnarray}
u\sum_{l=1}^{v-1}r^l\equiv0(\mod 2n)\quad{\rm{and}}\quad v^2-1\equiv 0(\mod m).
\end{eqnarray}
Then $\pi_3^2=1$ holds if and only if
\begin{eqnarray}\label{f3.9}
u\sum_{l=1}^vr^l-ur\equiv0(\mod 2n)\quad{\rm{and}}\quad v^2-1\equiv 0(\mod m).
\end{eqnarray}
\vskip 3mm
{\it Subcase 1.2: $n$ is odd.}
\vskip 3mm
In this case, $c^{a^n}=(c^{a^{n-1}})^a=c^ta^{(n-1)(1-r)+2sr}.$

Suppose   that $v=1$.
Then $c^ta^{(n-1)(1-r)+2sr}=c^{b^2}=a^{-u}a^uc=c,$
which implies $t=1$. So $X=G\rtimes C$.

Suppose   that $v\neq 1$. Then
$$c^ta^{(n-1)(1-r)+2sr}=c^{b^2}=a^{-u}(a^uc^v)^v=c^v(a^uc^v)^{v-1}=c^{v^2}a^{u\sum_{l=1}^{v-1}r^l},$$
that is
\begin{eqnarray}
u\sum_{l=1}^{v-1}r^l\equiv(n-1)(1-r)+2sr(\mod 2n)\quad{\rm{and}}\quad v^2\equiv t(\mod m).
\end{eqnarray}
Then $\pi_3^2=\Inn(a^n)$ holds if and only if
\begin{eqnarray}\label{f3.9.1}
u\sum_{l=1}^vr^l-ur\equiv(n-1)(1-r)+2sr(\mod 2n)\quad {\rm{and}}\quad v^2-t\equiv 0(\mod m).
\end{eqnarray}
\vskip 3mm
{\it Case 2: $u$ is odd.}
\vskip 3mm
If $t=1$, then $c$ normalises $\lg a\rg$, which implies $\lg a\rg\lhd X$.
By Lemma~\ref{a}, we get $G\lhd X$.
Then $v=1$.
So assume $t\ne 1$ and we shall get a contradiction.

Let $S=\lg a^2, c\rg $. Since $u$ is odd again, we know that $\lg a^2\rg \le S_X< S$.
Since $|X:S|=4$, we have $\ox=X/S_X=\lg \olc, \ola\rg.\lg \olb\rg \lessapprox S_4$. The only possibility is
$\o(\olc)=2$ and $\olx \cong D_8$ so that $m$ is even and $v$ is odd.
Then $t$ is odd, as $t^2\equiv 1{\pmod m}$.
Moreover, we have $\lg a^2, c^2\rg=S_X\lhd X$.

Consider $\ox=X/\lg a^2\rg=\lg\ola,\olc\rg\rtimes\lg\olb\rg$,
where $\ola^{\olb}=\ola,\,\olb^2=1,\,\olc^{\ola}=\olc^t$ and $\olc^{\olb}=\ola\olc^v$.
Let  $\pi_3$ be defined as above.
Since the induced action of $\pi_3$ preserves  $\olc^{\ola}=\olc^t$, we have
$(\ola\olc^v)^{\ola}=(\ola\olc^v)^t$, that is
$$\ola\olc^{tv}=\ola\olc^v((\ola\olc^v)^2)^{\frac{t-1}2}=\ola\olc^v(\olc^{tv+v})^{\frac{t-1}2}
=\ola\olc^{v+\frac{v(t+1)(t-1)}2},$$
which implies
$$tv\equiv v+\frac{v(t+1)(t-1)}2(\mod m).$$
Noting  $t^2\equiv1{\pmod m}$, $t\ne 1$  and $(v, m)=1$ is odd, we get
\begin{eqnarray}\label{f3.10}
t\equiv 1+\frac m2(\mod m).
\end{eqnarray}

Let $X_1=GS_X=\lg a, b\rg \lg c^2\rg $. By Eq(\ref{f3.10}), we have    $(c^2)^a=(a^{2s}c^t)^2=a^{2s(1+r^{-1})}c^2$, which implies $c^2$ that normalises $\lg a\rg$.
By Lemma~\ref{a}, we get $G\lhd X_1$.

If $n$ is odd, then $X_1=\lg a,b\rg\rtimes\lg c^2\rg\lhd X$, which implies $\lg a^n\rg\lhd X_1$.
Note that $a^n$ is an involution and $\lg c^2\rg_{X_1}=1$, then $Z(X_1)=\lg a^n\rg$.
Then $\lg a^n\rg\char X_1\lhd X$, which implies $a^n\in G_X$, that is $\lg a\rg\lhd X$.
By Lemma~\ref{a} again, we get $G\lhd X$, which implies $t=v=1$, a contradiction.
So in what follows, we assume that  $n$ is even.

By $G\lhd X_1$, we get $b^{c^2}\in G$.
Since $b^{c^2}=c^{-2}a^n(b^{-1}c^2b)b=c^{-2}a^n(a^uc^v)^2b=c^{v(t+1)-2}a^xb,$
for some $x$, we get  $v(t+1)-2\equiv0 (\mod m)$.
By combing Eq(\ref{f3.10})  we get
\begin{eqnarray}\label{f3.11}
v\equiv 1\pm \frac m4 (\mod \frac m2), \quad 4\di m.
\end{eqnarray}

Since  $$\olc=\olc^{\olb^2}=\ola(\ola\olc^v)^v=\olc^v(\ola\olc^v\ola\olc^v)^{\frac{v-1}2}
=\olc^{v+(tv+v)\frac{v-1}2},$$
\f that is
\begin{eqnarray}\label{f3.12}
(v-1)(\frac{v(t+1)}2+1) \equiv 0(\mod m).
\end{eqnarray}
Then  Eq(\ref{f3.11}) and  Eq(\ref{f3.12}) may give $\frac m2\equiv 0(\mod m)$, a contradiction.
\vskip 3mm
(4) Insure $\lg c\rg_X=1$:
If $t=1$, then $v=1$ and $1-2sr\equiv r(\mod n)$  by Eq(\ref{f3.4}).
For any integer $w$,
$$(c^w)^a=(a^{2s}c)^w=c^wa^{2s\sum_{l=1}^{w}r^l}\quad{\rm{and}}\quad
(c^w)^b=(a^uc)^w=c^{w}a^{u\sum_{l=1}^{w}(1-2sr)^l}.$$
Since $\lg c\rg_X=1$,  we know that
$2s\sum_{l=1}^{w}r^l\equiv0\equiv u\sum_{l=1}^{w}(1-2sr)^l(\mod 2n)\Leftrightarrow w\equiv0(\mod m)$.

If $t\neq1$, then $u$ is even, for any integer $w$,
$$(c^w)^a=(a^{2s}c^t)^w=c^{tw}a^{2s\sum_{l=1}^{w}r^l}\quad{\rm{and}}\quad (c^w)^b=(a^uc^v)^w=c^{vw}a^{u\sum_{l=1}^{w}r^l}.$$
Since $\lg c\rg_X=1$,  we know that
$2s\sum_{l=1}^{w}r^l\equiv  0\equiv u\sum_{l=1}^{w}r^l(\mod 2n)\Leftrightarrow w\equiv0(\mod m)$.

Summarizing  Eq(\ref{f3.1})-Eq(\ref{f3.9.1}), we get  the parameters $(m,n, r, s, t, u, v)$ as shown in the lemma. Moreover, since we  do every  above  group extension  by sufficient and necessary conditions,     for any given parameters satisfying the equations, there exists  $X=X(Q)$.
\qed
\subsection{$M=\lg a^2\rg\lg c\rg $ and $X/M_X\cong D_8$}
\begin{lem}\label{QD8}
Suppose  that $X=X(Q)$, $M=\lg a^2\rg\lg c\rg$, $X/M_X\cong D_8$ and $\lg c\rg_X=1$. Then
$$X=\lg a,b,c|R,(a^2)^{c^2}=a^{2r},(c^2)^a=a^{2s}c^{2t},(c^2)^b=a^{2u}c^{2}, a^c=bc^{2w}\rg,$$
where either $w=0$ and $r=s=t=u=1$; or
  $$\begin{array}{ll}
  &w\neq0,\,s=u^2\sum_{l=0}^{w-1}r^l+\frac{un}2,\,t=2wu+1,\\
  &r^{2w}-1\equiv (u\sum_{l=1}^{w}r^l+\frac n2)^2-r\equiv0(\mod n),\\
  &s\sum_{l=1}^tr^{l}+sr\equiv 2sr-u\sum_{l=1}^tr^l+ur\equiv 1-r(\mod n),\\
  &2w(1+uw)\equiv nw\equiv 2w(r-1)\equiv0(\mod\frac m2),\\
  &2^{\frac{1+(-1)^u}2}\sum_{l=1}^{i}r^l\equiv 0(\mod n)\Leftrightarrow i\equiv0(\mod\frac m2).
   \end{array}$$
\end{lem}
\demo  Under the hypothesis,   $M_X=\lg a^2\rg \rtimes\lg c^2\rg$. Set $2n=\o(a)$ and $m=\o(c)$.
If $n$ is odd, then $\lg\ola,\olb\rg\cong\ZZ_4$, a contradiction. So  both $n$ and $m$ are even.
Since $X/M_X=\lg\ola,\olb\rg\lg\olc\rg\cong D_8$, we can choose $\olb$
such that  the form of $X/M_X$ is the following:
$\ola^{\olc}=\olb$ and $\olb^{\olc}=\ola.$
Set $c_1:=c^2$ and $X_1=GM_X=\lg a,b\rg\lg c_1\rg$.
Noting $\lg a\rg\lg c_1\rg\leq X_1$ and $\lg c_1\rg_{X_1}=1$, by  Lemma~\ref{Qac}, we get
$$X_1=\lg a,b,c_1|R, (a^2)^{c_1}=a^{2r},c_1^a=a^{2s}c_1^{t}, c_1^b=a^{2u}c_1^v\rg$$
whose
\begin{eqnarray}\label{f4.1}
  \begin{array}{ll}
  &r^{t-1}\equiv r^{v-1}\equiv1(\mod n),\, t^2\equiv v^2\equiv1(\mod\frac m2),\\
  &2s\sum_{l=1}^tr^{l}+2sr\equiv 2sr+2s\sum_{l=1}^{v}r^l-u\sum_{l=1}^tr^l+ur
  \equiv 2(1-r)(\mod 2n),\\
  &u(\sum_{l=0}^{v-1}r^l-1)\equiv 0(\mod 2n),\\
  &s\sum_{l=1}^{i}r^l\equiv u\sum_{l=1}^{i}r^l\equiv 0(\mod n)
  \Leftrightarrow i\equiv 0(\mod \frac m2).
  \end{array}
\end{eqnarray}
Moreover, since $n$ is even and $\lg c_1\rg_{X_1}=1$, we get that $b^2=a^n$ is the unique involution of $Z(X_1)$. Then $a^n\in Z(X)$, that is $[b^2,c]=[a^n,c]=1$.
Now $X=X_1.\lg c\rg$. Set $a^c=bc_1^w$. Then $X$ may be defined by $R$ and
\begin{eqnarray}\label{m2}
(a^2)^{c_1}=a^{2r},\,c_1^a=a^{2s}c^{2t},\,c_1^b=a^{2u}c^{2v},\,a^c=bc_1^w.
\end{eqnarray}

If $w\equiv0(\mod\frac m2)$, then $\o(a)=\o(a^c)=\o(b)=4$, which implies $G\cong Q_8$,
and one can check $X$ is isomorphic to the following form:
$$X=\lg a,b,c|a^4=c^4=1,b^2=a^2,a^b=a^{-1},a^c=b,b^c=a^{-1}\rg,$$
that is the former part of Lemma~\ref{QD8}.
So in that follows, we assume $w\not\equiv 0(\mod\frac m2)$.

Firstly, we get $b^c=a^{c^2}c^{-2w}=c_1^{w+t-1}a^{2sr-1+n}.$
Set $\pi\in \Aut(X_1) :$
$a\to bc_1^{w}$, $b\to a^{1-2sr}c_1^{1-t-w}$ and $c_1\to c_1$.
We need to carry out the following  seven steps:

(i) $\o(\pi(b))=4:$
Since $b^2\in Z(X)$, we only show $(b^c)^2=a^n$:
$$(c^{2(w+t-1)}a^{2sr-1+n})^2=c^{2w(t+1)}a^{2sr^{w+1}+2s\sum_{l=1}^{w+t-1}r^{l}+2sr-2}=a^n,$$
that is
\begin{eqnarray}\label{f4.2}
w(t+1)\equiv 0(\mod\frac m2)\quad{\rm{and}}\quad
sr^{w+1}+s\sum_{l=1}^{w+t-1}r^{l}+sr-1\equiv \frac n2(\mod n),
\end{eqnarray}
which implies
\begin{eqnarray}\label{f4.3}
r^{2w}\equiv r^{w(t+1)}\equiv1(\mod n).
\end{eqnarray}

(ii)  $\o(\pi(a))=2n$:
$$
(bc^{2w})^{2n}=(c^{2w(v+1)}a^{2u\sum_{l=w+1}^{2w}r^l+n})^n
= c^{2nw(v+1)}a^{2nu\sum_{l=w+1}^{2w}r^l}=1,
$$
that is
\begin{eqnarray}\label{f4.4}
nw(v+1)\equiv0(\mod\frac m2).
\end{eqnarray}

(iii) $\pi$ preserves $(a^2)^{c_1}=a^{2r}$:

$$((a^2)^{c^2})^c=c^{2w(v+1)}a^{2ur\sum_{l=w+1}^{2w}r^l+n}\quad{\rm{and}}\quad
(a^{2r})^c=c^{2wr(v+1)}a^{2ur\sum_{l=w+1}^{2w}r^l+n},$$
that is
\begin{eqnarray}\label{f4.5}
w(v+1)(r-1)\equiv 0(\mod\frac m2).
\end{eqnarray}

(iv) $\pi $ preserves $c_1^a=a^{2s}c_1^{t}$:
$$ \begin{array}{lcl}
 ((c^2)^a)^c&=&(c^2)^{bc^{2w}}=(a^{2u}c^{2v})^{c^{2w}}=a^{2ur^w}c^{2v}=c^{2v}a^{2ur^{w+1}},\\
 (a^{2s}c^{2t})^c&=&(c^{2w(v+1)}a^{2u\sum_{l=w+1}^{2w}r^l+n})^sc^{2t}
 =c^{2ws(v+1)+2t}a^{2sru\sum_{l=w+1}^{2w}r^l+ns},
 \end{array}$$
that is
\begin{eqnarray}\label{f4.6}
v\equiv ws(v+1)+t(\mod\frac m2)\quad{\rm{and}}\quad
u\equiv su\sum_{l=1}^{w}r^l+\frac{ns}2(\mod n).
\end{eqnarray}

(v) $\pi$ preserves $c_1^b=a^{2u}c_1^{v}$:
$$ \begin{array}{lcl}
 ((c^2)^b)^c&=&(c^2)^{a^{1-2sr}c^{2-2t-2w}}=c^{2t}a^{2sr^{2-w}}\\
 (a^{2u}c^{2v})^c&=&(c^{2w(v+1)}a^{2u\sum_{l=w+1}^{2w}r^l+n})^uc^{2v}=c^{2wu(v+1)+2v}
                 a^{2u^2\sum_{l=w+2}^{2w+1}r^l+un},
 \end{array}$$
that is,
\begin{eqnarray}\label{f4.7}
t\equiv wu(v+1)+v(\mod\frac m2)\quad{\rm{and}}\quad
s\equiv u^2\sum_{l=0}^{w-1}r^l+\frac{un}2(\mod n).
\end{eqnarray}

(vi) $\pi^2=\Inn(c_1)$:
recall $\Inn(c_1)(a)=a^{1-2sr}c_1^{1-t},\,\Inn(c_1)(a^2)=a^{2r}$
and $\Inn(c_1)(b)=c_1^{v-1}a^{2ur}b$.
$$a^{1-2sr}c^{2-2t}=\Inn(c_1)(a)=\pi^2(a)=b^cc^{2w}=a^{1-2sr}c^{2-2t-2w+2w},$$
as desired;
$$a^{2r}=\Inn(c_1)(a^2)=\pi^2(a^2)=(c^{2w(v+1)}a^{2u\sum_{l=w+1}^{2w}r^l+n})^c
 =c^{2w(v+1)(1+uw+\frac n2)}a^{2(u\sum_{l=1}^wr^l+\frac n2)^2},$$
that is
\begin{eqnarray}\label{f4.8}
w(v+1)(1+uw+\frac n2)\equiv 0(\mod \frac m2)\quad{\rm{and}}\quad
r\equiv (u\sum_{l=1}^{w}r^l+\frac n2)^2 (\mod n),
\end{eqnarray}
and noting Eq (\ref{f4.5}) and (\ref{f4.6}), we get $w(v+1)(r-1)\equiv ws(v+1)+t-v\equiv0(\mod\frac m2)$ and
$u\equiv su\sum_{l=1}^{w}r^l+\frac{ns}2(\mod n).$
Then
$$c^{2(v-1)}a^{2ur}b=\Inn(c_1)(b)=\pi^2(b)=c^{2(t-1+wsr(v+1))}a^{2usr\sum_{l=1}^{w}r^l+(s+1)n}b,$$
as desired.

(vii) Insure $\lg c\rg_X=1$:
Since $\lg c\rg_X\leq M$, we get $\lg c\rg_X\leq M_X=\lg a^2\rg\lg c^2\rg$.
Then $\lg c\rg_X=\cap_{x\in X}C^x=\cap_{x\in G}C^x=\cap_{x\in G}\lg c^2\rg^x=\lg c^2\rg_{X_1}=1$.
Recall $2s\sum_{l=1}^{i}r^l\equiv u\sum_{l=1}^{i}(1-s(\sum_{l=1}^tr^l+r))^l\equiv 0(\mod 2n)
  \Leftrightarrow i\equiv 0(\mod \frac m2).$

Now we are  ready to determine the parameters by summarizing  Eq(\ref{f4.1})-Eq(\ref{f4.8}).
Firstly, we shall show $v=1$.

Suppose that $u$ is odd. Then by Eq(\ref{f4.8}), we get $(u,n)=(\sum_{l=1}^{w}r^l,\frac n2)=1$ as $(r,n)=1$.
Moreover, if $\frac n2$ is odd, then $\sum_{l=1}^{w}r^l$ is even as $r$ is odd.
Then by Eq(\ref{f4.7}), we get $s\equiv u^2\sum_{l=0}^{w-1}r^l+\frac{n}2(\mod n)$,
which implies $(s,n)=1$.
Then we have $\sum_{l=1}^{i}r^l\equiv0(\mod n)\Leftrightarrow i\equiv0(\mod \frac m2)$ by (vii) and $\sum_{l=0}^{v-1}r^l\equiv1(\mod n)$ from Eq(\ref{f4.1}).
Then $v\equiv 1(\mod \frac m2)$.

Suppose that $u$ is even. Then by Eq(\ref{f4.8}), we get
$(u,\frac n2)=(\sum_{l=1}^{w}r^l,\frac n2)=1$.
Then by Eq(\ref{f4.7}), we get $s\equiv u^2\sum_{l=0}^{w-1}r^l(\mod n)$,
which implies $s$ is even. Then by Eq(\ref{f4.6}), we get $u\equiv su\sum_{l=1}^{w}r^l(\mod n).$
Then we have $\sum_{l=1}^{i}r^l\equiv0(\mod\frac n2)\Leftrightarrow i\equiv0(\mod \frac m2)$ by (vii)
and $\sum_{l=0}^{v-1}r^l\equiv1(\mod\frac n2)$ from Eq(\ref{f4.1}).
Then $v\equiv 1(\mod \frac m2)$.

Inserting $v=1$ in Eq(\ref{f4.1})-Eq(\ref{f4.8}), we get
$s=u^2\sum_{l=0}^{w-1}r^l+\frac n2$ and  $t=2wu+1$ in Eq(\ref{f4.7});
$nw\equiv0(\mod\frac m2)$ in Eq(\ref{f4.2}), (\ref{f4.4}) and (\ref{f4.8}); and
$2w(r-1)\equiv2w(1+uw)\equiv0(\mod\frac m2)$ in Eq(\ref{f4.5}) and (\ref{f4.8}).
All these are summarized in the lemma.
\qed

\subsection{$M=\lg a^2\rg\lg c\rg$ and $X/M_X\cong A_4$}
\begin{lem}\label{QA4}
Suppose  that $X=X(Q)$, $M=\lg a^2\rg\lg c\rg$,   $X/M_X\cong A_4$ and $\lg c\rg_X=1$. Then
$$X=\lg a,b,c|R, (a^2)^{c}=a^{2r},(c^3)^a=a^{2s}c^3,(c^3)^b=a^{2u}c^{3}, a^c=bc^{\frac {im}2},b^c=a^xb\rg,$$
where $n\equiv 2(\mod 4)$ and either $i=s=u=0$ and $r=x=1$; or
$i=1$, $6\di m$, $r^{\frac m2}\equiv-1(\mod n)$ with $\o(r)=m$, $s\equiv \frac{r^{-3}-1}2(\mod\frac n2)$, $u\equiv\frac{r^3-1}{2r^2}(\mod\frac n2)$ and $x\equiv -r+r^2+\frac n2 (\mod n)$.
\end{lem}

\demo
Under the hypothesis, $M_X=\lg a^2\rg\rtimes\lg c^3\rg$. Set $2n=\o(a)$ and $m=\o(c)$.
If $n$ is odd, then $\lg\ola,\olb\rg\cong\ZZ_4$, a contradiction.
So  $n$ is even and $3\di m$.
Since $X/M_X=\lg\ola,\olb\rg\lg\olc\rg\cong A_4$, we can choose $\olb$
such that  the form of $X/M_X$ is the following:
$\ola^{\olc}=\olb$ and $\olb^{\olc}=\ola\olb.$
Set $c_1:=c^3$ and $X_1=GM_X=\lg a,b\rg\lg c_1\rg$.
By  Lemma~\ref{Qac}, we get
$$X_1=\lg a,b,c^3|R, (a^2)^{c_1}=a^{2r},(c_1)^a=a^{2s}c_1^{t}, (c_1)^b=a^{2u}c_1^{v}\rg$$
whose
\begin{eqnarray}\label{f5.1}
\begin{array}{ll}
  &r^{t-1}-1\equiv r^{v-1}-1\equiv u(\sum_{l=0}^{v-1}r^l-1)\equiv0(\mod n),\\
  &s\sum_{l=1}^tr^{l}+sr\equiv sr+s\sum_{l=1}^{v}r^l-u\sum_{l=1}^tr^l+ur
  \equiv 1-r(\mod n),\\
  &t^2-1\equiv v^2-1\equiv0(\mod\frac m3),\\
  &s\sum_{l=1}^{i}r^l\equiv u\sum_{l=1}^{i}r^l\equiv0(\mod n)\Leftrightarrow i\equiv0(\mod\frac m3).
  \end{array}
\end{eqnarray}
Moreover, since $n$ is even and $\lg c_1\rg_{X_1}=1$, we get that $b^2=a^n$ is the unique involution of $Z(X_1)$, that is  $[b^2,c]=[a^n,c]=1$. Now $X=X_1.\lg c\rg$. Set $a^c=bc_1^{w}$.
Then $X$ may be defined by $R$ and
\begin{eqnarray}\label{m3}
(a^2)^{c_1}=a^{2r},(c_1)^a=a^{2s}c_1^{t}, (c_1)^b=a^{2u}c_1^{v},\,a^c=bc_1^{w},b^c=a^{1+2x}bc_1^{y}.
\end{eqnarray}

If $w\equiv0(\mod\frac m3)$, then $\o(a)=\o(a^c)=\o(b)=4$, which implies $G\cong Q_8$,
and one can check $X$ is isomorphic to the following form:
$$X=\lg a,b,c|a^4=c^3=1,b^2=a^2,a^b=a^{-1},a^c=b,b^c=ab\rg,$$
that is the former part of Lemma~\ref{QA4}.
So in that follows, we assume $w\not\equiv 0(\mod\frac m2)$.

What we should do is to  determine the parameters  $r,s,t,u,v,w,x$ and $y$ by analysing  the last  extension  $X_1.\lg c\rg $, where  $a^c=bc_1^{w}$ and $b^c=a^{1+2x}bc_1^{y}$.
Set $\pi\in \Aut(X_1): a\to bc_1^{w}, \, b\to a^{1+2x}bc_1^{y}, \quad c_1\to c_1.$
We need to carry out the following  seven steps:

(i) $\o(\pi(b))=4$:
Since $b^2\in Z(X)$, we only show $(b^c)^2=a^n$:
$$(a^{1+2x}bc^{3y})^2=c^{3vty+3y}a^{2ur^y\sum_{l=1}^{ty}r^l+2xr^y(r^y-1)-2sr^y\sum_{l=1}^{y}r^l+n}=a^n,$$
that is
\begin{eqnarray}\label{f5.2}
y(tv+1)\equiv0(\mod\frac m3)\quad{\rm{and}}\quad
u\sum_{l=1}^{ty}r^l+xr^y-x-s\sum_{l=1}^yr^l\equiv0(\mod n),
\end{eqnarray}
which implies
\begin{eqnarray}
r^{2y}\equiv r^{y(tv+1)}\equiv1(\mod n).
\end{eqnarray}

(ii) $\o(\pi(ab))=4$: Since $(ab)^2=a^n\in Z(X)$, we only show $((ab)^c)^2=a^n$:
$$\begin{array}{lcl}
a^n&=&((ab)^c)^2=(c^{3vw+3yt}a^{2ur^y\sum_{l=1}^wr^l-2r^yx+2s\sum_{l=1}^yr^l-1+n})^2\\
&=&c^{3(vw+yt+tvw+y)}a^{2((r^{w+y}+1)(ur^y\sum_{l=1}^wr^l-r^yx+s\sum_{l=1}^yr^l)+
s\sum_{l=1}^{vw+yt}r^l-1)},
\end{array}$$
that is
\begin{eqnarray}\label{f5.3}
\begin{array}{ll}
&vw+yt+tvw+y\equiv0(\mod\frac m3);\\
&(r^{w+y}+1)(ur^y\sum_{l=1}^wr^l-r^yx+s\sum_{l=1}^yr^l)+
 s\sum_{l=1}^{vw+yt}r^l-1\equiv\frac n2(\mod n),
\end{array}
\end{eqnarray}
which implies $r^{2w}\equiv1(\mod\frac n2)$.

(iii) $\o(\pi(a))=2n$:
$$(bc^{3w})^n=(bc^{3w}bc^{3w})^{n}=(c^{3w(v+1)}a^{2u\sum_{l=w+1}^{2w}r^l+n})^{n}
=c^{3nw(v+1)}a^{2nur^w\sum_{l=1}^{w}r^l}=1,$$
that is
\begin{eqnarray}\label{f5.4}
nw(v+1)\equiv0(\mod\frac m3).
\end{eqnarray}

(iv) $\pi$ preserves $(a^2)^{c^3}=a^{2r}$:
$$\begin{array}{ll}
&((a^2)^{c^3})^c=(c^{3w(v+1)}a^{2u\sum_{l=w+1}^{2w}r^l+n})^{c^3}=c^{3w(v+1)}a^{2ur\sum_{l=w+1}^{2w}r^l+n},\\
&(a^{2r})^c=(c^{3w(v+1)}a^{2u\sum_{l=w+1}^{2w}r^l+n})^r=c^{3wr(v+1)}a^{2ur\sum_{l=w+1}^{2w}r^l+n},
\end{array}$$
that is
\begin{eqnarray}\label{f5.5}
w(v+1)(r-1)\equiv0(\mod\frac m3).
\end{eqnarray}

(v) $\pi$ preserves $(c^3)^a=a^{2s}c^{3t}$:
$$\begin{array}{ll}
&((c^3)^a)^c=(c^3)^{bc^{3w}}=(a^{2u}c^{3v})^{c^{3w}}=c^{3v}a^{2ur^{w+1}},\\
&(a^{2s}c^{3t})^c=(c^{3w(v+1)}a^{2u\sum_{l=w+1}^{2w}r^l+n})^sc^{3t}
=c^{3ws(v+1)+3v}a^{2sur^{w+1}\sum_{l=1}^{w}r^l+sn},
\end{array}$$
that is
\begin{eqnarray}\label{f5.6}
v\equiv ws(v+1)+t(\mod\frac m3)\quad{\rm{and}}\quad
u\equiv su\sum_{l=1}^{w}r^l+\frac{sn}2(\mod n).
\end{eqnarray}

(vi) $\pi$ preserves $(c^3)^b=a^{2u}c^{3v}$:
$$\begin{array}{ll}
&((c^3)^b)^c=(c^3)^{a^{1+2x}b}=(c^3a^{2x(1-r)})^{ab}=c^{3tv}a^{2(u\sum_{l=1}^tr^l-sr+x(r-1))},\\
&(a^{2u}c^{3v})^c=(c^{3w(v+1)}a^{2u\sum_{l=w+1}^{2w}r^l+n})^uc^{3v}
=c^{3wu(v+1)+3v}a^{2u^2r\sum_{l=w+1}^{2w}r^l+un},
\end{array}$$
that is
\begin{eqnarray}\label{f5.7}
\begin{array}{ll}
&tv\equiv wu(v+1)+v(\mod\frac m3),\\
&u\sum_{l=1}^tr^l-sr+x(r-1)\equiv u^2r^{w+1}\sum_{l=1}^{w}r^l+\frac{un}2(\mod n).
\end{array}
\end{eqnarray}

(vii) $\pi^3=\Inn(c_1)$:
Recall $\Inn(c_1)(a)=a^{1-2sr}c^{3-3t},\,\Inn(c_1)(a^2)=a^{2r}$
and $\Inn(c_1)(b)=c^{3(v-1)}a^{2ur}b$.
$$\begin{array}{ll}
&a^{1-2sr}c^{3-3t}=\Inn(c_1)(a)=\pi^3(a)=(a^{1+2x}bc^{3(w+y)})^c\\
&=a^{n-1}c^{3vtw(1+x(v+1))+3(w+2y)}a^{2ur^{w}\sum_{l=1}^{tw(1+x(v+1))}r^l
-r^{w}(2s\sum_{l=1}^{w(1+x(v+1))}r^l+2ux\sum_{l=w+1}^{2w}r^l+xn+2x)},
\end{array}$$
that is
\begin{eqnarray}\label{f5.8}\begin{array}{ll}
&1-t\equiv t(wv+wxv+wx)+w+2y(\mod\frac m3),\\
&r^{w}(u\sum_{l=1}^{tw(1+x(v+1))}r^l-s\sum_{l=1}^{w(1+x(v+1))}r^l-ux\sum_{l=w+1}^{2w}r^l-\frac {(x+1)n}2-x)\\
&\equiv 1-sr(\mod n);
\end{array}\end{eqnarray}
$$\begin{array}{lcl}
a^{2r}&=&\Inn(c_1)(a^2)=\pi^3(a^2)=(c^{3w(v+1)}a^{2ur^w\sum_{l=1}^wr^l+n})^{c^2}\\
&=&c^{3w(v+1)+3uw^2(v+1)(1+uw)}a^{2r^w(u\sum_{l=1}^{w}r^l)^3}a^{n},
\end{array}$$
that is
\begin{eqnarray}\label{f5.9}
w(v+1)+uw^2(v+1)(uw+1)\equiv0(\mod\frac m3)\quad{\rm{and}}\quad
r\equiv r^w(u\sum_{l=1}^wr^l)^3+\frac n2(\mod n),
\end{eqnarray}
which implies $(u,\frac n2)=(\sum_{l=1}^wr^l,\frac n2)=1$ as $(r,n)=1$,
and moreover, if $u$ is even, then $\frac n2$ is odd as $r$ is odd.
Noting Eq(\ref{f5.6}), that is $u\equiv su\sum_{l=1}^{w}r^l+\frac{sn}2(\mod n)$, we get that
$(s,\frac n2)=1$ and
$$\begin{array}{lcl}
c^{3(v-1)}a^{2ur}b&=&\Inn(c_1)(b)=\pi^3(b)=(c^{3(w+t-1)}a^{2sr-1+n})^c\\
&=&c^{3(t-1+ws(v+1))}a^{2usr\sum_{l=1}^{w}r^l+sn}b,
\end{array}$$
that is
\begin{eqnarray}\label{f5.10}
v\equiv t+ws(v+1)(\mod\frac m3).
\end{eqnarray}

(viii) Insure $\lg c\rg_X=1$:
Since $\lg c\rg_X\leq M$, we get $\lg c\rg_X\leq M_X=\lg a^2\rg\lg c^3\rg$.
Since $\lg c\rg_X=\cap_{x\in X}C^x=\cap_{x\in G}C^x=\cap_{x\in G}\lg c^3\rg^x=\lg c^3\rg_{X_1}=1$ and
$s\sum_{l=1}^{i}r^l\equiv u\sum_{l=1}^{i}r^l\equiv0(\mod n)\Leftrightarrow i\equiv0(\mod\frac m3).$
Noting that $(s,\frac n2)=(u,\frac n2)=1$ and both $u$ and $s$ are even only if $\frac n2$ is odd,
we have $2^{\frac{1+(-1)^u}2}\sum_{l=1}^{i}r^l\equiv 0(\mod  n)
\Leftrightarrow i\equiv0(\mod\frac m3)$.
\vskip 3mm
Now we are  ready to determine the parameters by summarizing  Eq(\ref{f5.1})-Eq(\ref{f5.10}).
Then we shall divide it into three steps:
\vskip 3mm
{\it Step 1: $t=v=1$, $w=\frac m6$, $r^w\equiv-1(\mod n)$ and $s\equiv \frac{1-r}{2r}(\mod\frac n2)$.}
\vskip 3mm
Since $(r,n)=(u,\frac n2)=1$ (after Eq(\ref{f5.9})),  we get from   Eq(\ref{f5.1}) that
$2^{\frac{1+(-1)^u}2}\sum_{l=1}^{v-1}r^l\equiv 0(\mod n).$
By (viii), $2^{\frac{1+(-1)^u}2}\sum_{l=1}^{i}r^l\equiv 0(\mod n)
\Leftrightarrow i\equiv0(\mod\frac m3)$,
which means  $v\equiv 1(\mod \frac m3)$.

Inserting $v=1$ in Eq(\ref{f5.1})-Eq(\ref{f5.10}), we get that
$2w(wu+1)\equiv0(\mod\frac m3)$ and $t\equiv1+2wu\equiv1-2ws(\mod\frac m3)$
in Eq(\ref{f5.3}), (\ref{f5.6}) and (\ref{f5.7}).
Then $2w\equiv0(\mod\frac m3)$ by Eq(\ref{f5.9}),
which implies $w=\frac m6$ as $w\not\equiv 0(\mod\frac m3)$.
Inserting $w=\frac m6$ in Eq(\ref{f5.1})-Eq(\ref{f5.10}) again, we get
$t\equiv1(\mod\frac m3)$ in Eq(\ref{f5.7}),
$s\equiv \frac{1-r}{2r}(\mod\frac n2)$ in Eq(\ref{f5.1})
and $r^w\equiv-1(\mod n)$ in Eq(\ref{f5.1}) and (\ref{f5.9}).
\vskip 3mm
{\it Step 2: $y=0$}
\vskip 3mm
Since $2y\equiv0(\mod\frac m3)$ in Eq(\ref{f5.2}), we know that $y$ is either $0$ or $\frac m6$.

Suppose that $y=\frac m6=w$. Then by Eq(\ref{f5.3}) we get that $\frac n2$ is odd,
and with Eq(\ref{f5.2}) and (\ref{f5.3}), we get
$2x\equiv (u-s)\sum_{l=1}^{w}r^l\equiv\frac n2-1(\mod n)$.
By Eq(\ref{f5.6}) and $(u,\frac n2)=1$, we get $2s\sum_{l=1}^{w}r^l\equiv1+\frac {sn}2(\mod n)$,
then $u\sum_{l=1}^{w}r^l\equiv0(\mod \frac n2)$, contradict to
$(u,\frac n2)=(\sum_{l=1}^{w}r^l,\frac n2)=1$.
So $y=0$.
\vskip 3mm
{\it Step3: Determine $u$ and $x$.}
\vskip 3mm
By Eq(\ref{f5.3}), we get $s\sum_{l=1}^{w}r^l\equiv1+\frac n2(\mod n).$
Then $\frac{(s+u)n}2\equiv0(\mod n)$ in Eq(\ref{f5.6}), which implies $u\equiv s(\mod2)$.
By Eq(\ref{f5.7}), we get $x(r-1)\equiv (s-u)r-u^2r\sum_{l=1}^{w}r^l+\frac{un}2(\mod n)$.
If $\frac n2$ is even, then $u\sum_{l=1}^{w}r^l$ is even.
But in Eq(\ref{f5.9}), we get $r\equiv -(u\sum_{l=1}^wr^l)^3+\frac n2(\mod n)$
which implies that $u\sum_{l=1}^{w}r^l$ is odd as both $r$ and $\frac n2$ are odd,
a contradiction. So $\frac n2$ is odd. Then we get $u\sum_{l=1}^wr^l$ is even in Eq(\ref{f5.9})
and $u$ is even in Eq(\ref{f5.7}). Then $\sum_{l=1}^{i}r^l\equiv 0(\mod \frac n2)
\Leftrightarrow i\equiv0(\mod\frac m3)$.
Recall $s\equiv \frac{r^{-1}-1}2(\mod\frac n2)$ in Eq(\ref{f5.1})
and $s\sum_{l=1}^{w}r^l\equiv1+\frac n2(\mod n)$ in Eq(\ref{f5.3}).
Since $(\sum_{l=1}^{w}r^l,\frac n2)=1$ and $x(r-1)\equiv (s-u)r-u^2r\sum_{l=1}^{w}r^l(\mod n)$,
we get $2x\equiv u\sum_{l=1}^{w}r^l+(u\sum_{l=1}^{w}r^l)^2-1+\frac n2(\mod n)$.
And by Eq(\ref{f5.9}), we get $-r\equiv (u\sum_{l=1}^{w}r^l)^3+\frac n2(\mod n).$
Take $l=-u\sum_{l=1}^{w}r^l+\frac n2$, then $l\equiv-\frac{2ru}{1-r}(\mod\frac n2)$,
$r\equiv l^3(\mod n)$, $u\equiv\frac{l^3-1}{2l^2}(\mod\frac n2)$ and
$1+2x\equiv -l+l^2+\frac n2 (\mod n)$.
Let us re-write $l$ as $r$ and $1+2x$ as $x$ for the sake of formatting.
Then $s\equiv \frac{r^{-3}-1}2(\mod\frac n2)$, $u\equiv\frac{r^3-1}{2r^2}(\mod\frac n2)$
and $x\equiv -r+r^2+\frac n2 (\mod n)$.
\qed
\vskip 3mm
In fact, if we add the conditions $t=1$ and $w\ne0$ and delete $\lg c\rg_X=1$ in the above calculation,
then we can get the following:
\begin{lem}\label{A4.1}
With the notation, suppose   that $t=1$ and $w\ne0$. Then
$$X=\lg a,b,c|R,
(a^2)^{c}=a^{2r},(c^3)^a=a^{2s}c^3,(c^3)^b=a^{2u}c^{3},a^c=bc^{\frac m2},b^c=a^xb\rg,$$
where $n\equiv 2(\mod 4)$, $m\equiv 0(\mod 6)$, $r^{\frac m2}\equiv-1(\mod n)$,
$s\equiv \frac{r^{-3}-1}2(\mod\frac n2)$, $u\equiv\frac{r^3-1}{2r^2}(\mod\frac n2)$ and $x\equiv -r+r^2+\frac n2(\mod n)$.
\end{lem}

\subsection{$M=\lg a^4\rg\lg c\rg$ and $X/M_X\cong S_4$}
\vskip 3mm
\begin{lem}\label{Qa4c3}
Suppose  that $X=X(Q)$, $M=\lg a^4\rg\lg c\rg$, $X/M_X\cong S_4$ and $\lg c\rg_X=1$. Then
$X=\lg a,b,c|R, (a^4)^{c}=a^{4r}, c_1^{a^2}=a^{4s}c_1,c_1^b=a^{4u}c_1,
(a^2)^c=bc^{\frac {im}2},b^c=a^{2x}b,c^a=a^{2(1+2z)}c^{1+\frac{jm}3}\rg,$
where either
\begin{enumerate}
  \item[\rm(1)] $i=0$, $r=j=1,x=3,s=u=z=0$; or
  \item[\rm(2)] $i=1$, $n\equiv 4(\mod 8)$, $6\di m$,
  $r^{\frac m2}\equiv-1(\mod \frac n2)$, $\o(r)=m$,
  $s\equiv \frac{r^{-3}-1}2(\mod\frac n4)$, \\
  $u\equiv\frac{r^3-1}{2r^2}(\mod\frac n4)$, $x\equiv -r+r^2+\frac n4(\mod \frac n2)$,
  $1+2z\equiv \frac{1-r}{2r}(\mod\frac n2)$, $j\in\{1,2\}$.
\end{enumerate}
\end{lem}
\demo
Under the hypothesis, $M_X=\lg a^4\rg\rtimes\lg c^3\rg$. Set $2n=\o(a)$ and $m=\o(c)$.
Then $n$ is even and $3\di m$.
If $\frac n2$ is odd, then $\lg\ola,\olb\rg\cong Q_8$, a contradiction.
So $\frac n2$ is even.
Since $X/M_X=\lg\ola,\olb\rg\lg\olc\rg\cong S_4$, we can choose $\olb$
such that  the form of $X/M_X$ is the following:
$(\ola^2)^{\olc}=\olb,\,\olb^{\olc}=\ola^2\olb$ and $(\olc)^{\ola}=\ola^2\olc^2.$
Take $a_1=a^2$ and $c_1=c^3$.
Then we set $a_1^c=bc_1^w, b^c=a_1^xbc_1^{y},c^a=a_1^{1+2z}c^{2+3d}$,
where $x$ is odd.

Suppose $w\equiv0(\mod\frac m3)$.
Note that $\o(a_1)=\o(a_1^c)=\o(b)=4$, which implies $G\cong Q_{16}$.
Thus one can check $X$ can only have the following form:
$X=\lg a,b,c|a^8=c^3=1,b^2=a^4,a^b=a^{-1},b^c=a_1^3b,c^a=a_1c^2\rg.$
So in what follows, we assume $w\not\equiv0(\mod \frac m3)$.

Then consider $X_1=GM_X=\lg a,b\rg\lg c_1\rg$.
Noting $\lg a\rg\lg c_1\rg\leq X_1$ and $\lg c_1\rg_{X_1}=1$, by Lemma~\ref{Qa^2},
we know $\lg a_1\rg\lhd X_1$, which implies that $c_1$ normalises $\lg a_1\rg$.
Take $X_2=\lg a_1,b\rg\lg c\rg$.
Then we get $X_2=(\lg a_1,b\rg\lg c_1\rg).\lg c\rg$.
Note that $c_1$ normalises $\lg a_1\rg$ in $X_2$.  Thus by Lemma~\ref{A4.1}, we get
$$X_2=\lg a_1,b,c|R,
(a_1^2)^{c}=a_1^{2r},c_1^{a_1}=a_1^{2s}c_1,c_1^b=a_1^{2u}c_1,a_1^c=bc^{\frac m2},b^c=a^xb\rg,$$
where
\begin{eqnarray}\label{f6.1}
\begin{array}{ll}
&n\equiv 4(\mod 8), m\equiv 0(\mod 6)\\
&r^{\frac m2}\equiv-1(\mod \frac n2), s\equiv \frac{r^{-3}-1}2(\mod\frac n4),
u\equiv\frac{r^3-1}{2r^2}(\mod\frac n4), x\equiv -r+r^2+\frac n4(\mod \frac n2).
\end{array}
\end{eqnarray}

Note $X=X_2.\lg a\rg$.
So $X$ may be defined by $R$ and
\begin{eqnarray}\label{m4}
(a_1^2)^{c}=a_1^{2r}, c_1^{a_1}=a_1^{2s}c_1,c_1^b=a_1^{2u}c_1,
a_1^c=bc^{\frac m2},b^c=a_1^xb,c^a=a_1^{1+2z}c^{2+3d}.
\end{eqnarray}
What we should to  determine the parameters  $r,z$ and $d$ by analyse the last one extension.
\vskip 3mm
$X_2.\lg a\rg $, where $c^a=a_1^{1+2z}c^{2+3d}$.
Set $\pi\in \Aut(X_1) :$
$a_1\to a_1$, $b\to a_1^{-1}b$ and $c\to a_1^{1+2z}c^{2+3d}$, where $d$ is odd.
We need to check the following eight equalities:

(i) $\pi$ preserves $(a_1^2)^{c}=a_1^{2r}$:
$$a_1^{2r}=((a_1^2)^{c})^a=(a_1^2)^{c^{2+3d}}=a_1^{2r^{2+3d}},$$
that is
\begin{eqnarray}\label{f6.2}
r^{1+3d}-1\equiv0(\mod\frac n2).
\end{eqnarray}
Since $r^{\frac m2}\equiv-1(\mod \frac n2)$, we get $\sum_{l=1}^{1+3d}r^{3l}\equiv0(\mod\frac n2).$

(ii) $\pi$ preserves $c_1^{a_1}=a_1^{2s}c_1$, that is $a_1^{c_1}=a_1^{1-2sr^3}$:
$$a_1^{1-2sr^3}=(a_1^{c_1})^a=a_1^{c_1^{2+3d}}=a_1^{(1-2sr^3)^{2+3d}},$$
that is
\begin{eqnarray}\label{f6.3}
(1-2sr^3)^{1+3d}-1\equiv0(\mod n).
\end{eqnarray}

(iii) $\pi$ preserves $c_1^b=a_1^{2u}c_1$, that is $b^{c_1}=a_1^{2ur^3}b$:
$$(b^{c_1})^a=a_1^{2ur^3-3+6sr^3+r^4-2r^2-3r-4z(r+r^2+r^3)+\frac n2}b\quad{\rm{and}}\quad
(a_1^{2ur^3}b)^a=a_1^{2ur^3-1}b,$$
that is
\begin{eqnarray}\label{f6.4}
6sr^3+r^4-2r^2-3r-4z(r+r^2+r^3)+\frac n2\equiv 2(\mod n).
\end{eqnarray}

(iv) $\pi$ preserves $a_1^c=bc^{\frac m2}$:
$$\begin{array}{ll}
&(a_1^c)^a=a_1^{c^{2+3d}}=(a_1^xbc^{\frac m2})^{c_1^d}=a_1^{x(1-2sr^3)^d+2u\sum_{l=1}^dr^{3l}}bc^{\frac m2},\\
&(bc^{\frac m2})^a=a_1^{-1}b(c^{\frac m2})^a=
a_1^{(2z(r+r^2+r^3)+\frac n2-2sr^3+r^2-x(1-2sr^3)^{d+\frac m6}-2ur^2\sum_{l=0}^{\frac m6-1}r^{3l})\sum_{l=0}^{\frac m6-1}r^{3l}-1}
bc^{\frac m2},
\end{array}$$
that is
\begin{eqnarray}\label{f6.5}
\begin{array}{ll}
&(r^2-2sr^3-x(1-2sr^3)^{d+\frac m6}-2ur^2\sum_{l=0}^{\frac m6-1}r^{3l}+2z(r+r^2+r^3)+\frac n2)\sum_{l=0}^{\frac m6-1}r^{3l}\equiv \\
&x(1-2sr^3)^d+2u\sum_{l=1}^dr^{3l}+1(\mod n).
\end{array}
\end{eqnarray}

(v) $\pi$ preserves $b^c=a_1^xb$:
$$\begin{array}{ll}
&(b^c)^a=(a_1^{-1}b)^{a_1^{1+2z}c^{2+3d}}=a_1^{x-1-2(1+2z)r+\frac n2+2u\sum_{l=1}^dr^{3l}}b,\\
&(a_1^xb)^a=a_1^{x-1}b,
\end{array}$$
that is
\begin{eqnarray}\label{f6.6}
-2(1+2z)r+\frac n2+2u\sum_{l=1}^dr^{3l}\equiv 0(\mod n).
\end{eqnarray}

(vi) $\o(\pi(c))=m$:
$$(c^a)^m=c_1^{(2+3d)\frac m3}a_1^{(-2sr^3+r^2-x(1-2sr^3)^{d+\frac m6}-2ur^2\sum_{l=0}^{\frac m6-1}r^{3l}+2z(r+r^2+r^3)+\frac n2)\sum_{l=0}^{\frac m3}r^{3l}}.$$
Since $r^{\frac m6}\equiv-1(\mod\frac n2)$, we get $\sum_{l=0}^{\frac m3}r^{3l}\equiv0(\mod \frac n2)$.
Note that $r^2-x(1-2sr^3)^{d+\frac m6}$ is even.
Then $$(-2sr^3+r^2-x(1-2sr^3)^{d+\frac m6}-2ur^2\sum_{l=0}^{\frac m6-1}r^{3l}+2z(r+r^2+r^3)+\frac n2)\sum_{l=0}^{\frac m3}r^{3l}\equiv 0(\mod n),$$ as desired.

(vii) $\pi^2=\Inn(a_1)$:
Recall $c_1^{a_1}=a_1^{2s}c_1$ and $a_1^c=bc^{\frac m2}$,
then we know $\Inn(a_1)(c)=c^{1+\frac m2}a_1^{-1+\frac n2}b$ and $\Inn(a_1)(c_1)=a_1^{2s}c_1$.
$$\begin{array}{ll}
&c^{1+\frac m2}a_1^{-1+\frac n2}b=\Inn(a_1)(c)=\pi^2(c)=a_1^{1+2z}(c^{2+3d})^a
=a_1^{1+2z}(a_1^{1+2z}c^{2+3d})^2(c_1^a)^d\\
&=c^{(3d+2)^2+\frac m2}a_1^{-3+r^{-1}-2(r+2zr+z)+\frac{n(2+\sum_{l=0}^{d-1}r^{3l}+(1-2sr^3)^d)}4
+(2sr^3-r^2-2z(r+r^2+r^3)-1-r)\sum_{l=0}^{d-1}r^{3l}}b,
\end{array}$$
that is,
\begin{eqnarray}\label{f6.7}
\begin{array}{ll}
&(1+d)(1+3d)\equiv 0(\mod\frac m3),\\
&r^{-1}-2(r+2zr+z)+(2sr^3-(2zr+1)(1+r+r^2))\sum_{l=0}^{d-1}r^{3l}\equiv \\
&2+\frac{n(\sum_{l=0}^{d-1}r^{3l}+(1-2sr^3)^d)}4(\mod n).
\end{array}
\end{eqnarray}
And
$$
a_1^{2s}c_1=\Inn(a_1)(c_1)=\pi^2(c_1)=c_1a_1^{(2z+1)(r+r^2+r^3)(\sum_{l=0}^{1+3d}r^{3l}+1)+\frac n2},
$$
that is,
\begin{eqnarray}\label{f6.8}
((2z+1)(r+r^2+r^3)+\frac n4)(\sum_{l=1}^{1+3d}r^{3l}+2)\equiv 2sr^3(\mod n).
\end{eqnarray}

(viii) Insure $\lg c\rg_X=1$:
Since $\lg c\rg_X\leq M$, we get $\lg c\rg_X\leq M_X=\lg a_1^2\rg\lg c_1\rg$, which implies
$\lg c\rg_X=\cap_{x\in X}C^x=\cap_{x\in G}C^x=\cap_{x\in G}\lg c_1\rg^x=\lg c_1\rg_{X_1}$.
Then it is suffer to insure $\lg c^3\rg_{X_1}=1$.

Recall
$$X_1=\lg a,b,c_1|R,(a_1)^{c_1}=a_1^{r^3}, c_1^b=a_1^{2u}c_1,
(c_1)^a=c_1^{2+3d}a_1^{i}\rg,$$
where $r^{\frac m2}\equiv-1(\mod \frac n2),\,2u\equiv\frac{r^3-1}{r^2}(\mod\frac n2)$
and $i\equiv\frac {1-r^3}2+\frac n4(\mod\frac n2)$.
Since $(u,\frac n4)=1$, we get $(i,\frac n4)=1$.
Noting $\lg c_1\rg_{X_1}=1$, by Lemma~\ref{Qac}, we get
$$\sum_{l=1}^{j}r^{3l}\equiv0(\mod\frac n4)\Leftrightarrow j\equiv0(\mod\frac m3).$$
Note that $\sum_{i=1}^{1+3d}r^i\equiv0(\mod\frac n2)$. Thus $1+3d\equiv0(\mod\frac m3)$.
Since $1+3d\ne0$, we get $1+3d$ is either $\frac m3$ or $\frac {2m}3$.
By Eq(\ref{f6.6}), we get $1+2z\equiv \frac{1-r}{2r}(\mod\frac n2)$.
\qed
\subsection{$M=\lg a^3\rg\lg c\rg $ and $X/M_X\cong S_4$}
\begin{lem}\label{Qa3c4}
Suppose that $X=X(Q)$, $M=\lg a^3\rg\lg c\rg$, $X/M_X\cong S_4$ and $\lg c\rg_X=1$. Then
\begin{eqnarray}\label{m5.0}
X=\lg a,b,c|R,a^{c^4}=a^r, b^{c^4}=a^{1-r}b, (a^3)^{c^{\frac m4}}=a^{-3}, a^{c^{\frac m4}}=bc^{\frac{3m}4}\rg,
\end{eqnarray}
where  $m\equiv 4(\mod8)$ and $r$ is of order $\frac m4$ in $\ZZ_{2n}^*$.
\end{lem}

In this case,  $M_X=\lg a^3\rg\lg c^4\rg$.
Set $a^3=a_1$ and $c^4=c_1$ so that $M_X=\lg a_1\rg \lg c_1\rg $.
Set $\o(a)=2n$ and $\o(c)=m$, where $n\equiv0(\mod3)$ and $m\equiv0(\mod4)$.
Then in  Lemma~\ref{Qa_1}, we shall show $\lg a_1\rg \lhd X$ and
in Lemma\ref{Qa3c4.1}, we shall get the classification of $X$.
\begin{lem}\label{Qa_1}
$\lg a_1\rg\lhd X$.
\end{lem}
\demo
Let $X_1=M_XG$.
Since $\lg a\rg\lg c_1\rg\leq X_1$ and $\lg c_1\rg_{X_1}=1$,
the subgroup $X_1$ has been given in Lemma~\ref{Qac}:
\begin{eqnarray} \label{m5.1}
X_1=\lg a,b,c_1| R,\,(a^2)^{c_1}=a^{2r},\,(c_1)^a=a_1^{2s}c_1^{t},\,(c_1)^b=a_1^{u}c_1^v\rg,
\end{eqnarray}
where
$$\begin{array}{ll}
&r^{t-1}\equiv r^{v-1}\equiv1(\mod n),\, t^2\equiv 1(\mod\frac m4),\\
&6s\sum_{l=1}^tr^{l}+6sr\equiv 6sr+6s\sum_{l=1}^{v}r^l-3u\sum_{l=1}^tr^l+3ur
  \equiv 2(1-r)(\mod 2n),\\
&3u(\sum_{l=0}^{v-1}r^l-1)\equiv 0(\mod 2n),\,{\rm and}\,v^2\equiv 1(\mod \frac m4),\,
  {\rm if}\,2\di n,\\
&3u\sum_{l=1}^vr^l-ur\equiv 6sr+r-1(\mod 2n),\, {\rm and}\,v^2\equiv t(\mod \frac m4),\,
  {\rm if}\,2\nmid n,\\
&2\di u,\,{\rm if}\,t\ne1,\\
&2s\sum_{l=1}^{w}r^l\equiv u\sum_{l=1}^{w}(1-s(\sum_{l=1}^tr^l+r))^l\equiv 0(\mod \frac{2n}3)
  \Leftrightarrow w\equiv 0(\mod \frac m4).
 \end{array}$$

Now $X=\lg X_1, c\rg $. Since $X/M_X=\lg\ola,\olb\rg\rtimes\lg\olc\rg\cong S_4$,
the only possibility under our conditions is:
\begin{eqnarray}\label{m5.2}
\ola^3=\olc^4=\olb^2=1, \ola^{\olb}=\ola^{-1},\, \ola^{\olc}=\ola^i\olb \olc^3 ,
\end{eqnarray}
where $i\in \ZZ_3$.
Observing Eq(\ref{m5.1}) and Eq(\ref{m5.2}), we may  relabel $a^ib$ by $b$.
Then in the perimage $X$, Eq(\ref{m5.2}) corresponds to
\begin{eqnarray}\label{m5.3}
a^3=a_1, b^2=a^n, c^4=c_1, a^c=bc^{3+4w}.
\end{eqnarray}
Set
\begin{eqnarray}\label{m5.4}
(a_1)^c=a_1^{z}c_1^{d},
\end{eqnarray}
necessarily, $z$ is odd, as $\o(a_1)$ is even.
Then $X$ is uniquely determined by Eq(\ref{m5.1}), Eq(\ref{m5.3}) and Eq(\ref{m5.4}).
To show $\lg a_1\rg\lhd X$, for the contrary, we assume $d\ne 0$.
Then we need to deal with two cases according to the parameter $t$  of $X_1$, separately.
\vskip 3mm
{\it Case 1: $t=1$}
\vskip 3mm
In this case, $v=1$ and $1-6sr-r\equiv0(\mod n)$ by $X_1$.
Set $r_1=1-6sr$. Then $r_1\equiv1(\mod6),\,a^{c_1}=a^{r_1}$ and $b^{c_1}=a_1^{ur_1}b$.
By Eq(\ref{m5.1}), Eq(\ref{m5.3}) and Eq(\ref{m5.4}),
one can check $b^c=a^{2+3x}bc^{4y}$ for some $x$ and $y$.

Since $c$ preserves  $a_1^{c_1}=a_1^{r_1}$, there exist some $x$ such that
$$((a_1)^{c_1})^c=a_1^{zr_1}c_1^{d}=c_1^{d}a_1^{zr_1^{1+d}}\quad{\rm{and}}\quad
(a_1^{r_1})^c=(a_1^{z}c_1^{d})^{r_1}=c_1^{dr_1}a_1^{3z\sum_{l=1}^{r_1}r_1^{dl}},$$
which gives
\begin{eqnarray}\label{f7.1}
d\equiv dr_1(\mod\frac m4).
\end{eqnarray}
Since $c$ preserves $b^{c_1}=a_1^{ur_1}b$, we get
$$(b^{c_1})^c=a^{{2+3x}r_1+3ur_1}bc^{4y}\quad{\rm{and}}\quad
(a_1^{ur_1}b)^c=c_1^{dur_1}
a_1^{z\sum_{l=1}^{ur_1}r_1^{dl}}a^{2+3x}bc^{4y},$$
which gives
\begin{eqnarray}\label{f7.2}
du\equiv 0(\mod\frac m4).
\end{eqnarray}
Since $a_1^{c_1}=a_1^{c^4}=a_1^{x_1}c_1^{d(z^3+z^2+z+1)}$ for some $x_1$, we get
\begin{eqnarray}\label{f7.3}
d(z^3+z^2+z+1)\equiv 0(\mod\frac m4).
\end{eqnarray}
By Eq(\ref{m5.3}), we get  $ac=cbc^{3+4w}$. Then
$$a_1^{ac}=a_1^z c_1^d\quad {\rm{and}}\quad a_1^{cbc^{3+4w}}=a_1^{x_1}c_1^{d-dz(z^2+z+1)},$$
which gives
\begin{eqnarray}\label{f7.4}
dz(z^2+z+1)\equiv0(\mod\frac m4).
\end{eqnarray}
With Eq(\ref{f7.3}), we know $d\equiv0(\mod\frac m4),$ contradicting with $d=0$.
\vskip 3mm
{\it Case 2: $t\ne 1$}
\vskip 3mm
In this case, we have $t\ne1$.
Suppose  that $\lg a_1^2\rg\lhd X$.
Then  in what follows we shall show  $\lg a_1^2\rg\lhd X$. If so,
then  by  considering $\ox=X/\lg a_1^2\rg$ and $\lg \overline{c_1}\rg\lhd\ox$,
 one  may get   $t=1$,  a contradiction,  as  $\oc\leq C_{\ox}(\lg \overline{c_1}\rg)=\ox$.

By Eq(\ref{m5.1}),  $u$ is even and $r\equiv1(\mod3)$.
By using Eq(\ref{m5.1}), Eq(\ref{m5.3}) and Eq(\ref{m5.4}),
one may derive  $b^c=a^{2+6x}bc_1^y$ for some $x$ and $y$, omitting the details.

Since $c$ preserves  $(a_1^2)^{c_1}=a_1^{2r}$, there exist some $x$ such that
$$((a_1^2)^{c_1})^c=(a_1^zc_1^da_1^zc_1^d)^{c_1}=a_1^{x}c_1^{d(t+1)}\quad{\rm and}\quad
(a_1^{2r})^c=(a_1^zc_1^d)^{2r}
=a_1^{x}c_1^{d(t+1)r},$$
which gives
\begin{eqnarray}\label{f7.1.1}
d(t+1)(r-1)\equiv0(\mod\frac m4).
\end{eqnarray}
Since $c$ preserves $(c_1)^b=a_1^{u}c_1^v$, we get
$$(c_1^b)^c=c_1^{a^{2+6x}bc_1^y}=a_1^{x}c_1^v,\,
(a_1^{u}c_1^v)^c=(a_1^zc_1^d)^{u}c_1^v=a_1^{x}c_1^{\frac{du(t+1)}2+v},$$
which gives
\begin{eqnarray}\label{f7.2.1}
\frac{du(t+1)}2\equiv0(\mod\frac m4).
\end{eqnarray}
Since $c$ preserves $c_1^a=a_1^{2s}c_1^t$, we get $c_1^{bc^{2+4w}}=a_1^{2s}c_1^t$.
Then
$$c_1^{bc^{2+4w}}=(a_1^{u}c_1^{v})^{c^{2+4w}}=((a_1^zc_1^d)^zc_1^d)^{ur^w}c_1^{v}=a_1^xc_1^{v},$$
which gives
\begin{eqnarray}\label{f7.3.1}
v\equiv t(\mod\frac m4).
\end{eqnarray}
By Eq(\ref{m5.3}) again, we get  $ac^2=cbc^{4(w+1)}$. Then
$$(a_1^2)^{ac^2}=a_1^xc_1^{d(t+1)(z+1)}\,{\rm{and}}\,(a_1^2)^{cbc_1^{w+1}}=a_1^xc_1^{d(t+1)},$$
which gives
\begin{eqnarray}\label{f7.4.1}
dz(t+1)\equiv0(\mod\frac m4).
\end{eqnarray}
Since $\o(a_1^c)=\frac {2n}3$, we get $\frac{dn(t+1)}3\equiv0(\mod\frac m4)$.
With $(\frac n3,z)=1$ and Eq(\ref{f7.4.1}), we get $d(t+1)\equiv0(\mod\frac m4)$.
Then $(a_1^2)^c=(a_1^zc_1^d)^2=a_1^xc_1^{d(t+1)}=a_1^x$ for some $x$, which implies $\lg a_1^2\rg\lhd X$,
as desired.
\qed
\begin{lem} \label{Qa3c4.1}
The group $X$ is given by Eq(\ref{m5.0}).
\end{lem}
\demo
By lemma, $\lg a_1\rg \lhd X$,  that is $(a_1)^c=a_1^z$ by Eq(\ref{m5.4}).
Since $\lg a^2\rg \lhd X_1$, we get $\lg a\rg \lhd X_1$ and so $G\lhd X_1$, that is  $t=v=1$ in Eq(\ref{m5.1}).
Then by Eq(\ref{m5.1}), (\ref{m5.3}) and (\ref{m5.4}), we can set
$$X=\lg a,b,c\di R,a^{c_1}=a^{r_1},b^{c_1}=a_1^{ur_1}b,(a_1)^c=a_1^z,a^c=bc^{3+4w}\rg,$$
where
$$\begin{array}{ll}
&r_1=1-6sr,\,r_1^{\frac m4}-1\equiv2(r_1-r)\equiv0(\mod 2n),\\
&2s\sum_{l=1}^jr^l\equiv0\equiv u\sum_{l=1}^jr_1^l(\mod\frac{2n}3)\Leftrightarrow j\equiv0(\mod\frac m4).
\end{array}$$
Note $2s\sum_{l=1}^jr^l\equiv0(\mod\frac{2n}3)\Leftrightarrow r_1^j-1\equiv0(\mod 2n)$.

In what follows, we shall divide the proof into two steps:
\vskip 3mm
{\it Step 1: Show $m\equiv4(\mod8)$.}
\vskip 3mm
Set $\lg c\rg=\lg c_2\rg\times\lg c_3\rg$,
where $\lg c_2\rg$ is a $2-$group and $\lg c_3\rg$ is the $2'-$Hall subgroup of $\lg c\rg$.
Then  $\lg c_1\rg=\lg c_2^4\rg\times\lg c_3\rg$.
To show $m\equiv4(\mod8)$, we only show $c_2^4=1$.

Consider $\ox=X/\lg a_1\rg=\lg\ola,\olb\rg\lg\olc\rg$.
Then one can check $C_{\ox}(\lg\olc_1\rg)=\ox$, which implies
$\lg\overline{c_1}\rg\le Z(\ox)$ and $\ox/\lg \overline{c_1}\rg\cong S_4$.
Note that $\lg\olc_3\rg\leq\lg\olc_3\rg(\lg\olc_2^4\rg\lg\ola\rg)\leq\ox$
where $(|\ox:\lg\olc_3\rg(\lg\olc_2^4\rg\lg\ola\rg)|,|\lg\olc_3\rg|)=1$.
Thus by Proportion~{\ref{complement}}, we get that $\lg\olc_3\rg$ has a complement in $\ox$,
which implies $X=(\lg a,b\rg\lg c_2\rg)\rtimes\lg c_3\rg.$

Consider $X_2=\lg a,b\rg\lg c_2\rg$, where $\lg c_2\rg_{X_2}=1$ and $\lg a_1\rg\lhd X_2$,
and $\overline{X_2}=X_2/\lg a_1\rg=\lg\ola,\olb\rg\lg \overline{c_2}\rg$.
Note that $\lg\overline{c_2}^4\rg\lhd \overline{X_2}$.
Then one can check $C_{\ox_2}(\lg\overline{c_2}^4\rg)=\ox_2$,
which implies that $\ox_2$ is the central expansion of $S_4$.
By Lemma~{\ref{Schur}}, we get the Schur multiplier of $S_4$ is $\ZZ_2$,
and then $\o(c_2)$ is either $4$ or $8$.
Suppose that $\o(c_2)=8$.
Then $c_2^4$ normalises $G$, and we set $a^{c_2^4}=a^i$, where $i\equiv1(\mod3).$
Note that $\lg a\rg\leq C_{X_2}(\lg a_1\rg)\lhd X_2$.
Then $\lg a,bc_2,c_2^2\rg\leq C_{X_2}(\lg a_1\rg),$
which implies $\lg a_1\rg\times\lg c_2^4\rg\lhd X_2$.
Since $i\equiv1(\mod\frac {2n}3)$ and $i^2\equiv 1(\mod 2n)$,
we get $i=1$, which implies $[a,c_2^4]=1$.
Then $\lg a,c_2^2\rg\leq C_{X_2}(\lg a_1\rg\times\lg c_2^4\rg)\lhd X_2$,
which implies $bc_2\in C_{X_2}(\lg a_1\rg\times\lg c_1\rg)$. So $(c_2^4)^b=c_2^4$.
Then $c_2^4\lhd X_2$, a contradiction.
So $\o(c_2)=4$, which implies $\lg c_1\rg=\lg c_3\rg$.
Then $X=(\lg a,b\rg\lg c_2\rg)\rtimes\lg c_1\rg.$
\vskip 3mm
{\it Step 2:  Determine the parameters $r_1,u,w$ and $z$.}
\vskip 3mm
In $X_2=\lg a,b\rg\lg c_2\rg$, we know $a^{c_2}=bc_2^3$ by Eq(\ref{m5.3}).
Consider $\lg a\rg\leq C_{X_2}(\lg a_1\rg)\lhd X_2$,
then $C_{X_2}(\lg a_1\rg)$ is either $\lg a,bc_2,c_2^2\rg$ or $X_2$.

Suppose that $\lg c\rg_X(\lg a_1\rg)=X$. Then we know $a_1=b^2$ as $[a_1,b]=1$, that is $n=3$ and $a_1=a_1^{-1}$.
Then one can check $m=4$ and $z,r=1$, as desired.

Suppose that $\lg c\rg_X(\lg a_1\rg)=\lg a,bc_2,c_2^2\rg$.
Then $a_1=a_1^{bc_2}=(a_1^{-1})^{c_2}$, which implies $z=-1$.
In $X_1=\lg a,b\rg\rtimes\lg c_1\rg$, we know $X_1=\lg a,b,c_1\di R,a^{c_1}=a^{r_1},b^{c_1}=a_1^{ur_1}b\rg.$
Since $c_2$ preserves  $a^{c_1}=a^{r_1}$, we get
$$(bc_2^3)^{c_1}=a_1^{ur_1}bc_2^3\quad{\rm{and}}\quad
(a^{r_1})^{c_2}=(a_1^{\frac{r_1-1}3})^{c_2}a^{c_2}=a_1^{\frac{1-r_1}3}bc_2^3,$$
which gives
$$ur_1\equiv \frac{1-r_1}3(\mod\frac {2n}3).$$
Recall that $r_1^j-1\equiv0\equiv 3u\sum_{l=1}^jr_1^l(\mod 2n)\Leftrightarrow j\equiv0(\mod\frac m4)$.
Then we get that $r_1^j-1\equiv0(\mod 2n)\Leftrightarrow j\equiv0(\mod\frac m4)$,
which implies $\o(r_1)=\frac m4$.
For the purpose of formatting uniformity, replacing $r_1$ by $r$, then we get Eq(\ref{m5.0}),
as desired.
\qed
\section{Proof of Theorem~\ref{main4}}
To prove Theorem~1.4, let  $\lg c\rg _X=1$ and set $R:=\{a^{n}=b^2=c^m=1,\,a^b=a^{-1}\}$.
Then we shall deal with the five cases in Theorem 1.1 in the following five subsections, separately.

\subsection{$M=\lg a\rg \lg c\rg $}
\begin{lem}\label{Dac}
Suppose  that $X=X(D)$, $M=\lg a\rg \lg c\rg $ and $\lg c\rg _X=1$. Then
\begin{eqnarray}\label{Dm1}
X=\lg a,b,c| R, (a^2)^c=a^{2r},  c^a=a^{2s}c^t,c^b=a^uc^v\rg,
\end{eqnarray}
 $$\begin{array}{ll}
  &2(r^{t-1}-1)\equiv 2(r^{v-1}-1)\equiv u(\sum_{l=0}^{v-1}r^l-1)\equiv0(\mod n),\,
  t^2\equiv v^2\equiv1(\mod m)\\
  &2s\sum_{l=1}^tr^{l}+2sr\equiv 2sr+2s\sum_{l=1}^{v}r^l-u\sum_{l=1}^tr^l+ur\equiv 2(1-r)(\mod n),\\
  &{\rm if}\, t\ne 1,\, {\rm then}\, u\equiv0(\mod2),\\
  &2s\sum_{l=1}^{w}r^l\equiv u\sum_{l=1}^{w}(1-s(\sum_{l=1}^tr^{l}+r))^l\equiv0(\mod n)\Leftrightarrow w\equiv0(\mod m).
  \end{array}$$
Moreover, $G_X=\lg a^2\rg $ if $tv,t,v\ne 1$;
$\lg a^2, b\rg $ if $v=1$ and $u$ is even but $t\ne 1$;
$\lg a^2,ab\rg$ if $tv=1$ but $t,v\ne1$; and  $\lg a, b\rg$ if $t=v=1$, respectively.
\end{lem}
\demo Noting  $\lg a^2\rg\lhd X$ and $M=\lg a\rg\lg c\rg\le X$,  we have $X$ may be obtained  by three cylic extension of groups  in order:
$$\lg a^2\rg\rtimes \lg c\rg, \quad  (\lg a^2\rg\rtimes \lg c\rg).\lg a\rg \quad  {\rm and}\quad
((\lg a^2\rg\rtimes \lg c\rg).\lg a\rg) \rtimes \lg b\rg .$$
So $X$ has the presentation as in Eq(\ref{Dm1}).
What we should to  determine the parameters  $r,s,t,u$ and $v$ by analysing three extensions.
\vskip 3mm
(1) $\lg a^2\rg\rtimes \lg c\rg$, where $(a^2)^c=a^{2r}$. Set $\pi_1\in \Aut(\lg a^2\rg)$ such that $(a^2)^{\pi_1}=a^{2r}$.
\vskip 3mm
This extension
 is valid if and only if  $\o((a^{2})^{\pi_1})=\o(a^2)$ and $\pi_1^{m}=\Inn(a^2)$, that is
 \begin{eqnarray}\label{Df3.1}
 2(r^m-1)\equiv 0(\mod n).
 \end{eqnarray}

(2) $(\lg a^2\rg\rtimes \lg c\rg).\lg a\rg$, where $c^a=a^{2s}c^t$. Set $\pi_2\in \Aut((\lg a^2\rg\rtimes \lg c\rg)$: $a\to a$ and $c\to a^{2s}c^t$.
\vskip 3mm

(i) $\pi_2$   preserves  $(a^2)^c=a^{2r}$:
\begin{eqnarray}\label{Df3.2}
2(r^{t-1}-1)\equiv  0(\mod n).
\end{eqnarray}

(ii)  $\o(\pi_2(c))=m$:
  $$(a^{2s}c^t)^m=c^{tm}a^{2s\sum_{l=1}^mr^{tl}}=c^{tm}a^{2s\sum_{l=1}^m r^{l}}=1,$$
\f that is
\begin{eqnarray}\label{Df3.3}
2s\sum_{l=1}^{m} r^{l}\equiv  0(\mod n).
\end{eqnarray}

(iii) $\pi_2^2=\Inn(a^2)$:
$$ca^{2-2r}=\Inn(a^2)(c)=\pi_2^2(c)=(a^{2s}c^t)^{a}
=c^{t^2}a^{2sr+2s\sum_{l=1}^{t} r^{l}},$$
that is
\begin{eqnarray}\label{Df3.4}
t^2-1\equiv 0(\mod m) \quad {\rm and}\quad 2(s\sum_{l=1}^{t} r^{l}+rs+r-1)\equiv 0(\mod n).
\end{eqnarray}
\vskip 3mm
(3) $((\lg a^2\rg\rtimes \lg c\rg).\lg a\rg)\rtimes \lg b\rg $, where $c^b=a^uc^v$.
Set $\pi_3\in \Aut((\lg a^2\rg\rtimes \lg c\rg).\lg a\rg) :$ $a\to a^{-1}$ and $c\to a^uc^v$.
We divide two cases, separately.

\vskip 3mm
{\it Case 1: $u$ is even.}
\vskip 3mm

(i) $\pi_3$ preserves  $(a^2)^c=a^{2r}$:
\begin{eqnarray}\label{Df3.5}
2(r^{v-1}-1)\equiv 0(\mod n).
\end{eqnarray}

(ii) $\o(\pi_3(c))=m$:
$$1=(a^uc^v)^m=c^{vm}a^{u\sum_{l=1}^{m}r^l},$$
that is
\begin{eqnarray}\label{Df3.6}
u\sum_{l=1}^{m}r^l\equiv0(\mod n).
\end{eqnarray}

(iii)  $\pi_3$ preserves $c^a=a^{2s}c^t$:  that is
$$a^uc^v=(a^{-2s}(a^uc^v)^t)^a=a^{-2s}c^{t^2v}a^{(ru+2s\sum_{l=1}^{v}r^l)\sum_{l=0}^{t-1}r^l},$$
which implies
\begin{eqnarray}\label{Df3.7}
r(u+2s)\equiv (ru+2s\sum_{l=1}^{v}r^l)\sum_{l=0}^{t-1}r^l (\mod n).
\end{eqnarray}
With Eq(\ref{Df3.4}), Eq(\ref{Df3.7}) is if and only if
\begin{eqnarray}\label{Df3.8}
u\sum_{l=1}^{t}r^l-ur-2s\sum_{l=1}^{v}r^l-2sr+2(1-r)\equiv 0 {\pmod n}.
\end{eqnarray}

(iv) $\pi_3^2=1$:
Suppose  that $v=1$. Then $c=c^{b^2}=a^{-u}a^uc=c$, as desired.
Suppose that $v\neq 1$. Then
$$c=c^{b^2}=a^{-u}(a^uc^v)^v=c^v(a^uc^v)^{v-1}=c^{v^2}a^{u\sum_{l=1}^{v-1}r^l},$$
that is
\begin{eqnarray}\label{Df3.9}
u\sum_{l=1}^{v-1}r^l\equiv0(\mod n)\,{\rm{and}}\,v^2-1\equiv 0(\mod m).
\end{eqnarray}
Then $\pi_3^2=1$ is if and only if
\begin{eqnarray}\label{Df3.9.1}
u\sum_{l=1}^vr^l-ur\equiv0(\mod n)\quad{\rm{and}}\quad v^2-1\equiv 0(\mod m).
\end{eqnarray}

\vskip 3mm
{\it Case 2: $u$ is odd.}
\vskip 3mm
If $t=1$, then by Lemma~\ref{a}, we get $G\lhd X$, which implies $v=1$.
So assume $t\ne 1$ and we shall get a contradiction.

Let $S=\lg a^2, c\rg $. Since $u$ is odd again, we know that $\lg a^2\rg \le S_X\lnapprox S$.
Since $|X:S|=4$, we have $\ox=X/S_X=\lg \olc, \ola\rg \rtimes \lg \olb\rg \lessapprox S_4$. The only possibility is
$\o(\olc)=2$ and   $\olx \cong D_8$ so that $m$ is even and $v$ is odd.
Then $t$ is odd, as $t^2\equiv 1{\pmod m}$.
Moreover, we have $\lg a^2, c^2\rg=S_X\lhd X$.

Consider $\ox=X/\lg a^2\rg=\lg\ola,\olc\rg\rtimes\lg\olb\rg$,
where $\ola^{\olb}=\ola,\,\olc^{\ola}=\olc^t$ and $\olc^{\olb}=\ola\olc^v$.

Let  $\pi_3$ be defined as above.
Since the induced action of $\pi_3$ preserves  $\olc^{\ola}=\olc^t$, we have
$(\ola\olc^v)^{\ola}=(\ola\olc^v)^t$, that is
$$\ola\olc^{tv}=\ola\olc^v((\ola\olc^v)^2)^{\frac{t-1}2}=\ola\olc^v(\olc^{tv+v})^{\frac{t-1}2}
=\ola\olc^{v+\frac{v(t+1)(t-1)}2},$$
which implies
$$tv\equiv v+\frac{v(t+1)(t-1)}2{\pmod m}.$$
\f Noting  $t^2\equiv1{\pmod m}$, $t\ne 1$  and $(v, m)=1$ is odd, we get
\begin{eqnarray}\label{Df3.10}
t\equiv 1+\frac m2(\mod m).
\end{eqnarray}

Let $X_1=GS_X=\lg a, b\rg \lg c^2\rg $.
By Eq(\ref{Df3.10}), we have $(c^2)^a=(a^{2s}c^t)^2=a^{2s(1+r^{-1})}c^2$,
which implies $c^2$ normalises $\lg a\rg$.
Then we get $G\lhd X_1$ and so $b^{c^2}\le G$.
Since $b^{c^2}=c^{-2}(bc^2b)b=c^{-2}(a^uc^v)^2b=c^{v(t+1)-2}a^xb,$ for some $x$,
we get  $v(t+1)-2\equiv0 (\mod m)$.
By combing Eq(\ref{Df3.10})  we get
\begin{eqnarray}\label{Df3.11}
v\equiv 1\pm \frac m4\, (\mod \frac m2), \quad 4\di m.
\end{eqnarray}
Since $\olc=\olc^{\olb^2}=\ola(\ola\olc^v)^v=\olc^v(\ola\olc^v\ola\olc^v)^{\frac{v-1}2}
=\olc^{v+(tv+v)\frac{v-1}2},$ we get
\begin{eqnarray}\label{Df3.12}
(v-1)(\frac{v(t+1)}2+1) \equiv 0(\mod m).
\end{eqnarray}
Then  Eq(\ref{Df3.11}) and  Eq(\ref{Df3.12}) may give $\frac m2\equiv 0(\mod m)$, a contradiction.
\vskip 3mm
(4) Insure $\lg c\rg_X=1$:
When $u$ is even, for any integer $w$, we get
$$(c^w)^a=(a^{2s}c^t)^w=c^{tw}a^{2s\sum_{l=1}^{w}r^l}\quad{\rm{and}}\quad (c^w)^b=(a^uc^v)^w=c^{vw}a^{u\sum_{l=1}^{w}r^l}.$$
Since $\lg c\rg_X=1$,  we know that
$2s\sum_{l=1}^{w}r^l\equiv  0\equiv u\sum_{l=1}^{w}r^l(\mod n)$ is if and only if $w\equiv0(\mod m)$.

When $u$ is odd, we know $t=v=1$. Then $2(1-2sr)\equiv 2r(\mod n)$  by Eq(\ref{Df3.4}).
For any integer $w$,
$$(c^w)^a=(a^{2s}c)^w=c^wa^{2s\sum_{l=1}^{w}r^l}\quad{\rm{and}}\quad
(c^w)^b=(a^uc)^w=c^{w}a^{u\sum_{l=1}^{w}(1-2sr)^l}.$$
Since $\lg c\rg_X=1$,  we know that
$2s\sum_{l=1}^{w}r^l\equiv0\equiv u\sum_{l=1}^{w}(1-2sr)^l(\mod n)$ is if and only if $w\equiv0(\mod m)$.

Summarizing  Eq(\ref{Df3.1})-Eq(\ref{Df3.9.1}), we get  the parameters $(m,n, r, s, t, u, v)$ as shown in the lemma.

\qed

\subsection{$M=\lg a^2\rg\lg c\rg $ and $X/M_X\cong D_8$}
\begin{lem}\label{DD8}
Suppose  that $X=X(D)$,  $M=\lg a^2\rg\lg c\rg$,   $X/M_X\cong D_8$ and $\lg c\rg_X=1$. Then
$$X=\lg a,b,c|R,(a^2)^{c^2}=a^{2r},(c^2)^a=a^{2s}c^{2t},(c^2)^b=a^{2u}c^{2}, a^c=bc^{2w}\rg,$$
where either $w=s=u=0$ and $r=t=1$; or
$$\begin{array}{ll}
  &w\neq0,\,s=u^2\sum_{l=0}^{w-1}r^l,\,t=1+2wu,\\
  &nw\equiv2w(r-1)\equiv2w(1+uw)\equiv0(\mod\frac m2),\\
  &r^{2w}-1\equiv(u\sum_{l=1}^{w}r^l)^2-r\equiv (r^w+1)(1+s\sum_{l=0}^{w-1}r^l) \equiv0(\mod\frac n2),\\
  &\sum_{l=1}^{i}r^l\equiv0(\mod\frac n2)\Leftrightarrow i\equiv0(\mod\frac m2).
  \end{array}$$
\end{lem}
\demo Under the hypothesis, $M_X=\lg a^2\rg \rtimes\lg c^2\rg$. Set $n=\o(a)$ and $m=\o(c)$.
Then  both $n$ and $m$ are even.
Since $X/M_X=\lg\ola,\olb\rg\lg\olc\rg\cong D_8$, we can choose $\olb$
such that  the form of $X/M_X$ is the following:
$\ola^{\olc}=\olb$ and $\olb^{\olc}=\ola.$
Set $c_1:=c^2$ and $X_1=GM_X=\lg a,b\rg\lg c_1\rg$.
Noting $\lg a\rg\lg c_1\rg\leq X_1$ and $\lg c_1\rg_{X_1}=1$, by  Lemma~\ref{Dac}, we get
$$X_1=\lg a,b,c_1|R, (a^2)^{c_1}=a^{2r},c_1^a=a^{2s}c_1^{t}, c_1^b=a^{2u}c_1^v\rg,$$
where
\begin{eqnarray}\label{Df4.1}
\begin{array}{ll}
&r^{t-1}-1\equiv r^{v-1}-1\equiv u\sum_{l=0}^{v-1}r^l-u\equiv0(\mod\frac n2),\,
t^2\equiv v^2\equiv 1(\mod\frac m2),\\
&s\sum_{l=1}^{t} r^{l}+sr\equiv sr+s\sum_{l=1}^{v}r^l-u\sum_{l=1}^tr^l+ur\equiv 1-r(\mod\frac n2),\\
&s\sum_{l=1}^{i}r^l\equiv u\sum_{l=1}^{i}r^l\equiv0(\mod\frac n2)\Leftrightarrow i\equiv0(\mod\frac m2).
\end{array}
\end{eqnarray}
Now $X=X_1.\lg c\rg$. Set $a^c=bc_1^w$. Then $X$ may be defined by $R$ and
\begin{eqnarray}\label{Dm2}
(a^2)^{c_1}=a^{2r},\,c_1^a=a^{2s}c^{2t},\,c_1^b=a^{2u}c^{2v},\,a^c=bc_1^w.
\end{eqnarray}

If $w\equiv0(\mod\frac m2)$, then $\o(a)=\o(a^c)=\o(b)=2$, which implies $X\cong D_8$.
Then $w=s=u=0$ and $r=t=1$, as desired.
So in that follows, we assume $w\not\equiv 0(\mod\frac m2)$.

Firstly, we get $ b^c=a^{c^2}c_1^{-w}=a^{1-2sr}c_1^{1-t-w}.$
Set $\pi\in \Aut(X_1) :$
$a\to bc_1^{2w}$, $b\to a^{1-2sr}c_1^{1-t-w}$ and $c_1\to c_1$.
We need to carry out the following  seven steps:
\vskip 3mm
(i) $\o(\pi(b))=2:$
$$(c^{2(w+t-1)}a^{2sr-1})^2=c^{2w(t+1)}a^{2sr^{w+1}+2s\sum_{l=1}^{w+t-1}r^{l}+2sr-2}=1,$$
that is
\begin{eqnarray}\label{Df4.2}
w(t+1)\equiv 0(\mod\frac m2)\quad{\rm{and}}\quad
sr^{w+1}+s\sum_{l=1}^{w+t-1}r^{l}+sr-1\equiv 0(\mod\frac n2),
\end{eqnarray}
which implies
\begin{eqnarray}\label{Df4.3}
r^{2w}\equiv r^{w(t+1)}\equiv1(\mod\frac n2).
\end{eqnarray}

\vskip 3mm

(ii)  $\o(\pi(a))=n$:
$$(bc^{2x})^n=(c^{2w(v+1)}a^{2u\sum_{l=w+1}^{2w}r^l})^{\frac n2}
  = c^{nw(v+1)}a^{un\sum_{l=w+1}^{2w}r^l}=1,$$
that is
\begin{eqnarray}\label{Df4.4}
\frac n2w(v+1)\equiv0(\mod\frac m2).
\end{eqnarray}

(iii) $\pi$ preserves $(a^2)^{c_1}=a^{2r}$:
$$((a^2)^{c^2})^c=c^{2w(v+1)}a^{2ur\sum_{l=w+1}^{2w}r^l}\quad{\rm{and}}\quad
 (a^{2r})^c=c^{2wr(v+1)}a^{2ur\sum_{l=w+1}^{2w}r^l},$$
that is
\begin{eqnarray}\label{Df4.5}
w(t+1)(r-1)\equiv 0(\mod\frac m2).
\end{eqnarray}

(iv) $\pi $ preserves $c_1^a=a^{2s}c_1^{t}$:
$$((c^2)^a)^c=c^{2v}a^{2ur^{w+1}}\quad{\rm{and}}\quad
 (a^{2s}c^{2t})^c=c^{2ws(v+1)+2t}a^{2sru\sum_{l=w+1}^{2w}r^l},$$
that is
\begin{eqnarray}\label{Df4.6}
v\equiv ws(v+1)+t(\mod\frac m2)\quad{\rm{and}}\quad
u\equiv su\sum_{l=1}^{w}r^l(\mod \frac n2).
\end{eqnarray}

(v) $\pi$ preserves $c_1^b=a^{2u}c_1^v$:
$$ \begin{array}{lcl}
 ((c^2)^b)^c&=&(c^2)^{a^{1-2sr}c^{2-2t-2w}}=c^{2t}a^{2sr^{2-w}},\\
 (a^{2u}c^{2v})^c&=&(c^{2w(v+1)}a^{2u\sum_{l=w+1}^{2w}r^l})^uc^{2v}=c^{2wu(v+1)+2v}
                 a^{2u^2\sum_{l=w+2}^{2w+1}r^l},
 \end{array}$$
that is,
\begin{eqnarray}\label{Df4.7}
t\equiv wu(v+1)+v(\mod\frac m2)\quad{\rm{and}}\quad
s\equiv u^2\sum_{l=0}^{w-1}r^l(\mod\frac n2).
\end{eqnarray}

(vi) $\pi^2=\Inn(c_1)$:
Recall $\Inn(c_1)(a)=a^{1-2sr}c_1^{1-t},\,\Inn(c_1)(a^2)=a^{2r}$
and $\Inn(c_1)(b)=c_1^{v-1}a^{2ur}b$.

$$a^{1-2sr}c^{2-2t}=\Inn(c_1)(a)=\pi^2(a)=b^cc^{2w}=a^{1-2sr}c^{2-2t-2w+2w},$$
as desired;

$$a^{2r}=\Inn(c^2)(a_1)=\pi^2(a^2)=(c^{2w(v+1)}a^{2u\sum_{l=w+1}^{2w}r^l})^c
 =c^{2w(v+1)(1+uw)}a^{2(u\sum_{l=1}^wr^l)^2},$$
that is
\begin{eqnarray}\label{Df4.8}
w(v+1)(1+uw)\equiv 0(\mod \frac m2)\quad{\rm{and}}\quad
r\equiv (u\sum_{l=1}^{w}r^l)^2 (\mod \frac n2),
\end{eqnarray}
which implies $(u,\frac n2)=1$ and $(\sum_{l=1}^{w}r^l,\frac n2)=1$  as $(r,\frac n2)=1$;
and

$$c^{2(v-1)}a^{2ur}b=\Inn(c_1)(b)=\pi^2(b)
=c^{2(w+t-1)+2w(v+1)(sr-1)+2vw}a^{2sur\sum_{l=1}^{w}r^l}b,$$
that is,
\begin{eqnarray}\label{Df4.9}
v\equiv t+wvsr+wsr(\mod\frac m2)\quad{\rm{and}}\quad
s\sum_{l=1}^{w}r^l\equiv1(\mod\frac n2),
\end{eqnarray}
which implies $(u,\frac n2)=1$.

(vii) Insure $\lg c\rg_X=1$:
Since $\lg c\rg_X\leq M$, we get $\lg c\rg_X\leq M_X=\lg a^2\rg\lg c^2\rg$.
Since $\lg c\rg_X=\cap_{x\in X}C^x=\cap_{x\in G}C^x=\cap_{x\in G}\lg c^2\rg^x=\lg c^2\rg_{X_1}=1$ and
$u\sum_{l=1}^{i}r^l\equiv0\equiv s\sum_{l=1}^{i}r^l(\mod\frac n2)\Leftrightarrow
i\equiv0(\mod \frac m2)$, noting  $(s,\frac n2)=(u,\frac n2)=1$, we have
$\sum_{l=1}^{i}r^l\equiv0(\mod\frac n2)\Leftrightarrow i\equiv0(\mod \frac m2)$.
\vskip 3mm
Now we are  ready to determine the parameters by summarizing  Eq(\ref{Df4.1})-Eq(\ref{Df4.9}).
Since $(r,\frac n2)=(u,\frac n2)=1$ (after Eq(\ref{Df4.8})),  we get from  Eq(\ref{Df4.1}) that
$\sum_{l=1}^{v-1}r^l\equiv0(\mod\frac n2).$
By (vii), $\sum_{l=1}^{i}r^l\equiv0(\mod\frac n2)\Leftrightarrow i\equiv0(\mod \frac m2)$,
which means  $v\equiv 1(\mod \frac m2)$.
Inserting $v=1$ in Eq(\ref{Df4.1})-Eq(\ref{Df4.9}), we get that
$s=u^2\sum_{l=0}^{w-1}r^l$ and  $t=2wu+1$ in Eq(\ref{Df4.7});
$(r^w+1)(1+s\sum_{l=0}^{w-1}r^l) \equiv0(\mod\frac n2)$ in Eq (\ref{Df4.2})
$nw\equiv0(\mod\frac m2)$ in Eq(\ref{Df4.4}); and
$2w(r-1)\equiv2w(1+uw)\equiv0(\mod\frac m2)$ in Eq(\ref{Df4.5}) and (\ref{Df4.8}).
All these are summarized in the lemma.
\qed

\subsection{$M=\lg a^2\rg\lg c\rg$ and $X/M_X\cong A_4$}
\begin{lem}\label{DA4}
Suppose  that $X=X(D)$, $M=\lg a^2\rg\lg c\rg,\,X/M_X\cong A_4$ and $\lg c\rg_X=1$. Then
$$X=\lg a,b,c|R, a^{c^3}=a^r,(c^3)^b=a^{2u}c^{3}, a^c=bc^{\frac {im}2},b^c=a^xb\rg,$$
where $n\equiv 2(\mod 4)$ and either $i=u=0$ and $r=x=1$;
or $i=1$,  $6\di m$, $l^{\frac m2}\equiv-1(\mod\frac n2)$ with $\o(l)=m$, $r=l^3$, $u=\frac{l^3-1}{2l^2}$ and $x\equiv -l+l^2+\frac n2(\mod n)$.
\end{lem}
\demo
Under the hypothesis, $M_X=\lg a^2\rg\rtimes\lg c^3\rg$. Set $n=\o(a)$ and $m=\o(c)$.
Then $n$ is even and $3\di m$.
Since $X/M_X=\lg\ola,\olb\rg\lg\olc\rg\cong A_4$, we can choose $\olb$
such that  the form of $X/M_X$ is the following:
$\ola^{\olc}=\olb$ and $\olb^{\olc}=\ola\olb.$
Set $c_1:=c^3$ and $X_1=GM_X=\lg a,b\rg\lg c_1\rg$.
By  Lemma~\ref{Dac}, we get
$$X_1=\lg a,b,c^3|R, (a^2)^{c_1}=a^{2r},(c_1)^a=a^{2s}c_1^{t}, (c_1)^b=a^{2u}c_1^{v}\rg$$
whose
\begin{eqnarray}\label{Df5.1}
\begin{array}{ll}
&r^{t-1}-1\equiv r^{v-1}-1\equiv u(\sum_{l=0}^{v-1}r^l-1)\equiv0(\mod\frac n2),\,
t^2\equiv v^2\equiv1(\mod\frac m3),\\
&s\sum_{l=1}^tr^{l}+sr\equiv sr+s\sum_{l=1}^{v}r^l-u\sum_{l=1}^tr^l+ur\equiv 1-r(\mod \frac n2),\\
&s\sum_{l=1}^{i}r^l\equiv u\sum_{l=1}^{i}r^l\equiv0(\mod n)\Leftrightarrow i\equiv0(\mod\frac m3).
\end{array}
\end{eqnarray}
Now $X=X.\lg c\rg$. Set $a^c=bc_1^{w}$.
Then $X$ may be defined by $R$ and
\begin{eqnarray}\label{Dm3}
(a^2)^{c_1}=a^{2r},(c_1)^a=a^{2s}c_1^{t}, (c_1)^b=a^{2u}c_1^{v},\,a^c=bc_1^{w},b^c=a^{1+2x}bc_1^{y}.
\end{eqnarray}

If $w\equiv0(\mod\frac m3)$, then $\o(a)=\o(a^c)=\o(b)=2$, which implies $X\cong A_4$.
Then $n=2,\,i=u=s=0$ and $r=x=t=v=1$, as desired.
So in that follows, we assume $w\not\equiv 0(\mod\frac m3)$.

To determine  the parameters  $r,s,t,u,v,w,x$ and $y$, we only to consider the   last extension
$X_1.\lg c\rg $ in Eq(\ref{Dm3}), where $a^c=bc_1^{w}$ and $b^c=a^{1+2x}bc_1^{y}$.
Set $\pi\in \Aut(X_1): a\to bc_1^{w}, \, b\to a^{1+2x}bc_1^{y}, \quad c_1\to c_1.$
We need to carry out the following  eight steps:

(i) $\o(\pi(b))=2$:
$$(a^{1+2x}bc^{3y})^2=ba^{-(1+2x)}c^{3y}a^{1+2x}bc^{3y}
=c^{3y(tv+1)}a^{2r^y(u\sum_{l=1}^{ty}r^l+xr^y-x-s\sum_{l=1}^{y}r^l)}=1,$$
that is
\begin{eqnarray}\label{Df5.2}
y(tv+1)\equiv0(\mod\frac m3)\quad{\rm{and}}\quad
u\sum_{l=1}^{ty}r^l+xr^y-x-s\sum_{l=1}^yr^l\equiv0(\mod\frac n2),
\end{eqnarray}
which implies $r^{2y}\equiv r^{y(tv+1)}\equiv1(\mod\frac n2);$

(ii) $\o(\pi(ab))=2$:
$$\begin{array}{lcl}
1&=&(c^{3(vw+yt)}a^{2ur^y\sum_{l=1}^wr^l-2r^yx+2s\sum_{l=1}^yr^l-1})^2\\
 &=&c^{3(vw+yt+tvw+y)}
 a^{2((r^{w+y}+1)(ur^y\sum_{l=1}^wr^l-r^yx+s\sum_{l=1}^yr^l)+
 s\sum_{l=1}^{vw+yt}r^l-1)},
\end{array}$$
that is
{\small{\begin{eqnarray}\label{Df5.3}
\begin{array}{ll}
&vw+yt+tvw+y\equiv0(\mod\frac m3);\\
&(r^{w+y}+1)(ur^y\sum_{l=1}^wr^l-r^yx+s\sum_{l=1}^yr^l)+
 s\sum_{l=1}^{vw+yt}r^l\equiv1(\mod\frac n2),
\end{array}\end{eqnarray}}}
which implies $r^{2w}\equiv1(\mod\frac n2);$

(iii) $\o(\pi(a))=n$:
$$
(bc^{3w})^n=(bc^{3w}bc^{3w})^{\frac n2}=(c^{3w(v+1)}a^{2u\sum_{l=w+1}^{2w}r^l})^{\frac n2}
=c^{3\frac n2w(v+1)}a^{nur^w\sum_{l=1}^{w}r^l}=1,
$$
that is
\begin{eqnarray}\label{Df5.4}
\frac n2w(v+1)\equiv0(\mod\frac m3).
\end{eqnarray}

(iv) $\pi$ preserves $(a^2)^{c_1}=a^{2r}$:
$$\begin{array}{ll}
&((a^2)^{c^3})^c=(c^{3w(v+1)}a^{2u\sum_{l=w+1}^{2w}r^l})^{c^3}=c^{3w(v+1)}a^{2ur\sum_{l=w+1}^{2w}r^l},\\
&(a^{2r})^c=(c^{3w(v+1)}a^{2u\sum_{l=w+1}^{2w}r^l})^r=c^{3wr(v+1)}a^{2ur\sum_{l=w+1}^{2w}r^l},
\end{array}$$
that is
\begin{eqnarray}\label{Df5.5}
w(v+1)(r-1)\equiv0(\mod\frac m3).
\end{eqnarray}

(v) $\pi$ preserves $c_1^a=a^{2s}c_1^{t}$:
$$\begin{array}{ll}
&((c^3)^a)^c=(c^3)^{bc^{3w}}=(a^{2u}c^{3v})^{c^{3w}}=c^{3v}a^{2ur^{w+1}},\\
&(a^{2s}c^{3t})^c=(c^{3w(v+1)}a^{2u\sum_{l=w+1}^{2w}r^l})^sc^{3t}
=c^{3ws(v+1)+3v}a^{2sur^{w+1}\sum_{l=1}^{w}r^l},
\end{array}$$
that is
\begin{eqnarray}\label{Df5.6}
v\equiv ws(v+1)+t(\mod\frac m3)\quad{\rm{and}}\quad
u\equiv su\sum_{l=1}^{w}r^l(\mod\frac n2).
\end{eqnarray}

(vi) $\pi$ preserves $c_1^b=a^{2u}c_1^{v}$:
$$\begin{array}{ll}
&((c^3)^b)^c=(c^3)^{a^{1+2x}b}=(c^3a^{2x(1-r)})^{ab}=c^{3tv}a^{2(u\sum_{l=1}^tr^l-sr+x(r-1))},\\
&(a^{2u}c^{3v})^c=(c^{3w(v+1)}a^{2u\sum_{l=w+1}^{2w}r^l})^uc^{3v}
=c^{3wu(v+1)+3v}a^{2u^2r\sum_{l=w+1}^{2w}r^l},
\end{array}$$
that is
\begin{eqnarray}\label{Df5.7}
\begin{array}{ll}
&t\equiv wu(v+1)+1(\mod\frac m3);\\
&u\sum_{l=1}^tr^l-sr+x(r-1)\equiv u^2r^{w+1}\sum_{l=1}^{w}r^l(\mod\frac n2).
\end{array}
\end{eqnarray}

(vii) $\pi^3=\Inn(c^3)$: Recall $\Inn(c_1)(a)=a^{1-2sr}c_1^{1-t},\,\Inn(c_1)(a^2)=a^{2r}$
and $\Inn(c_1)(b)=c_1^{v-1}a^{2ur}b$.
$$\begin{array}{ll}
&a^{1-2sr}c^{3-3t}=\Inn(c_1)(a)=\pi^3(a)=(a^{1+2x}bc^{3(w+y)})^c\\
&a^{-1}c^{3vt(w+wxv+wx)+3(w+2y)}
a^{2r^w(u\sum_{l=1}^{t(w+wxv+wx)}r^l-s\sum_{l=1}^{w+wxv+wx}r^l-ux\sum_{l=w+1}^{2w}r^l-x)},
\end{array}$$
that is
\begin{eqnarray}\label{Df5.8}
\begin{array}{ll}
&1-t\equiv vt(w+wxv+wx)+w+2y(\mod\frac m3)\\
&1-sr\equiv
r^w(u\sum_{l=1}^{t(w+wxv+wx)}r^l-s\sum_{l=1}^{w+wxv+wx}r^l-ux\sum_{l=w+1}^{2w}r^l-x)(\mod \frac n2);
\end{array}\end{eqnarray}
$$a^{2r}=\Inn(c_1)(a^2)=\pi^3(a^2)=(a^2)^{c^3}
=c^{3w(v+1)+3uw^2(v+1)(uw+1)}a^{2r^w(u\sum_{l=1}^wr^l)^3},
$$
that is
\begin{eqnarray}\label{Df5.10}
w(v+1)+uw^2(v+1)(uw+1)\equiv0(\mod\frac m3)\quad{\rm{and}}\quad
r\equiv r^w(u\sum_{l=1}^wr^l)^3(\mod\frac n2),
\end{eqnarray}
which implies $(u,\frac n2)=(\sum_{l=1}^wr^l,\frac n2)=1$ as $(r,\frac n2)=1$, and

$$c^{3(v-1)}a^{2ur}b=\Inn(c_1)(b)=\pi^3(b)=(c^{3(w+t-1)}a^{2sr-1})^c
=c^{3(t-1)+3wsr(v+1)}a^{2sur\sum_{l=1}^wr^l}b,$$
that is
\begin{eqnarray}\label{Df5.9}
v\equiv t+wsr(v+1)(\mod\frac m3)\quad{\rm{and}}\quad
1\equiv s\sum_{l=1}^wr^l(\mod\frac n2),
\end{eqnarray}
which implies $(s,\frac n2)=1$.

(viii) Insure $\lg c\rg_X=1$:
Since $\lg c\rg_X\leq M$, we get $\lg c\rg_X\leq M_X=\lg a^2\rg\lg c^3\rg$.
Since $\lg c\rg_X=\cap_{x\in X}C^x=\cap_{x\in G}C^x=\cap_{x\in G}\lg c^3\rg^x=\lg c^3\rg_{X_1}=1$ and
$u\sum_{l=1}^{i}r^l\equiv0\equiv s\sum_{l=1}^{i}r^l(\mod\frac n2)\Leftrightarrow
i\equiv0(\mod \frac m3)$, noting  $(s,\frac n2)=(u,\frac n2)=1$, we have
$\sum_{l=1}^{i}r^l\equiv0(\mod\frac n2)\Leftrightarrow i\equiv0(\mod \frac m3)$.
\vskip 3mm
Now we are  ready to determine the parameters by summarizing  Eq(\ref{Df5.1})-Eq(\ref{Df5.9}).

Since $(r,\frac n2)=(u,\frac n2)=1$ (after Eq(\ref{Df5.10})),  we get from   Eq(\ref{Df5.1}) that
$\sum_{l=1}^{v-1}r^l\equiv0(\mod\frac n2).$
By (viii), $\sum_{l=1}^{i}r^l\equiv0(\mod\frac n2)\Leftrightarrow i\equiv0(\mod \frac m3)$,
which means  $v\equiv 1(\mod \frac m3)$.
Inserting $v=1$ in Eq(\ref{Df5.1})-Eq(\ref{Df5.9}), we get
$2w(wu+1)\equiv0(\mod\frac m3)$ and $v\equiv1+2wu(\mod\frac m3)$
in Eq(\ref{Df5.3}), (\ref{Df5.7}) and (\ref{Df5.9}),
which implies $2wu\equiv0(\mod\frac m3)$.
Then $2w\equiv0(\mod\frac m3)$ as $nw\equiv0(\mod\frac m3)$ and $(u,\frac n2)=1$,
which implies $w=\frac m6$ as $w\not\equiv 0(\mod\frac m3)$.
Inserting $w=\frac m6$ in   Eq(\ref{Df5.1})-Eq(\ref{Df5.9}) again,
we get that $r^w\equiv -1(\mod \frac n2)$
in Eq(\ref{Df5.1}) and (\ref{Df5.9}) and $t\equiv1(\mod\frac m3)$ in Eq(\ref{Df5.7}).

Since $2y\equiv0(\mod\frac m3)$ in Eq(\ref{Df5.2}), we know $y$ is either $0$ or $\frac m6$.
If $y=\frac m6=w$, then with Eq(\ref{Df5.2}) and (\ref{Df5.3}), we get $s\equiv u(\mod\frac n2)$.
Then $2sr\equiv 1-r(\mod\frac n2)$ in Eq(\ref{Df5.1}), which implies $(r-1,\frac n2)=1$,
and then $r=1$, contradicting with $r^w\equiv-1(\mod\frac n2)$.
So $y=0$.

By Eq(\ref{Df5.7}), we get $2x\equiv u\sum_{l=1}^{w}r^l+(u\sum_{l=1}^{w}r^l)^2-1(\mod\frac n2),$
which implies $\frac n2$ is odd, then $s\equiv \frac{r^{-1}-1}2(\mod\frac n2)$ in Eq(\ref{Df5.1}).
And by Eq(\ref{Df5.10}), we get $-r\equiv (u\sum_{l=1}^{w}r^l)^3(\mod\frac n2).$
Take $l=-\frac{2ru}{1-r}$, then $r=l^3$, $u=\frac{l(r-1)}{2r}$ and $1+2x\equiv -l+l^2(\mod\frac n2)$.
Since $1+2x$ is odd, we get $1+2x\equiv -l+l^2+\frac n2(\mod n)$.
\qed
\vskip 3mm
In fact, if we add the conditions $t=1$ and $w\ne0$ and delete $\lg c\rg_X=1$ in the above calculation,
then we can get the following:
\begin{lem}\label{DA4.1}
With the notation, suppose   that $t=1$ and $w\ne0$. Then
$$X=\lg a,b,c|R, a^{c^3}=a^r,(c^3)^b=a^{2u}c^{3}, a^c=bc^{\frac {m}2},b^c=a^xb\rg,$$
where $n\equiv 2(\mod 4)$, $m\equiv 0(\mod 6)$, $l^{\frac m2}\equiv-1(\mod\frac n2)$,
$r=l^3$, $u=\frac{l^3-1}{2l^2}$ and $x\equiv -l+l^2+\frac n2(\mod n)$.
\end{lem}

\subsection{$M=\lg a^4\rg\lg c\rg$, $X/M_X\cong S_4$ and $\lg c\rg_X=1$}
\vskip 3mm
\begin{lem}\label{Da4c3}
Suppose  that $X=X(D)$, $M=\lg a^4\rg\lg c\rg,\,X/M_X\cong S_4$ and $\lg c\rg_X=1$. Then
$X=\lg a,b,c|R, (a^2)^{c^3}=a^{2r}, (c^3)^b=a^{\frac{2(l^3-1)}{l^2}}c^3,
(a^2)^c=bc^{\frac {im}2}, b^c=a^{2(-l+l^2+\frac n4)}b,c^a=a^{2+4z}c^{1+\frac{km}3}\rg,$
where either $i=z=0$ and $k=l=1$; or
$i=1$, $n\equiv 4(\mod 8),\, m\equiv 0(\mod 6)$,
  $l^{\frac m2}\equiv-1(\mod\frac n4)$ with $\o(l)=m$, $r=l^3$, $z=\frac {1-3l}{4l}$,
  $k\in\{ 1,2\}$ and $\sum_{i=1}^{j}r^i\equiv0(\mod\frac n2)
  \Leftrightarrow j\equiv0(\mod\frac m3)$.
\end{lem}
\demo
Under the hypothesis, $M_X=\lg a^4\rg\rtimes\lg c^3\rg$. Set $n=\o(a)$ and $m=\o(c)$.
Then $4\di n$ and $3\di m$.
Since $X/M_X=\lg\ola,\olb\rg\lg\olc\rg\cong S_4$, we can choose $\olb$
such that  the form of $X/M_X$ is the following:
$(\ola^2)^{\olc}=\olb,\,\olb^{\olc}=\ola^2\olb$ and $(\olc)^{\ola}=\ola^2\olc^2.$
Take $a_1=a^2$ and $c_1=c^3$.
Then we set $a_1^c=bc_1^w, b^c=a_1^xbc_1^{y},c^a=a_1^{1+2z}c^{2+3d}$,
where $x$ is odd.

Suppose $w\equiv0(\mod\frac m3)$. Note $\o(a^2)=\o((a^2)^c)=\o(b)=2$, then $X\cong S_4$.
Then $l=1,i=z=d=0$ as desired.
So in what follows, we assume $w\not\equiv0(\mod \frac m3)$.

Then consider $X_1=GM_X=\lg a,b\rg\lg c_1\rg$.
Noting $\lg a\rg\lg c_1\rg\leq X_1$ and $\lg c_1\rg_{X_1}=1$, by Lemma~\ref{Da^2},
we know $\lg a_1\rg\lhd X_1$, which implies that $c_1$ normalises $\lg a_1\rg$.
Take $X_2=\lg a_1,b\rg\lg c\rg$.
Then we get $X_2=(\lg a_1,b\rg\lg c_1\rg).\lg c\rg$.
Note that  $c_1$ normalises $\lg a_1\rg$ in $X_2$. Then by Lemma~\ref{DA4.1}, we get
$$X=\lg a_1,b,c|R,
(a_1^2)^{c}=a_1^{2r},c_1^{a_1}=a_1^{2s}c_1,c_1^b=a_1^{2u}c_1,a_1^c=bc^{\frac m2},b^c=a^xb\rg,$$
where
\begin{eqnarray}\label{Df6.1}
\begin{array}{ll}
&n\equiv 4(\mod 8),\, m\equiv 0(\mod 6),\, v=1,\,w=\frac m6,\, y=0, \\
&r=l^3,\,s=\frac{1-l^3}{2l^3},\,u=\frac{l^3-1}{2l^2},\,1+2x=-l+l^2+\frac n4,\,l^{\frac m2}\equiv-1(\mod\frac n4).
\end{array}
\end{eqnarray}

Note $X=X_2.\lg a\rg$.
So $X$ may be defined by $R$ and
\begin{eqnarray}\label{Dm4}
(a_1^2)^{c}=a_1^{2r}, c_1^{a_1}=a_1^{2s}c_1,c_1^b=a_1^{2u}c_1,
a_1^c=bc^{\frac m2},b^c=a_1^xb,c^a=a_1^{1+2z}c^{2+3d}.
\end{eqnarray}
What we should to  determine the parameters  $r,z$ and $d$ by analyse the last one extension.
\vskip 3mm
$X_2.\lg a\rg $, where $c^a=a^{2+4z}c^{2+3d}$.
Set $\pi\in \Aut(X_1) :$
$a^2\to a^2$, $b\to a^{-2}b$ and $c\to a^{2+4z}c^{2+3d}$.
\vskip 3mm
We need to check the following seven equalities:

(i) $\pi$ preserves $(a^2)^{c^3}=a^{2r}$:
$$a^{2r}=((a^2)^{c^3})^a
=(a^2)^{c^{3(2+3d)}}
=a^{2r^{2+3d}},$$
that is
\begin{eqnarray}\label{Df6.2}
r^{1+3d}-1\equiv0(\mod\frac n2).
\end{eqnarray}

(ii) $\pi$ preserves $(c^3)^b=a^{\frac{2(l^3-1)}{l^2}}c^{3}$:
$$\begin{array}{lcl}
((c^3)^b)^a&=&(c^{3(2+3d)}a^{2r(1+2z)+4r^{1+2d}zl+r^d(2lr^{1+d}+2l^2+\frac n2+4l^2z)})^{ba^2}\\
&=&c^{3(2+3d)}a^{2l(l^3-1)+2r(1+2z)+4r^{1+2d}zl+r^d(2lr^{1+d}+2l^2+\frac n2+4l^2z)},
\end{array}$$
that is
\begin{eqnarray}\label{Df6.3}
1-r\equiv
2r(1+2z)+4r^{1+2d}zl+r^d(2lr^{1+d}+2l^2+4l^2z)
(\mod\frac n2).
\end{eqnarray}

(iii) $\pi$ preserves $(a^2)^c=bc^{\frac m2}$:
$$\begin{array}{lcl}
((a^2)^c)^a&=&(a^2)^{c^{2+3d}}=(a^{2(-l+l^2+\frac n4)}bc^{\frac m2})^{c^{3d}}
=a^{2r^d(l^2+\frac n4)-2l}bc^{\frac m2}\\
(bc^{\frac m2})^a&=&ba^2(c^{3(2+3d)}a^{2r(1+2z)+4r^{1+2d}zl+r^d(2lr^{1+d}+2l^2+\frac n2+4l^2z)})^{\frac m6}\\
&=&a^{\frac{2}{1-r}(2r(1+2z)+4r^{1+2d}zl+r^d(2lr^{1+d}+2l^2+\frac n2+4l^2z))-2}bc^{\frac m2},
\end{array}$$
that is,
\begin{eqnarray}\label{Df6.4}
\begin{array}{ll}
&\frac{2}{1-r}(2r(1+2z)+4r^{1+2d}zl+r^d(2lr^{1+d}+2l^2+\frac n2+4l^2z))-2\equiv\\
&2r^d(l^2+\frac n4)-2l(\mod n).
\end{array}
\end{eqnarray}
With Eq(\ref{Df6.3}), we get
$$r^dl\equiv1(\mod\frac n4).$$

(iv) $\pi$ preserves $b^c=a^{2(-l+l^2+\frac n4)}b$:
$$(b^c)^a=(ba^2)^{a^{2+4z}c^2}=a^{2l^2-8zl-8l+\frac n2}b\quad{\rm{and}}\quad
(a^{2(-l+l^2+\frac n4)}b)^a=a^{2(-l+l^2)-2+\frac n2}b,$$
that is
\begin{eqnarray}\label{Df6.6}
z\equiv\frac {1-3l}{4l}(\mod\frac n4).
\end{eqnarray}

(v) $\o(\pi(c))=m$:
$$1=(a^{2+4z}c^{2+3d})^m=c^{m(2+3d)}
 a^{\sum_{i=0}^{\frac m3-1}r^i(2r(1+2z)+4r^{1+2d}zl+r^d(2lr^{1+d}+2l^2+\frac n2+4l^2z))}.$$
Note that $r^{\frac m6}\equiv-1(\mod\frac n4)$.
Then $\sum_{i=0}^{\frac m3-1}r^i\equiv0(\mod\frac n4)$.
Note that $2r+\frac n2+2r^dl^2(l^{2+3d}+1)\equiv0(\mod4)$
which implies $4|(2r(1+2z)+4r^{1+2d}zl+r^d(2lr^{1+d}+2l^2+\frac n2+4l^2z))$.
Then
$$\sum_{i=0}^{\frac m3-1}r^i(2r(1+2z)+4r^{1+2d}zl+r^d(2lr^{1+d}+2l^2+\frac n2+4l^2z))\equiv0(\mod n),$$
as desired.

(vi) $\pi^2=\Inn(a^2)$: Recall $(a^2)^{c^3}=a^{2r}$ and $(a^2)^c=bc^{\frac m2}$,
then we know $\Inn(a^2)(c)=c^{1+\frac m2}a^{-2}b$ and $\Inn(a^2)(c^3)=c^3a^{2-2r}$.
$$
c^{1+\frac m2}ba^2=\Inn(a^2)(c)=\pi^2(c)=(a^{2+4z}c^{2+3d})^a=c^{1+3(1+d)(1+3d)+\frac m2}ba^{2+l(1-r^dl)+\frac n2\sum_{i=0}^{d}r^i},
$$
that is,
\begin{eqnarray}\label{Df6.7}
\begin{array}{ll}
&(1+d)(1+3d)\equiv 0(\mod\frac m3),\\
&0\equiv l(1-r^dl)+\sum_{i=0}^{d}r^i\frac n2(\mod n);
\end{array}
\end{eqnarray}
and
$$c^3a^{2-2r}=\Inn(a^2)(c^3)=\pi^2(c^3)=
c^{3(2+3d)^2}a^{(1+\sum_{i=0}^{1+3d}r^i)(1-2l^2-l^3+2r^{1+d}+\frac n2)},$$
that is,
\begin{eqnarray}\label{Df6.8}
1\equiv(2+3d)^2(\mod\frac m3)\quad{\rm{and}}\quad
0\equiv (\sum_{i=0}^{1+3d}r^i-1)(1-r)(\mod n).
\end{eqnarray}

(vii) Insure $\lg c\rg_X=1$:
Since $\lg c\rg_X\leq M$, we get $\lg c\rg_X\leq M_X=\lg a^4\rg\lg c^3\rg$, which implies
$\lg c\rg_X=\cap_{x\in X}C^x=\cap_{x\in G}C^x=\cap_{x\in G}\lg c^3\rg^x=\lg c^3\rg_{X_1}$.
Then it is suffer to insure $\lg c^3\rg_{X_1}=1$, where
$X_1=\lg a,b,c|R,(a^2)^{c^3}=a^{2l^3}, (c^3)^b=a^{\frac{2(l^3-1)}{l^2}}c^3,
(c^3)^a=c^{3(2+3d)}a^{1-l^3+\frac n2}\rg$.
Then by Lemma~\ref{Dac}, we get
$\sum_{i=1}^{j}r^i\equiv0(\mod\frac n2)\Leftrightarrow j\equiv0(\mod\frac m3)$.
By Eq(\ref{Df6.2}) and $r^{\frac m6}\equiv-1(\mod\frac n4)$, we get
$\sum_{i=1}^{1+3d}r^i\equiv0(\mod\frac n2)$.
Then $1+3d\equiv0(\mod\frac m3)$.
Then $0\equiv l(1-r^dl)(\mod n)$ with Eq(\ref{Df6.7}).
\qed

\subsection{$M=\lg a^3\rg\lg c\rg $ and $X/M_X\cong S_4$}
\begin{lem}\label{Da3c4}
Suppose  that $X=X(D)$, $M=\lg a^3\rg\lg c\rg,\,X/M_X\cong S_4$ and $\lg c\rg_X=1$. Then
\begin{eqnarray}\label{Dm5.0}
X=\lg a,b,c|R,a^{c^4}=a^r, b^{c^4}=a^{1-r}b, (a^3)^{c^{\frac m4}}=a^{-3}, a^{c^{\frac m4}}=bc^{\frac{3m}4}\rg,
\end{eqnarray}
where  $m\equiv 4(\mod8)$ and $r$ is of order $\frac m4$ in $\ZZ_{2n}^*$.
\end{lem}

In this case,  $M_X=\lg a^3\rg\lg c^4\rg$.
Set $a^3=a_1$ and $c^4=c_1$ so that $M_X=\lg a_1\rg \lg c_1\rg $.
Set $\o(a)=n$ and $\o(c)=m$, where $n\equiv0(\mod3)$ and $m\equiv0(\mod4)$.
Then in  Lemma~\ref{Da_1}, we shall show $\lg a_1\rg \lhd X$ and
in Lemma~\ref{Da3c4.1}, we shall get the classification of $X$.
\begin{lem}\label{Da_1}
$\lg a_1\rg\lhd X$.
\end{lem}
\demo  Let $X_1=M_XG$.
 Since $\lg a\rg\lg c_1\rg\leq X_1$ and $\lg c_1\rg_{X_1}=1$, the subgroup  $X_1$  has been given   in Lemma~\ref{Dac}:
 \begin{eqnarray} \label{Dm5.1} X_1=\lg a,b,c_1| \lg a, b\rg, c_1^{\frac m4}=1,
  (a^2)^{c_1}=a^{2r}, (c_1)^a=a^{6s}c_1^{t}, (c_1)^b=a_1^{u}c_1^v\rg ,
 \end{eqnarray}
where
\begin{enumerate}
  \item[\rm(1)] if $t=1$, then $v=1$, $6s\sum_{l=1}^{\frac m4}r^{l}\equiv2(1-6sr)-2r\equiv0(\mod n)$
  and $6s\sum_{l=1}^{w}r^l\equiv0\equiv 3u\sum_{l=1}^{w}(1-6sr)^l(\mod n)$ is if and only if
  $w\equiv0(\mod\frac m4)$;
  \item[\rm(2)] if $t\ne1$, then both $n$ and $u$ are even, and $r,s,t,u,v$ are given by
  $$\begin{array}{ll}
  &r^{t-1}-1\equiv r^{v-1}-1\equiv \frac {3u}2(\sum_{l=0}^{v-1}r^l-1)\equiv0(\mod \frac n2),\\
  &3s\sum_{l=1}^tr^{l}+3sr\equiv 3sr+3s\sum_{l=1}^{v}r^l-\frac {3u}2\sum_{l=1}^tr^l+\frac {3ur}2
  \equiv 1-r(\mod \frac n2),\\
  &t^2-1\equiv v^2-1\equiv 0(\mod\frac m4),\\
  &3s\sum_{l=1}^{i}r^l\equiv\frac {3u}2\sum_{l=1}^{i}r^l\equiv0(\mod\frac n2)\Leftrightarrow
  i\equiv0(\mod\frac m4).
  \end{array}$$
\end{enumerate}
Now $X=\lg X_1, c\rg $ and we need to write the relation of $c$ with $X_1$.
Since $X/M_X=\lg\ola,\olb\rg\rtimes\lg\olc\rg\cong S_4$, checked by Magma, under  our condition,
the only possibilities are the following:
\begin{eqnarray}\label{Dm5.2}
\ola^3=\olc^4=\olb^2=1, \ola^{\olb}=\ola^{-1}, \ola^{\olc}=\ola^i\olb \olc^3 ,
\end{eqnarray}
where $i\in \ZZ_3$.
Observing Eq(\ref{Dm5.1}) and Eq(\ref{Dm5.2}), we may  relabel $a^ib$ by $b$.
Then in $X$, Eq(\ref{Dm5.2}) corresponds to
 \begin{eqnarray}\label{Dm5.3}
 a^3=a_1, b^2=1, c^4=c_1, a^c=bc^{3+4w}.
 \end{eqnarray}
 Moreover, the conjugacy of $\olc$ on $M_X$ is needed:
\begin{eqnarray}\label{Dm5.4}
(a_1)^c=a_1^{z}c_1^{d}.
\end{eqnarray}
Then the group $X$ is uniquely determined by Eq(\ref{Dm5.1}), Eq(\ref{Dm5.3}) and Eq(\ref{Dm5.4}).
For the contrary, we assume $d\ne 0$.
Then we need to deal with two cases, according to the parameter $t$  of $X_1$.
\vskip 3mm
{\it Case 1: $t=1$.}
\vskip 3mm
In this case, $v=1$ and $2(1-6sr)-2r\equiv0(\mod n)$ by $X_1$.
Set $r_1=1-6sr$. Then $r_1\equiv1(\mod3),\,a^{c_1}=a^{r_1}$ and $b^{c_1}=a_1^{ur_1}b$.
By Eq(\ref{Dm5.1}), Eq(\ref{Dm5.3}) and Eq(\ref{Dm5.4}),
one can check $b^c=a_1^{x}bc^{2+4y}$ for some $x$ and $y$.

Since $c$ preserves  $a_1^{c_1}=a_1^{r_1}$, we get
$$((a_1)^{c_1})^c=a_1^{zr_1}c_1^{d}=c_1^{d}a_1^{zr_1^{1+d}}\,{\rm{and}}\,
(a_1^{r_1})^c=(a_1^{z}c_1^{d})^{r_1}=c_1^{dr_1}a_1^{3z\sum_{l=1}^{r_1}r_1^{dl}},$$
which gives
\begin{eqnarray}\label{Df7.1}
d\equiv dr_1(\mod\frac m4).
\end{eqnarray}
Since  $c$ preserves $b^{c_1}=a_1^{ur_1}b$, we get
$$(b^{c_1})^c=a^{{2+3x}r_1+3ur_1}bc^{4y}\quad{\rm{and}}\quad
(a_1^{ur_1}b)^c=c_1^{dur_1}
a_1^{z\sum_{l=1}^{ur_1}r_1^{dl}}a^{2+3x}bc^{4y},$$
which gives
\begin{eqnarray}\label{Df7.2}
du\equiv 0(\mod\frac m4).
\end{eqnarray}
Since $a_1^{c_1}=a_1^{c^4}=a_1^{x_1}c_1^{d(z^3+z^2+z+1)}$ for some $x_1$, we get
\begin{eqnarray}\label{Df7.3}
d(z^3+z^2+z+1)\equiv 0(\mod\frac m4);
\end{eqnarray}
By last equation of Eq(\ref{Dm5.3}), we get  $ac=cbc^{3+4w}$. Then
$$a_1^{ac}=a_1^z c_1^d\,{\rm{and}}\,a_1^{cbc^{3+4w}}=a_1^{x_1}c_1^{d-dz(z^2+z+1)},$$
which gives
\begin{eqnarray}\label{Df7.4}
dz(z^2+z+1)\equiv0(\mod\frac m4).
\end{eqnarray}
With Eq(\ref{Df7.3}), we know $d\equiv0(\mod\frac m4),$
a contradiction.
\vskip 3mm
{\it Case 2: $t\ne 1$.}
\vskip 3mm
Suppose  that $t\ne1$.
Then if $\lg a_1^2\rg\lhd X$, consider $\ox=X/\lg a_1^2\rg$ and $\lg \overline{c_1}\rg\lhd\ox$.
Note $\oc\leq C_{\ox}(\lg \overline{c_1}\rg)=\ox$, then $t=1$, contradicting with $t\ne1$.
So in what follows, we shall show $\lg a_1^2\rg\lhd X$.

By Eq(\ref{Dm5.1}),  both $u$ and $n$ are even and $r\equiv1(\mod3)$.
Then $z$ is odd as $(a,\frac n3)=1$.

By Eq(\ref{Dm5.1}), Eq(\ref{Dm5.3}) and Eq(\ref{Dm5.4}),
one can check $b^c=a^{2+6x}bc^{4y}$ for some $x$ and $y$.

Since $c$ preserves  $(a_1^2)^{c_1}=a_1^{2r}$, we get
$$
((a_1^2)^{c_1})^c=(a_1^zc_1^da_1^zc_1^d)^{c_1}=a_1^{x_1}c_1^{d(t+1)},\,
(a_1^{2r})^c=(a_1^zc_1^d)^{2r}
=a_1^{x_2}c_1^{d(t+1)r},$$
which gives
\begin{eqnarray}
d(t+1)(r-1)\equiv0(\mod\frac m4).
\end{eqnarray}

Since  $c$ preserves $(c_1)^b=a_1^{u}c_1^v$, we get
$$
(c_1^b)^c=c_1^{a^{2+6x}bc^{4y}}=a_1^{x_3}c_1^v,\,
(a_1^{u}c_1^v)^c=(a_1^zc_1^d)^{u}c_1^v=a_1^{x_4}c_1^{\frac{du(t+1)}2+v},
$$
which gives
\begin{eqnarray}\label{Df7.5}
\frac{du(t+1)}2\equiv0(\mod\frac m4).
\end{eqnarray}

Since $c$ preserves  $c_1^a=a_1^{2s}c_1^t$, which implies $c_1^{bc^{2+4w}}=a_1^{2s}c_1^t$, we get
$$c_1^{bc^{2+4w}}=(a_1^{u}c_1^{v})^{c^{2+4w}}=((a_1^zc_1^d)^zc_1^d)^{ur^w}c_1^{v}
=a_1^xc_1^{v},$$
which gives
\begin{eqnarray}
v\equiv t(\mod\frac m4).
\end{eqnarray}

By last equation of Eq(\ref{Dm5.3}) again, we get  $ac^2=cbc^{4(w+1)}$. Then
$$(a_1^2)^{ac^2}=a_1^xc_1^{d(t+1)(z+1)}\,{\rm{and}}\,
(a_1^2)^{cbc_1^{w+1}}=a_1^xc_1^{d(t+1)},$$
which gives
\begin{eqnarray}\label{Df7.6}
dz(t+1)\equiv0(\mod\frac m4).
\end{eqnarray}
Since $\o(a_1^c)=\frac n3$, we get $\frac{dn(t+1)}6\equiv0(\mod\frac m4)$.
With $(\frac n6,z)=1$ and Eq(\ref{Df7.6}), we get $d(t+1)\equiv0(\mod\frac m4)$.
Then $(a_1^2)^c=(a_1^zc_1^d)^2=a_1^xc_1^{d(t+1)}=a_1^x$, which implies $\lg a_1^2\rg\lhd X$,
as desired.
\qed
\begin{lem}\label{DD2}
With the notation, suppose that $\lg a_1\rg\lhd X$, $\lg c\rg_X=1$ and $\lg c\rg$ is $2-$group. Then
$\o(c)=4$.
\end{lem}
\demo
Consider $\ox=X/\lg a_1\rg=\lg\ola,\olb\rg\lg\olc\rg$.
Note $\overline{M_X}=\lg\olc^4\rg$,
Then one can check $C_{\ox}(\lg\olc_1\rg)=\ox$, which implies $\ox$ is the central expansion of $S_4$.
Since the Schur multiplier of $S_4$ is $\ZZ_2$, we know $\o(c)$ is either $4$ or $8$.
For the contrary, we assume $\o(c)=8$.
Then consider $X_1=\lg a,b\rg\rtimes\lg c_1\rg$, where $c_1=c^4$ and $\o(c_1)=2$, again.
By Lemma~\ref{Dac}, we have
$$X_1=\lg a,b,c_1\di R,a^{c_1}=a^r,b^{c_1}=a_1^ub\rg,$$
where $r\equiv1(\mod3).$
Note that $\lg a\rg\leq \lg c\rg_X(\lg a_1\rg)\lhd X$.
Thus $\lg a,bc,c^2\rg\leq \lg c\rg_X(\lg a_1\rg),$ which implies $\lg a_1\rg\times\lg c_1\rg\lhd X$.
One can check $r=1$(as $3r\equiv3(\mod\frac n3)$ and $1\equiv r^2=2r-1(\mod\frac n3)$).
Then $\lg a,c^2\rg\leq \lg c\rg_X(\lg a_1\rg\times\lg c_1\rg)\lhd X$,
which implies $bc\in \lg c\rg_X(\lg a_1\rg\times\lg c_1\rg)$. So $c_1^b=c_1$.
Then $c_1\in \lg c\rg_X=1$, a contradiction.
\qed
\begin{lem}\label{Da3c4.1}
The group $X$ is given by Eq(\ref{Dm5.0}).
\end{lem}
\demo
By lemma, $\lg a_1\rg \lhd X$,  that is $(a_1)^c=a_1^z$ by Eq(\ref{Dm5.4}).
Since $\lg a^2\rg \lhd X_1$, we get $\lg a\rg \lhd X_1$ and so $G\lhd X_1$, that is  $t=v=1$ in Eq(\ref{Dm5.1}).
Then by Eq(\ref{Dm5.1}), (\ref{Dm5.3}) and (\ref{Dm5.4}), we can set
$$X=\lg a,b,c\di R,a^{c_1}=a^{r_1},b^{c_1}=a_1^{ur_1}b,(a_1)^c=a_1^z,a^c=bc^{3+4w}\rg,$$
where $r_1=1-6sr$, $r_1^{\frac m4}-1\equiv2(r_1-r)\equiv0(\mod n)$ and
$r_1^{i}-1\equiv0\equiv 3u\sum_{l=1}^ir_1^l(\mod n)$ is if and only if
$i\equiv0(\mod\frac m4)$.

Set $\lg c\rg=\lg c_2\rg\times\lg c_3\rg$,
where $\lg c_2\rg$ is $2-$group and $\lg c_3\rg$ is $2'-$Hall subgroup of $\lg c\rg$.
Then  $\lg c_1\rg=\lg c_2^4\rg\times\lg c_3\rg$.
And we shall show $c_2^4=1$.

Consider $\ox=X/\lg a_1\rg=\lg\ola,\olb\rg\lg\olc\rg$.
Then one can check $C_{\ox}(\lg\olc_1\rg)=\ox$, which implies
$\lg\overline{c_1}\rg\le Z(\ox)$ and $\ox/\lg \overline{c_1}\rg\cong S_4$.
Note $\lg\olc_3\rg\leq\lg\olc_3\rg(\lg\olc_2^4\rg\lg\ola\rg)\leq\ox$
where $(|\ox:\lg\olc_3\rg(\lg\olc_2^4\rg\lg\ola\rg)|,|\lg\olc_3\rg|)=1$,
then by Proportion~{\ref{complement}}, $\lg\olc_3\rg$ has a complement in $\ox$,
which implies $X=(\lg a,b\rg\lg c_2\rg)\rtimes\lg c_3\rg.$
Consider $X_2=\lg a,b\rg\lg c_2\rg$, where $\lg c_2\rg_{X_2}=1$ and $\lg a_1\rg\lhd X_2$.
Then by Lemma~\ref{DD2}, we get $|\lg c_2\rg|=4$, which implies $\lg c_1\rg=\lg c_3\rg$.
Then
$$X=(\lg a,b\rg\lg c_2\rg)\rtimes\lg c_1\rg.$$

In what follows, we shall determine $X$.
In $X_2=\lg a,b\rg\lg c_2\rg$, we know $a^{c_2}=bc_2^3$ by Eq(\ref{Dm5.3}).
Consider $\lg a\rg\leq C_{X_2}(\lg a_1\rg)\lhd X_2$.
Then $C_{X_2}(\lg a_1\rg)$ is either $\lg a,bc_2,c_2^2\rg$ or $X_2$.
Suppose   that $\lg c\rg_X(\lg a_1\rg)=X$. Then we know $a_1^2=1$ as $[a_1,b]=1$, that is $n=3$ or $6$.
Then one can check
$$\begin{array}{ll}
&X=X_2\cong S_4,\,{\rm{if}}\, n=3;\\
&X=X_2=\lg a,b,c|a^6=b^2=c^4=1,a^b=a^{-1},a^c=bc^3,b^c=a^2b\rg,\,{\rm{if}}\,n=6.
\end{array}$$

Suppose   that $\lg c\rg_X(\lg a_1\rg)=\lg a,bc_2,c_2^2\rg$.
Then $a_1=a_1^{bc_2}=(a_1^{-1})^{c_2}$, which implies $z=-1$.
In $X_1=\lg a,b\rg\rtimes\lg c_1\rg$, we know
$$X_1=\lg a,b,c_1\di R,a^{c_1}=a^{r_1},b^{c_1}=a_1^{ur_1}b\rg.$$
Since $c_2$ preserves  $a^{c_1}=a^{r_1}$, we get
$$(bc_2^3)^{c_1}=a_1^{ur_1}bc^3\,{\rm{and}}\,(a^{r_1})^{c_2}=a_1^{\frac{1-r_1}3}bc^3,$$
which gives
$$ur_1\equiv \frac{1-r_1}3(\mod\frac n3).$$
Then
$$X=\lg a,b,c\di R,a^{c_1}=a^{r_1},b^{c_1}=a_1^{\frac{1-r_1}3}b,(a_1)^{c_2}=a_1^{-1},a^{c_2}=bc_2^3\rg,$$
where $c_2=c^{\frac m4}$ and $\o(r_1)=\frac m4$.
\qed

\end{document}